\documentclass{amsart}
\usepackage{leftindex}
\usepackage{amsmath,amssymb,amscd,amsthm,amsfonts}
\usepackage{upgreek}
\usepackage{stmaryrd}

\usepackage[hyphens]{url}
\usepackage{hyperref}

\usepackage{enumerate}
\usepackage{mathtools}
\usepackage{makecell}

\newtheorem{thm}{Theorem}[section]

\newtheorem{lemma}[thm]{Lemma}
\newtheorem{prop}[thm]{Proposition}
\newtheorem{cor}[thm]{Corollary} 
\theoremstyle{definition}
\newtheorem{defn}[thm]{Definition} 
\theoremstyle{remark}
\newtheorem{eg}[thm]{Example}
\newtheorem{rmk}[thm]{Remark}
\numberwithin{equation}{section}

\newcommand{\spacedvert}{\,\vert\,}

\newcommand{\spe}{sp}
\newcommand{\spr}{sp}

\newcommand{\CC}{\mathbb{C}}

\newcommand{\GG}{\mathbb{G}}

\newcommand{\RR}{\mathbb{R}}
 
\newcommand{\ZZ}{\mathbb{Z}}

\newcommand{\bmu}{\boldsymbol\upmu}
\newcommand{\bpsi}{\boldsymbol\uppsi}

\newcommand{\bm}{\mathbf{m}}

\newcommand{\bp}{\mathbf{p}}

\newcommand{\bT}{\mathbf{T}}

\newcommand{\cH}{\mathcal{H}}

\newcommand{\cO}{\mathcal{O}}
\newcommand{\cP}{\mathcal{P}}

\newcommand{\cR}{\mathcal{R}}
\newcommand{\cS}{\mathcal{S}}

\newcommand{\cZ}{\mathcal{Z}}

\newcommand{\tA}{\widetilde{A}}

\newcommand{\tE}{\widetilde{E}}

\newcommand{\tG}{\widetilde{G}}
\newcommand{\tH}{\widetilde{H}}

\newcommand{\tK}{\widetilde{K}}

\newcommand{\tM}{\widetilde{M}}
\newcommand{\tN}{\widetilde{N}}

\newcommand{\tP}{\widetilde{P}}

\newcommand{\tV}{\widetilde{V}}

\newcommand{\ta}{\widetilde{a}}

\newcommand{\tg}{\widetilde{g}}
\newcommand{\tilh}{\widetilde{h}}

\newcommand{\tm}{\widetilde{m}}

\newcommand{\tilr}{\widetilde{r}}

\newcommand{\tw}{\widetilde{w}}
\newcommand{\tx}{\widetilde{x}}

\newcommand{\talpha}{\widetilde{\alpha}}

\newcommand{\half}{\frac{1}{2}}
\newcommand{\triv}{\mathbf{1}}

\DeclareMathOperator{\Rat}{\mathrm{Rat}}
\DeclareMathOperator{\Stab}{\mathrm{Stab}}

\DeclareMathOperator{\GL}{\mathrm{GL}}

\DeclareMathOperator{\SL}{\mathrm{SL}}

\DeclareMathOperator{\rU}{\mathrm{U}}
\DeclareMathOperator{\Sp}{\mathrm{Sp}}

\DeclareMathOperator{\Mp}{\mathrm{Mp}}

\DeclareMathOperator{\SO}{\mathrm{SO}}

\DeclareMathOperator{\rO}{\mathrm{O}}

\DeclareMathOperator{\Ind}{\mathrm{Ind}}
\DeclareMathOperator{\ind}{\mathrm{ind}}
\DeclareMathOperator{\Hom}{\mathrm{Hom}}
\DeclareMathOperator{\End}{\mathrm{End}}

\DeclareMathOperator{\Inn}{\mathrm{Inn}}

\DeclareMathOperator{\Rep}{\mathrm{Rep}}

\newcommand{\Id}{\mathrm{Id}}

\newcommand{\isom}{\simeq}

\DeclareMathOperator{\rank}{\mathrm{rank}}

\DeclareMathOperator{\Fr}{\mathrm{Fr}}

\DeclareMathOperator{\supp}{\mathrm{supp}}

\newcommand{\lmod}{\backslash}

\let\Re\undefined
\DeclareMathOperator{\Re}{\mathrm{Re}}
\let\Im\undefined
\DeclareMathOperator{\Im}{\mathrm{Im}}

\newcommand{\nr}{\mathrm{nr}} %

\newcommand{\lsup}[2]{\ensuremath\leftindex^{{#1}}_{} {{#2}}}

\newcommand{\form}[2]{\langle{#1},{#2}\rangle}

\def\smatc[#1,#2,#3,#4]{\ensuremath\bigl( \begin{smallmatrix}
#1&#2\\ #3&#4
\end{smallmatrix} \bigr)}

\def\bmatc[#1,#2,#3,#4]{\ensuremath
  \begin{pmatrix}
    #1&#2\\ #3&#4
  \end{pmatrix}}

\def\vecc[#1,#2]{\ensuremath \begin{pmatrix}
#1\\#2
\end{pmatrix}}

\newcommand{\tEBO}{\tE_{B_{\cO}}}
\newcommand{\tEKBO}{\tE_{K(B_{\cO})}}
\newcommand{\BO}{B_{\cO}}
\newcommand{\KBO}{K(B_{\cO})}

\def\forma[#1,#2]{\ensuremath\langle #1 , #2 \rangle}

\newcommand{\Rao}{\mathrm{Rao}}

\newcommand{\Xnr}{\mathfrak{X}^{\mathrm{\nr}}}
\newcommand{\red}{\mathrm{red}}

\newcommand{\Jord}{\mathrm{Jord}}
\newcommand{\alt}{\mathrm{alt}}
\newcommand{\MOD}{\mathcal{M}od}

\newcommand{\scrS}{\mathscr{S}} %

\usepackage[shortlabels]{enumitem}
\usepackage{sidenotes}
\usepackage[bb=stix,scr=rsfs]{mathalpha}
\usepackage{tikz-cd}

\newcommand{\Fz}{{F_{0}}} %
\newcommand{\Fs}{F} %
\begin{document}

\title[Intertwining Algebras and Affine Hecke Algebras]{Intertwining Algebras and Affine Hecke Algebras for Finite Central Extensions of Classical $p$-adic Groups with Application to Metaplectic Groups}
\author{Volker Heiermann}
\address{Aix-Marseille Universit\'e, CNRS, I2M, UMR 7373, 13453 Marseille, France}
\email{volker.heiermann@univ-amu.fr}

\author{Chenyan Wu}
\address{School of Mathematics and Statistics,
  the University of Melbourne,
  Victoria, 3010, Australia}
\email{chenyan.wu@unimelb.edu.au}

\begin{abstract}
  For a finite central extension $\tG$ of a classical $p$-adic reductive group, we consider the endomorphism algebra of some induced
  projective generator \`a la Bernstein of the category of smooth
  representations of $\tG$. In the case where the Levi subgroups decompose, we compute this algebra to get a result similar to the one previously obtained by the first author for classical $p$-adic groups, showing that this intertwining algebra is a twisted semi-direct product of an affine Hecke algebra with parameters by a
 twisted finite group algebra. We discuss also the  case with non-decomposed Levi subgroups.
  
We  then give an application to the category of genuine representations of a $p$-adic metaplectic group. Using results of C. M\oe glin relative to the Howe correspondence, we show that the Bernstein components of these groups are equivalent to tensor products of categories of unipotent representations of classical groups. This generalizes a previous result of the first author. It implies an equivalence of categories between the category of genuine representations of the $p$-adic metaplectic group and the direct sums of those of smooth representations of the corresponding split odd special orthogonal group and its pure inner form.
\end{abstract}

\maketitle{}

\section{Introduction}
\label{sec:introduction}

Fix a  non-archimedean local field $\Fz$.
Let $\GG$ be a reductive group over $\Fz$ and set $G = \GG(\Fz)$.
Let $\bm$ be a positive integer, $\bmu_{\bm}$ the group of $\bm$-th roots of unity and $\tG$ a central extension of $G$ by $\bmu_{\bm}$.

We will consider the category $\Rep(\tG)$ of smooth complex representations of $\tG$
and specialize to the case where $\GG $ is  a classical group,  meaning a general linear group, a symplectic group, a special or full orthogonal group or  a unitary group.

Let $M$ be a Levi subgroup of $G$.
By abuse of language, we also say that its preimage $\tM$ in $\tG$ is a Levi subgroup of $\tG$.
Let $(\widetilde{\sigma},\tE)$ be a supercuspidal representation of $\tM$. We will denote its equivalence class by the same symbol, giving an additional explanation only if it is not clear from the context what is meant.

Denote by $\Xnr(\tM)$ the group of unramified characters of $\tM$.
It can be identified with $\Xnr(M)$ and acts by torsion on $\Rep(\tM)$.
Let $\cO$ be the orbit of the equivalence class of $\widetilde{\sigma}$ by $\Xnr(\tM)$.
Let $W$ be the Weyl group of $G$ with respect to the maximal split torus $A_{0}$ in the fixed minimal Levi subgroup $M_{0}$ of $G$.
Write $\lsup{W}{\cO}$ for the union of $w\cO$, $w\in W$.
Let $P$ be a parabolic subgroup with $M$ as its Levi subgroup.
Define the category $\Rep_{\tG}(\cO)$ to be the full subcategory of $\Rep(\tG)$ whose objects are the representations $\pi$ such that all the Jordan--H\"older factors of the Jacquet modules $r^{\tG}_{\tP'}\pi$ are in $\lsup{W}{\cO}$ if $P'$ is associate to $P$ and are non-supercuspidal, if $P'$ is not associate to $P$.

By the theory of the Bernstein center \cite{Bernstein-MR771671}, which holds  for finite central extensions, $\Rep(\tG)$ decomposes into a direct sum of subcategories $\Rep_{\tG}(\cO)$ parameterized by equivalence classes of orbits $\cO $ as above.

Let $\tE_{1}$ be an irreducible sub-representation of $\widetilde{\sigma}|_{\tM^{1}}$, where $\tM^{1}$ is the intersection of the kernels of all unramified characters of $\tM$.
The induced representation $i_{\tP}^{\tG} \ind_{\tM^{1}}^{\tM}\tE_{1}$ is then a projective generator of the category $\Rep_{\tG}(\cO)$.
For reductive groups, this is a result of J.N. Bernstein whose proof given in \cite{Roche-Bernstein-decomp-MR1934305} generalizes to covering groups.
An immediate corollary of this is that the category $\Rep_{\tG}(\cO)$ is equivalent to the category of right $\End_{\tG}(i_{\tP}^{\tG} \ind_{\tM^{1}}^{\tM}\tE_{1})$-modules.

When  $\tG$ is a reductive classical group, it was shown in \cite{Heiermann-MR2827179} and \cite{Heiermann-MR2975415} that the latter  is the semi-direct product of a Hecke algebra with parameters by a finite group algebra and that the resulting equivalence of categories preserves parabolic induction, Jacquet functor, tempered and discrete series representations.

However, when $\tG$ is no longer a reductive group, the situation becomes more complicated. We will show that things go still quite well if the Levi subgroups of $\tG$ decompose and give an application in the case of the metaplectic group, where we obtain complete results analogous to \cite{Heiermann-MR3719526}.

If the Levi subgroups do not decompose, important ingredients for the approach in \cite{Heiermann-MR2827179} may fail, including multiplicity 1. Also the contribution of the $R$-group becomes more complicated to deal with. We do think that this can still be managed: the approach in \cite{Heiermann-MR2827179}   is general, as is already mentioned there.
Using results of A. Roche \cite{Roche-MR2508719}, M. Solleveld \cite{Solleveld-MR4432237} made this  more precise under the multiplicity $1$ assumption, adding some different perspective when multiplicity $1$ fails.

\section*{Acknowledgement}
\label{sec:acknowledgement}

We would like to thank the Institut de Math\'ematiques de Marseille (I2M) of Aix-Marseille University and CNRS for hosting the second author for three months. The second author had financial support from the School of Mathematics and Statistics of the University of Melbourne and she was able to be on leave thanks to the Special Studies Program of the University of Melbourne. We  thank also the Institut Henri Poincar\'e, Paris, where both authors stayed for two weeks in 2025 during the trimester program ``Representation Theory and $K$-theory''. The financial support for the stay of the first author was provided jointly by the organizers of the trimester and the Institut de Math\'ematiques de Marseille.
We are grateful to Wen-Wei Li and Fan Gao for the invaluable discussions on cover groups with decomposed Levis and again to Wen-Wei Li
for sending us his work which shows that most of the Brylinski--Deligne covers of classical groups have decomposed Levi subgroups.
Wen-Wei Li also  pointed us to the beautiful work of his student, Jiahe Chen, on the extended Langlands--Deligne parameters of the metaplectic groups.
The second author would like to thank Wee Teck Gan and Rui Chen for some enlightening discussions.

\setcounter{tocdepth}{1} 
\tableofcontents

\section{Preliminaries}
\label{sec:preliminaries}

Let $\Fz$ be a non-archimedean local field, $\cO_{F_{0}}$ its ring of integers, $\varpi_{\Fz}$ a uniformiser of $\Fz$ and  $q_\Fz$ the cardinality of the residue field.
Denote by $|\cdot |_{\Fz}$ the norm of $\Fz$ such that $|\varpi_{\Fz} |_{\Fz} = q_\Fz^{-1}$.

Let $G$ be the group of $\Fz$-points of a connected reductive group $\GG$ over $\Fz$, $\bm$ a positive integer and $\tG$ a central extension of $G$ by $\bmu_{\bm}$, so that we have a short exact sequence
\begin{align*}
  1 \rightarrow \bmu_{\bm} \rightarrow \tG \rightarrow G \rightarrow 1.
\end{align*}
Write $\bp: \tG \rightarrow G$ for the projection.
We may also consider finite central extensions $\tG$ of $G$ by a finite abelian group $A$, but we can reduce the study of their representations to the case of finite cyclic central extensions as follows.
If  $(\pi,V)$ is a representation of $\tG$, then we can decompose $V$ according to the action of $A$, $  V = \oplus_{\eta} V_{\eta}$ where $\eta$ runs over characters of $A$ and then each $V_{\eta}$ is a representation of $\tG$ and the action factors through $\tG / \ker(\eta)$, which is a finite cyclic central extension of $G$.
From now on, we focus on  central extensions of $G$ by $\bmu_{\bm}$.
Some of the fundamental results for reductive groups generalise to the central extensions and we refer the readers to \cite{Kaplan-Szpruch-MR4613279} for details.

To alleviate notation, when there is no confusion, we will abuse language and use $G$ when we mean $\GG$.
Later, we will specialize to the case where $G$ is   a symplectic, a full orthogonal, a special orthogonal, a unitary or a general linear group.
As the full even orthogonal group is not connected, remarks on how to treat them will be supplied later (see Section~\ref{sec:decom-levis}).

For a subgroup $H$ of $G$, we will use $\tH$ to denote the preimage of $H$ under the projection $\tG\rightarrow G$.

Fix a minimal parabolic subgroup $P_{0}$ of $G$.
Let $U_{0}$ be the unipotent radical of $P_{0}$.
Fix  a Levi subgroup $M_{0}$ of $P_{0}$ and  a maximal split torus $A_{0}$ of $M_{0}$.
We have $P_{0}=M_{0}U_{0}$.
Let $W := W^{G} := N_{G}(A_{0}) / Z_{G}(A_{0})$ be the Weyl group of $G$ defined relative to $A_{0}$.
It is equipped with the length function $\ell$ determined by $P_{0}$ (See \eqref{eq:length-Weyl}).
Let $K$ be a maximal compact subgroup of $G$ which is in good position relative to $A_{0}$ by Bruhat--Tits theory.

By \cite[Proposition (a) in Appendix 1]{MW-eis-book-MR1361168}, whose proof applies to the setting of  local fields, there exists a unique continuous $P_0'$-equivariant section
\begin{equation}\label{eq:section-unip}
  \iota_{U_0'}: U_0' \rightarrow \tG
\end{equation}
of the projection $\tG\rightarrow G$ where $P_{0}'$ is any parabolic subgroup with Levi subgroup $M_{0}$ and $U_{0}'$ is its unipotent radical;
thus we may regard $U_0'$ as a unipotent subgroup of $\tG$ via this section.
Remark that for the unipotent radicals $U'$ and $U''$ of two parabolic subgroups with Levi subgroup $M_0$, $\iota_{U'}$ agrees with $\iota_{U''}$ on $U'\cap U''$.
This can be seen as follows.
Consider the map of $U'\cap U''$ given by $\eta: u \mapsto \iota_{U'}(u)\iota_{U''}(u)^{-1}$.
This is a character of $U'\cap U''$ and locally constant, by continuity.
Thus there exists an open neighbourhood $C$ of $1$ in $U'\cap U''$ on which $\eta$ is trivial.
For any element $u$ in $U'\cap U''$, there exists  $a \in A_{0}$ such that $aua^{-1}\in C$.
Then $\iota_{U'}(aua^{-1}) = \iota_{U''}(aua^{-1})$ and by the equivariance of conjugation by parabolic subgroups, we get $\ta \iota_{U'}(u) \ta^{-1} = \ta \iota_{U''}(u) \ta^{-1}$
where $\ta$ is any lift of $a$ in $\tG$.
Thus we get $\iota_{U'}(u) = \iota_{U''}(u)$ for all $u\in U'\cap U''$.

In the sequel we will consider a fixed parabolic subgroup $P\supseteq P_0$ of $G$, $P=MU$, where $M$ is the Levi factor which contains $M_0$.
Write $\cP(M)$ for the set of parabolic subgroups of $G$ with Levi factor $M$ and write $\cP(\tM)$ for the preimages in $\tG$.

Let $A_{M}$ be the maximal split torus contained in the centre of $M$.
Let $A_{M}^{\dagger} := A_{M}^{\bm} = \{ a^{\bm} \spacedvert a\in A_{M} \}$.
It is an open and closed subgroup of $A_{M}$ of finite index and $\tA_{M}^{\dagger}$ is central in $\tM$ (\cite[Proposition~2.1.1 and its proof]{WWLi-trace-formula-II-MR3053009}).

Write $\Rat(M)$ for the group of rational characters of $M$.
Set $a_{M}^{*} = \Rat(M) \otimes_{\ZZ} \RR$ and  $a_{M} = \Hom_{\RR}(a_{M}^{*},\RR)$.
There is a natural pairing between $a_{M}^{*}$ and $a_{M}$.
For two Levi subgroups $M\subseteq M'$ of $G$, let $a_{M}^{M'}$ be the subspace of $a_{M}$ that is orthogonal to $a_{M'}^{*}$.

We write $\Sigma(A_{M})$ for the $\Fz$-roots of $A_{M}$ in $G$ and $\Sigma(A_{M},P)$ for the subset of roots in the unipotent radical $U$ of $P$.

We will put in the subscript `$\red$' to indicate that we take only the reduced roots.
In addition, when there is confusion, we will use a superscript to indicate the ambient group.
There is a bijection $\alpha \mapsto M_{\alpha}$ between $\Sigma_{\red}(A_{M},P)$ and the set of Levi subgroups of $G$ for which $M$  is a maximal proper Levi subgroup.

To $\lambda\in a_{M,\CC}^{*}:=\Rat(M) \otimes_{\ZZ} \CC$, we associate a character $\chi_{\lambda}$ of $M$ determined as follows, called an \it unramified character \rm of $M$.
For $\lambda = \alpha\otimes s$,  we set $\chi_{\lambda}(m) = |\alpha(m)|_{\Fz}^{s}$. Denote by $\Xnr(M)$ the group formed by these unramified characters of $M$
and by $M^{1}$  the intersection of  $\ker(\vert\chi \vert_{\Fz})$, $\chi \in\Rat(M)$. The quotient $M/M^1$ is a lattice of dimension equal to the rank of $M$. An element of $\Xnr(M)$ is nothing but a smooth complex character of $M$ that is trivial on $M^1$.

We let an unramified character $\widetilde{\chi }$ of $\tM $  be an unramified character $\chi $ of $M$ composed with the projection $\bp$, $\widetilde{\chi }=\chi\circ \bp$, and denote by $\Xnr(\tM)$ the group of unramified characters of $\tM $, by $\tM^{1}$ the preimage of $M^{1}$ in $\tG$. Thus, an unramified character of $\tM $ is nothing but a smooth complex character of $\tM $ which is trivial on $\tM^{1}$. In addition, $\tM/\tM^{1}$ is isomorphic to $M/M^1$ via the projection $\bp$.

The map $\Rat(M) \otimes_{\ZZ} \CC \rightarrow \Xnr(\tM)$ is surjective.
For a topological group $H$, we write $X(H)$ for $\Hom(H,\CC^{\times})$ where $\Hom$ means continuous group homomorphisms.

Let $B = B_{\tM}=\CC[\tM/\tM^1]$ be the ring of regular functions on the complex affine variety $\Xnr(\tM)$.
It is  integral and we write $K(B)$ for its fraction field.

Let $\cO$ be the set of equivalence classes of an $\Xnr(\tM)$-orbit of an irreducible supercuspidal representation of $\tM$.

For $w\in W^{G}$, we use $\tw$ to denote a lift of $w$ in $\tK$, unless otherwise specified.
We note that the conjugate action of any element $\tx \in\tG$ on $\tG$ depends only on its image in $G$.
A Weyl element $w\in W^{G}$ acts on the set $\cP(\tM)$ of preimages of parabolic subgroups of $G$ of which $M$ is a Levi factor
by conjugation by its lift $\tw$  in $\tK$ and we use the notation $w\tP'$ for $\tw \tP' \tw^{-1}$ for $\tP'\in\cP(\tM)$.
It is clear that $w\tP$ is the preimage of $wP$ in $\tG$.
We use a similar notation for the Weyl group-conjugates of Levi subgroups.

Let $W(M)$ be the set of representatives in $W^{G}$ of
\begin{equation*}
  \{ w\in W^{G} \spacedvert w^{-1}M = M \} / W^{M}
\end{equation*}
which are of minimal length in their class mod $W^{M}$.
This is a subgroup of $W^{G}$ and  can be defined alternatively as the subgroup of $W^{G}$  that consists of elements that stabilise $M$ and $M\cap P_{0}$.
For $\widetilde{\sigma}\in\cO$,  the equivalence class of representations $w\widetilde{\sigma}$ is well-defined for $w\in W(M)$.
Let
\begin{equation*}
  W(M,\cO) = \{ w\in W(M) \spacedvert w\cO = \cO \}.
\end{equation*}

For a parabolic subgroup $P'$ containing $M_{0}$, write $\overline{P'}$ for the opposite parabolic of $P'$.
The length function for $w\in W$ is given by
\begin{align}\label{eq:length-Weyl}
  \ell(w) = | \Sigma_{\red}(A_{0},P_{0}) \cap \Sigma_{\red}(A_{0},w \overline{P_{0}}) | .
\end{align}
For $w\in W(M)$, define
\begin{align*}
  \ell_{M}(w) = | \Sigma_{\red}(A_{M},P) \cap \Sigma_{\red}(A_{M},w \overline{P}) |.
\end{align*}

\begin{prop} [{\cite[Proposition~1.1]{Heiermann-MR2827179}}]
  For two elements $w$ and $w'$ of $W(M)$ satisfying $\ell_{M}(ww') =  \ell_{M}(w') + \ell_{M} (w)$, it is necessary and sufficient that $\ell(ww') = \ell(w) + \ell(w')$.
\end{prop}

\section{Standard Intertwining Operators and Harish-Chandra's $\mu$-function}
\label{sec:intertw-oper}

Let $P=MU, P'=MU'$ be two parabolic subgroups with Levi factor $M$,  $(\widetilde{\sigma},\tE)$  a  representation of finite length of $\tM$ and $\cO$  the $\Xnr(\tM)$-orbit of $\widetilde{\sigma}$ in $\Rep(\tM)$.
(Here $\widetilde{\sigma}$ is not supposed to be supercuspidal irreducible, but later we will do so.)
We follow  \cite[Sec.~2.4]{WWLi-trace-formula-II-MR3053009} to define the intertwining operator
\begin{align*}
  J_{\tP' | \tP} (\widetilde{\sigma}) : i_{\tP}^{\tG} \tE \rightarrow i_{\tP'}^{\tG} \tE.
\end{align*}
For $f\in i_{\tP}^{\tG} \tE$, set
\begin{align*}
  (J_{\tP' | \tP} (\widetilde{\sigma}) f) (\tx) = \int_{(U\cap U')\lmod U'} f(u' \tx ) du'.
\end{align*}
We may identify the space $i_{\tP}^{\tG} \tE$ with $\ind_{\tP\cap\tK}^{\tK} \tE$ via restriction to $\tK$.
For $\widetilde{\sigma}$ varying in $\cO$, we use the same representation space $\tE$ and thus, as long as the integral is absolutely convergent, we get a family of operators $J_{\tP' | \tP} (\widetilde{\sigma})$ which send $\ind_{\tP\cap\tK}^{\tK} \tE$ to $\ind_{\tP'\cap\tK}^{\tK} \tE$.
\cite[Th\'eor\`eme~2.4.1]{WWLi-trace-formula-II-MR3053009} shows that
there exists $R\in\RR$ such that if $\form{\Re(\chi)}{\alpha^{\vee}} > R$ for all $\alpha\in\Sigma(A_{M},P)\cap\Sigma(A_{M},\overline{P'})$, then $J_{\tP' | \tP} (\widetilde{\sigma}\otimes\chi)$ is defined by a convergent integral.
Furthermore, if $\widetilde{\sigma}$ is tempered, then $R$ can be taken to be $0$.
The operator $J_{\tP' | \tP}$  extends to a rational operator on all $\cO$, i.e. its matrix coefficients can be identified with elements in $K(B)$.

Let $w\in W^{G}$ and $\tw\in\tK$ a representative of $w$.
Define the operator
\begin{align*}
\lambda(\tw):  i_{\tP}^{\tG}(\tE) &\rightarrow i_{w\tP}^{\tG}(\tw \tE) \\
  f(\cdot)          &\mapsto f(\tw^{-1}\cdot).
\end{align*}
Here $\tw \tE$ is the space of representation  of $\tw\widetilde{\sigma}$ which has the same underlying vector space as $\tE$.

Set $d(P',P):= |\Sigma_{\red}(A_{M}, P') \cap \Sigma_{\red}(A_{M},\overline{P})|$. The usual composition and compatibility relations for the intertwining operators still hold for covering groups:

\begin{prop}[{\cite[Proposition~2.4.2]{WWLi-trace-formula-II-MR3053009}}]\label{prop:property-JPP}
 Let $(\widetilde{\sigma},\tE)$ be an admissible representation of $\tM$ of finite length.
 Let $P,P'$ be two parabolic subgroups of $G$ with Levi factor $M$.
  Then
  \begin{enumerate}
  \item $J_{\tP'|\tP}(\widetilde{\sigma})^{\vee} = J_{\tP|\tP'}(\widetilde{\sigma}^{\vee})$ where ${}^{\vee}$ means taking contragredient;
  \item $J_{\tP'|\tP}(\widetilde{\sigma}) = J_{\tP'|\tP''}(\widetilde{\sigma}) J_{\tP''|\tP}(\widetilde{\sigma})$, if $d(P',P) = d(P',P'') + d(P'',P)$;
  \item For $P''=M''U''$ containing $P$ and $P'$, with the identifications
    \begin{align*}
      i_{\tP}^{\tG}(\widetilde{\sigma}) =& i_{\tP''}^{\tG} i_{\tP\cap\tM''}^{\tM''} (\widetilde{\sigma})\\
      i_{\tP'}^{\tG}(\widetilde{\sigma}) =& i_{\tP''}^{\tG} i_{\tP'\cap\tM''}^{\tM''} (\widetilde{\sigma}),
    \end{align*}
 the operator $J_{\tP'|\tP}^{\tG}(\widetilde{\sigma})$ is  deduced from $J_{\tP'\cap\tM''|\tP\cap\tM''}^{\tM''}(\widetilde{\sigma})$ via the induction functor $i_{\tP''}^{\tG}(\cdot)$;
  \item Let $w\in W^{G}$ and $\tw\in\tK$ a representative of $w$.
    Then
    \begin{align*}
      J_{w\tP'|w\tP}(\tw \widetilde{\sigma}) = \lambda(\tw) J_{\tP'|\tP}(\widetilde{\sigma}) \lambda(\tw)^{-1}.
    \end{align*}
  \end{enumerate}
\end{prop}

Assume now that $\widetilde{\sigma}$ is square-integrable irreducible.
Then there exists a rational function $j$ on $\cO$ independent of the choice of $P\in\cP(M)$  such that 
\begin{align*}
  J_{\tP|\widetilde{\overline{P}}}(\widetilde{\sigma}')   J_{\widetilde{\overline{P}}|\tP}(\widetilde{\sigma}') = j(\widetilde{\sigma}'), \quad\text{for all $\widetilde{\sigma}'\in\cO$}.
\end{align*}
We have
\begin{align*}
  j^{\tG} = \prod_{\alpha\in\Sigma_{\red}(A_{M})/\pm} j^{\tM_{\alpha}}
\end{align*}
where we use superscripts to indicate the ambient group.
See, for example, paragraphs after \cite[Proposition~2.4.2]{WWLi-trace-formula-II-MR3053009}.

The Harish-Chandra $\mu$-function is defined by
\begin{align*}
  \mu := j^{-1}.
\end{align*}
In \cite[V.2.]{Waldspurger-plancherel-MR1989693}, there are some other factors in the definition of $\mu$.

\begin{prop}[{\cite[Proposition 2.4.3]{WWLi-trace-formula-II-MR3053009}}]\label{prop:mu-property}
  The Harish-Chandra $\mu$-function is a rational function on $\cO$.
  It is regular and non-negative real on $\Im(\cO)$.
  We have
  \begin{align*}
    \mu^{\tG} = \prod_{\alpha\in\Sigma_{\red}(A_{M})/\pm} \mu^{\tM_{\alpha}}.
  \end{align*}
  It is invariant under the action of $W^{G}$ and changing $\widetilde{\sigma}$ to $\widetilde{\sigma}^{\vee}$.
\end{prop}

We go back to our original setting, which means that  $\cO$ is assumed to contain a supercuspidal representation.

\begin{thm}[{\cite[Th\'eor\`eme~1.1, Lemme~1.2]{Heiermann-notes}}]\label{thm:mu-zero-irreducibility}
  Let $\alpha\in \Sigma(A_{M},P)$ and let $\widetilde{\sigma}\in\cO$ be a unitary supercuspidal irreducible representation of $\tM$.
  Write $W^{M_{\alpha}}$ for the Weyl group  of $M_{\alpha}$.
  Let $s=s_{\alpha}$ be the unique element of $W^{M_{\alpha}}$  of minimal length in the coset $W^{M_{\alpha}}/W^{M}$ such that $s(\overline{P}\cap M_{\alpha })\supseteq P_0\cap M_{\alpha }$.
  Then the following hold.
  \begin{enumerate}
  \item If $sM \neq M$ or $s \widetilde{\sigma} \not\isom \widetilde{\sigma}$, then $i_{\tP\cap\tM_{\alpha}}^{\tM_{\alpha}} (\widetilde{\sigma}\otimes\chi_{x \alpha})$ is irreducible for all $x\in \RR$ and
    the map $x \mapsto \mu^{\tM_{\alpha}}(\widetilde{\sigma}\otimes\chi_{x \alpha})$ is regular and nonzero for $x\in\RR$.
  \item Assume that $sM = M$ and $s \widetilde{\sigma} \isom \widetilde{\sigma}$.   Then the following hold.
    \begin{enumerate}
    \item There exists a unique $x_{0}\ge 0$ such that $i_{\tP\cap\tM_{\alpha}}^{\tM_{\alpha}} (\widetilde{\sigma}\otimes\chi_{x_{0} \alpha})$ is reducible and then $i_{\tP\cap\tM_{\alpha}}^{\tM_{\alpha}} (\widetilde{\sigma}\otimes\chi_{-x_{0} \alpha})$ is also reducible.
    \item The representation $i_{\tP\cap\tM_{\alpha}}^{\tM_{\alpha}} \widetilde{\sigma}$ is irreducible if and only if $\mu^{\tM_{\alpha}}(\widetilde{\sigma})=0$. Furthermore, if $\mu^{\tM_{\alpha}}(\widetilde{\sigma})\neq 0$, then $i_{\tP\cap\tM_{\alpha}}^{\tM_{\alpha}} \widetilde{\sigma}$ is a direct sum of two non-isomorphic irreducible representations.
    \item Let $x_{0}>0$. Then $(\mu^{\tM_{\alpha}})^{-1}$ is defined at $\widetilde{\sigma}\otimes\chi_{x_{0} \alpha}$. The representation  $i_{\tP\cap\tM_{\alpha}}^{\tM_{\alpha}} (\widetilde{\sigma}\otimes\chi_{x_{0} \alpha})$ is reducible if and only if $x_{0}$ is a pole of the function $x\mapsto\mu(\widetilde{\sigma}_{x\alpha})$.
      This pole is simple.
    \item Up to a twist by a character in $\Xnr(\tM_{\alpha})$, there exists at most one $\widetilde{\sigma}_{-} \in \cO$ of the form $\widetilde{\sigma}\otimes\chi_{x\alpha}$ for some $x\in\CC$ such that $\widetilde{\sigma}_{-}\not\isom\widetilde{\sigma}$ and $s\widetilde{\sigma}_{-} \isom \widetilde{\sigma}_{-}$.
    \item If the intertwining operator $J_{\widetilde{\overline{P}} \cap \tM_{\alpha}  | \tP\cap\tM_{\alpha}}^{\tM_{\alpha}}$ relative to $\tM_{\alpha}$ is regular on all $\widetilde{\sigma}\otimes\chi_{x \alpha}$ for $x\in\RR$, then it is bijective for all $x\in\RR$. %
    \end{enumerate}
    \item The poles of the intertwining operator $J_{\widetilde{\overline{P}} \cap \tM_{\alpha}  | \tP\cap\tM_{\alpha}}^{\tM_{\alpha}}$ relative to $\tM_{\alpha}$ are exactly the zeroes of  $\mu ^{\tM_{\alpha}}$. They are all simple.
    \item If $J_{\widetilde{\overline{P}} \cap \tM_{\alpha}  | \tP\cap\tM_{\alpha}}^{\tM_{\alpha}}$ is regular at $\widetilde{\sigma}$, then it is regular and bijective for all $\widetilde{\sigma}\otimes\chi_{x\alpha}$, $x\in \RR$.
  \end{enumerate}
\end{thm}

\begin{prop}[{c.f.  \cite[Proposition~1.3]{Heiermann-MR2827179}}] \label{prop:is-root-sys} 
    The set 
\begin{equation*}
  \Sigma_{\cO,\mu} = \{ \alpha \in \Sigma_{\red}(A_{M}) \spacedvert  \mu^{\tM_{\alpha}} \text{ has a zero on } \cO \}
\end{equation*}
is a root system.
  For $\alpha\in \Sigma_{\cO,\mu}$, denote by $s_{\alpha}$ the unique element of $W^{M_{\alpha}}(M,\cO)$ such that $s_{\alpha}(P\cap M_{\alpha}) = \overline{P}\cap M_{\alpha}$.
  Then the Weyl group $W_{\cO}$ of $\Sigma_{\cO,\mu}$ identifies with the subgroup of $W^{G}(M,\cO)$ generated by the reflections $s_{\alpha}$ for $\alpha\in \Sigma_{\cO,\mu}$.
  For each $\alpha\in \Sigma_{\cO,\mu}$, denote by $\alpha^{\vee}$ the unique element of $a_{M}^{M_{\alpha}}$ which satisfies $\form{\alpha}{\alpha^{\vee}}=2$.
  Then the set $\Sigma_{\cO,\mu}^{\vee} = \{ \alpha^{\vee} \spacedvert \alpha\in\Sigma_{\cO,\mu} \}$ is the set of coroots of $\Sigma_{\cO,\mu}$, the duality being that between $a_{M}$ and $a_{M}^{*}$, and
  the set $\Sigma(P) \cap \Sigma_{\cO,\mu}$  is the set of positive roots for a certain order on $\Sigma_{\cO,\mu}$.
\end{prop}
\begin{defn}
  Set
  \begin{equation*}
    \Sigma_{\cO,\mu}(P) = \Sigma(P) \cap \Sigma_{\cO,\mu}
  \end{equation*}
  and let $\Delta_{\cO,\mu}$ be the set of simple roots.
  Let $W_{\cO}$ be the Weyl group of $\Sigma_{\cO,\mu}$ and $\ell_{\cO}$ be the length function of $W_{\cO}$.
\end{defn}

\begin{rmk}\label{rmk:rep-inertial-orbit}
  As elements in $\Delta_{\cO,\mu}$ are linearly independent and $\mu^{\tM_{\alpha}}$ is invariant under the translation action of $\Xnr(\tM_{\alpha})$,  we can fix a unitary supercuspidal representation $\widetilde{\sigma}$ of $\tM$ whose equivalence class belongs to $\cO$ such that  $\mu^{\tM_{\alpha}} (\widetilde{\sigma}) = 0$ for every $\alpha\in \Delta_{\cO,\mu}$.
  
\end{rmk}

Recall the Harish-Chandra map \cite[top of p. 240]{Waldspurger-plancherel-MR1989693}
\begin{align*}
  H_{M}:M\rightarrow a_{M}
\end{align*}
determined by 
\begin{align*}
  q_{\Fz}^{- \form{H_{M}(m)}{\alpha}} = |\alpha(m)|_{\Fz},\quad\text{for all  $\alpha\in\Rat(M)$.}
\end{align*}
We define $H_{\tM}$ to be the composition of $H_{M}$ and the projection $\tM \rightarrow M$.

We have an embedding of $\Xnr(\tM_{\alpha})$ to $\Xnr(\tM)$ given by restriction.
We note that $\Xnr(\tM)$ is identified with $\Xnr(M)$.
Following \cite{Heiermann-MR2827179}, for $\alpha\in \Sigma_{\cO,\mu}$, there exists $h_{\alpha} \in M \cap M_{\alpha}^{1}$ such that $H_M(h_{\alpha})$ is a positive multiple of $\alpha^{\vee}$ and that for every $\chi\in\Xnr(M)$, $\chi(h_{\alpha})=1$ is equivalent to $\chi\in \Xnr(M_{\alpha})$.
Write $\tilh_{\alpha}$ for any lift in $\tM$ of $h_{\alpha}$.
Then the last condition can be formulated as follows:
for  $\chi\in\Xnr(\tM)$, $\chi(\tilh_{\alpha})=1$ is equivalent to $\chi\in \Xnr(\tM_{\alpha})$.
We note that $H_{M}(h_{\alpha})$ is uniquely determined, since $M \cap M_{\alpha}^{1} / M^{1}$ is a free $\ZZ$-module of rank $1$.

Set $\Stab(\cO) = \{\chi\in\Xnr(\tM)  \spacedvert \widetilde{\sigma}\otimes\chi \isom \widetilde{\sigma} \}$.
This is a finite subgroup of $\Xnr(\tM)$ that depends only on $\cO$.
Let $t_{\alpha}$ be the smallest positive integer such that $\chi(\tilh_{\alpha}^{t_{\alpha}}) = 1$ is equivalent to $\chi\in \Xnr(\tM_{\alpha})\Stab(\cO)$.

Let $b_{h_{\alpha}}$ be the map $\Xnr(\tM) \rightarrow \CC$, $\chi \mapsto \chi(\tilh_{\alpha})$.
This is an element in $B_{\tM} = \CC[\tM/\tM^1]$, the ring of polynomials of the complex affine variety $\Xnr(\tM)$.
Put $Y_{\alpha} = b_{h_{\alpha}}$ and $X_{\alpha} = Y_{\alpha}^{t_{\alpha}}$.
The following lemma  involves only the underlying reductive group  $G$.
\begin{lemma}[{\cite[Lemma~1.5]{Heiermann-MR2827179}}]
  For $w\in W(M)$, we have $\lsup{w}{Y_{\alpha}} = Y_{w\alpha}$.
  In particular, $\lsup{s_{\alpha}}{Y_{\alpha}} = Y_{\alpha}^{-1}$.
\end{lemma}

\begin{prop}[{\cite[Corollaire~1.4]{Heiermann-notes}}]\label{prop:expression-mu-function}
  Let $\alpha \in \Sigma_{\red}(P)$ and $s=s_{\alpha}$.
  Then there exist a  real constant $c_{s}' > 0$ and two real numbers $a_{s} \ge 0$ and $a_{s,-} \ge 0$ such that
  \begin{align*}
    \mu^{\tM_{\alpha}}(\widetilde{\sigma}\otimes\cdot) = c_{s}' \frac{(1-X_{\alpha}(\cdot)) (1-X_{\alpha}^{-1}(\cdot))}{(1-X_{\alpha}(\cdot)q_{\Fz}^{-a_{s}}) (1-X_{\alpha}^{-1}(\cdot) q_{\Fz}^{-a_{s}})}
    \frac{(1+X_{\alpha}(\cdot)) (1+X_{\alpha}^{-1}(\cdot))}{(1+X_{\alpha}(\cdot)q_{\Fz}^{-a_{s,-}}) (1+X_{\alpha}^{-1}(\cdot) q_{\Fz}^{-a_{s,-}})} .
  \end{align*}
  In particular, $\mu^{\tM_{\alpha}}$ is non-constant, if and only if at least one of $a_{s}$ or $a_{s,-}$ is nonzero.

\end{prop}

\begin{rmk}
  We may shift the base point $\widetilde{\sigma}$ such that for every $\alpha\in\Delta_{\cO,\mu}$, the real numbers $a_{s_{\alpha}}$ and $a_{s_{\alpha},-}$ in Proposition~\ref{prop:expression-mu-function}  satisfy $a_{s_{\alpha}}\ge a_{s_{\alpha},-}$.
  We will now assume that this is the case.
  This implies that $a_{s_{\alpha}} > 0$ for $\alpha\in\Delta_{\cO,\mu}$.
\end{rmk}

\begin{defn}
   Put
  \begin{equation*}
    R(\cO) = \{ w\in W(M,\cO) \spacedvert w\alpha \in \Sigma(P) \text{ for every } \alpha\in \Sigma(P) \cap \Sigma_{\cO,\mu}\}.
  \end{equation*}
\end{defn}

\begin{prop}\label{prop:WMO=RO-WO}
  The group  $R(\cO)$ is a subgroup of $W(M,\cO)$.
  We have
  \begin{align*}
    W(M,\cO) = R(\cO) \ltimes W_{\cO}.
  \end{align*}
\end{prop}
\begin{proof}
  The proof of \cite[Proposition~1.12]{Heiermann-MR2827179} carries over.
\end{proof}

\section{The case of decomposed Levi subgroups}
\label{sec:decom-levis}

We will from now on assume that $\GG$ is a symplectic, full or special orthogonal,  unitary or general linear group.
(These are exactly the general linear group and the quasisplit classical groups and their pure inner forms.)
If $\GG $ is a unitary group, $\Fs$ will denote the splitting field of the quasi-split inner form of $\GG$.
Otherwise, we put $\Fs=\Fz$ to get uniform statements.

Recall that $G=\GG(\Fz)$.
When we say that $G$ is of type $A$, $B$, $C$ or $D$, we mean by that that the relative root system of $G$ is of that type.
In particular, it is of type $C$ if and only if $\GG$ is a symplectic group or a quasi-split even unitary group and of type $D$ if and only if $\GG$ is a split full or special even orthogonal group.
Else, $G$ is of type $B$ except if $\GG$ is a general linear group.

Recall that $P$, $P_0$ and $M$, $M_0$ have been fixed.
If $\GG $ is an orthogonal group (which is not connected), then the parabolic subgroups $P$ and $P_0$ are assumed in addition to be cuspidal in the sense of \cite{Goldberg-Herb-MR1422314} and we add to the Weyl group of its connected component the element which induces a non-trivial outer isomorphism of order $2$ of the underlying special orthogonal group.

It follows from \cite[Lemma~4.7]{Ban-Jantzen-MR2017065} that the analog of Harish-Chandra's result remains valid for this Weyl group element, allowing only reducibility at $0$.

We have the following proposition which serves to set up some notation for the reduced roots and Weyl group elements.
\begin{prop}\label{prop:str-Weyl-group}
  \begin{enumerate}[leftmargin=*]
  \item After conjugating $\tM$, we can assume $M=\bp(\tM)$ to be of the form
    \begin{align*}
      \begin{array}{c}
       \GL_{d_{1}} (\Fs)\times\cdots \times\GL_{d_{1}} (\Fs)\times\cdots \times\GL_{d_{r}} (\Fs)\times\cdots \times\GL_{d_{r}}(\Fs)\times H_{d},
      \end{array}
    \end{align*}
    where $\GL_{d_{i}} (\Fs)$ occurs $k_{i}$ times\footnote{We count $\GL_{d_{i}} (\Fs)$ and $\GL_{d_{j}} (\Fs)$ differently, if $i\ne j$, even if $d_i=d_j$.}
for $i=1,\ldots, r$ and  $H_{d}$ is a reductive group of  absolute semisimple rank $d$ of the same type as $G$.\footnote{For $d=0$ the group $H_{d}$ does not appear except if $G$ is an odd unitary or the full odd orthogonal group.}
Suppose that $G$ is of type $D$.
Then, we have $d\ne 1$ and, for $d=0$, we may assume that the last simple root of $G$ is not a root of  $M$.\footnote{This may require an outer automorphism or even an isomorphism with another cover of $G$ inducing an outer automorphism of $G$.}

\item Identify $A_{M}$ with $\GL_1(\Fs)^{k_{1}} \times \GL_1(\Fs)^{k_{2}} \times \cdots \times \GL_1(\Fs)^{k_{r}}$.
For $i=1,\ldots , r$, $i'=1,\dots , r-1$ and $j = 1,\ldots, k_{i}$, set $\alpha^\mp_{i,j}$ (resp. ${\alpha_{i',k_{i'}}'} ^{\mp}$)  to be the rational character of $A_{M}$ which sends
       \begin{align*}
      x=(x_{1,1}, x_{1,2}, \ldots , x_{1,k_{1}},
      x_{2,1}, x_{2,2}, \ldots , x_{2,k_{2}},
      \ldots , x_{r,1}, x_{r,2}, \ldots , x_{r,k_{r}})  \in A_{M}
    \end{align*}
    to $x_{i,j} x_{i,j+1}^{\mp 1}$ if $j < k_{i}$ and to $x_{i,k_{i}}$ if $j=k_{i}$\footnote{Remark that $\alpha^+_{i,k_i}=\alpha^-_{i,k_i}$, but it is handy for uniform statements.}
(resp. to $x_{i',k_{i'}}x_{i'+1,1}^{\mp 1}$). We will frequently write $\alpha_{i,j}$ for $\alpha^-_{i,j}$ and $\alpha_{i',k_{i'}}'$ for ${\alpha '_{i',k_{i'}}}^-$.

    The $\alpha_{i,j}^\mp$, $j\ne k_i$, ${\alpha_{i',k_{i'}}'}^{\mp}$ are roots of $A_M$ in the unipotent radical $U$ of $P$ and, if $d\ne 0$, the $\alpha_{i,k_i}$ are always too.
    When these are roots, denote by $s_{i,j}^\mp$ (resp. ${s_{i',k_{i'}}'} ^{\mp}$) the element $s_{\alpha_{i,j}^\mp}$  (resp. $s_{{\alpha_{i',k_{i'}}'} ^{\mp}}$) of  $W^{M_{\alpha_{i,j}^\mp}}$ (resp. $W^{M_{{\alpha_{i',k_{i'}}'} ^\mp}}$) considered in Theorem~\ref{thm:mu-zero-irreducibility}. 
    
    The element $s_{i,k_i}$ corresponds, if defined, to the outer automorphisms of $\GL_{d_{i}} (\Fs)$. We let $s_{i,k_i}$ stand also for the action by the outer automorphism of $\GL_{d_{i}} (\Fs)$, if $\alpha_{i,k_i}$ is not a root.

    One has $s_{i,j}^{\mp}M=M$ and $s_{i,k_i}M=M$, except if $G$ is an even special orthogonal group, $d_i$ is odd, $d=0$, when the latter does not hold.
    
    In addition, ${s_{i,k_i}'} ^{\mp}M=M$ if and only if $d_i=d_{i+1}$.

    \item Suppose $j\ne k_i$. The action of $s_{i,j}$ just permutes the $j$th and the $j+1$th copy of $\GL_{d_{i}}(\Fs)$ in $\bp(\tM)$. The action of $s_{i,j}^{+}$ permutes also the $j$th and the $j+1$th copy of $\GL_{d_{i}} (\Fs)$  in $\bp(\tM)$, but it applies in addition on each of the two copies the outer automorphism of $\GL_{d_{i}} (\Fs)$ preserving the set of positive roots.

   \item Suppose now $d\ne 0$ or that $G$ is of type $B$.

   If $G$ is either not a special even orthogonal group or $d_i$ is even, then $s_{i,k_i}$ acts only on the last copy of  $\GL_{d_{i}} (\Fs)$ in $M=\bp(\tM)$ and this action is the same outer automorphism of $\GL_{d_{i}} (\Fs)$ as above.

   If $G$ is a special even orthogonal group and $d_i$ is odd,  then $s_{i,k_i}$ acts on the last copy of  $\GL_{d_{i}} (\Fs)$ and on $H_{d}$. The action on $\GL_{d_{i}} (\Fs)$ is again the outer automorphism above of $\GL_{d_{i}} (\Fs)$ and the action on $H_{d}$ is given by the outer automorphism of $H_{d}$ that sends positive roots to positive roots.

    \item Suppose now $d=0$ and $G$ is not of type $B$.

    If $G$ is of type $C$, then $\alpha_{i,k_i}$ is not a root, but $2\alpha_{i,k_i}$ is and $s_{2\alpha_{i,k_i}}$ acts as $s_{i,k_i}$ defined above.

    If $G$ is an even special or full orthogonal group and $d_i$ is even, then $2\alpha_{i,k_i}$ is a root and $s_{2\alpha_{i,k_i}}$ acts again as $s_{i,k_i}$ defined above. If $d_i$ is odd, then $s_{i,k_i}$ is not in the Weyl-group if $G$ is the special orthogonal group, but it is when $G$ is the full orthogonal group. 
     \item Suppose $d_i\ne d_j$ if $i\ne j$.
     \begin{enumerate}
     \item If either $G$ is of type $B$, or $d\ne 0$,  then $W(M)$ is generated by the  $s_{\alpha_{i,j}}$ and $s_{\alpha _{i,k_i}}$, $i=1,\dots, r$; $j=1,\dots , k_i-1$.
     \item If $d=0$ and $G$ is of type $C$ or if $G$ is a full even orthogonal group, then $W(M)$ is generated by the  $s_{i,j}$ and $s_{2\alpha _{i,k_i}}=s_{i,k_i}$, $i=1,\dots, r$, $j=1,\dots , k_i-1$.
     \item If $d=0$ and $G$ is an even special orthogonal group, the group  $W(M)$ is generated by the  $s_{i,j}$, $j=1,\dots , k_i-1$, the $s_{i,k_i}$, $i=1,\dots, r$, with $d_i$ even, the $s_{i,k_i-1}^+$ with $d_i$ odd, and by products of even numbers of $s_{i,k_i}$ with $d_i$ odd.
     \end{enumerate}
     \end{enumerate}
\end{prop}

\null
\begin{proof}
The form of the Levi subgroup is standard knowledge for $G$. As the projection $\bp:\tM\rightarrow M$ commutes with conjugation in $\tG$, one can always assume $\tM$ to have this form. In addition, conjugation in $\tG$ leaves $\tM$ invariant, if and only if it leaves $M$ invariant.

The rational characters $\alpha_{i,j}^{\mp }$, $j\ne k_i$ are clearly restrictions of roots relative of $G$ to $A_M$ and they are positive. If $d\ne 0$, then there is always a root of $G$ that restricts to $\alpha_{i,k_i}^{\mp }$, as $H_d$ is semisimple.

If $G$ is an even special orthogonal group, $d=0$ and $d_i$ odd, then the outer automorphism of $G$ is a factor of $s_{i,k_i}$, and it follows that the corresponding factor $\GL_{d_i}(\Fs)$ of $M$ is not preserved by $s_{i,k_i}$ and $s_{i,k_i}\tM\ne\tM$.

Let us now go to the properties in (3) and (4).

If $d\ne 0$ or $G$ is of type $B$, then the properties (3) are standard.
For (4), one only needs to distinguish the different cases for the even special orthogonal group, which are caused by its outer automorphism, which is an inner automorphism of the full orthogonal group.
The constraint for the special orthogonal group is that Weyl group elements can act only by even number of sign changes.
If $d_i$ is even, then the Weyl group element does not involve the outer automorphism of $H_d$.
However, if $d_i$ is odd, it does.

Suppose $d=0$ and $G$ not of type $B$.
The case $G$ of type $C$ is as above after replacing $\alpha _{i,k_i}$ by $2\alpha _{i,k_i}$.

When $G$ is of type $D$ and $d_i$ even or $G$ is the full orthogonal group, then by restriction we always find that $2\alpha _{i,k_i}$ is a root and its symmetry acts like $s_{i,k_i}$.
But, in case that $d_i$ is odd and $G$ an even special orthogonal group, the corresponding symmetry denoted here $s_{i,k_i}$ does not lie in the Weyl group of $G$ or preserve $M$ as already observed above.
However, products of even numbers of $s_{i,k_i}$ for such $i$ lie in the Weyl group of $G$ and preserve $M$.

For (6), it follows from what we just did that all these Weyl group elements lie in $W(M)$.
Let now be $w\in W(M)$.
Then $w$ must send each $\GL_{d_i}(\Fs)$ to a copy of itself applying possibly in addition an outer automorphism.
We will call length of $w$ the length of the permutation adding the number of outer automorphism involved and will do an induction on this length of $w$.
If the length of $w$ is $0$, then clearly $w$ is the identity.
To apply the induction, it will be sufficient to show that for an arbitrary $w\in W(M)$ we can find a symmetry as above (or possibly a product of symmetries) so that the product lowers the length of $w$.
If $w$ involves permutations, then the above clearly provides a symmetry which is a permutation lowering the length of $w$.
If $w$ involves only outer automorphism, then the same is true when $G$ is of type $B$, $C$, a full even orthogonal group or one of the outer automorphisms is on $\GL_{d_i}(\Fs)$ with $d_i$ even.
If $G$ is a special even orthogonal group and the outer automorphisms are all on $\GL_{d_i}(\Fs)$ with $d_i$ odd, then, as the corresponding outer automorphisms of $\GL_{d_i}(\Fs)$ involves an odd number of sign changes and cannot be realized by an element in the Weyl group of $G$, $w$ must apply at least to another $\GL_{d_j}(\Fs)$, $d_j$ odd, an outer automorphism.
But this can be realized by a product of a conjugate $s_{i,k_i}$ with a conjugate of $s_{j,k_j}$, lowering the length of $w$ as defined above.
\end{proof}

\begin{defn}\label{defn:prod-of-covers}
Let $\widetilde{M}$  be a standard Levi subgroup of $\widetilde{G}$ and write $$\bp(\widetilde{M})=\GL_{d_1,1}(\Fs)\times\cdots\times \GL_{d_r,r}(\Fs)\times H_d$$ with $H_d$ a classical group and $\GL_{d_i,i}(\Fs)$ isomorphic to $\GL_{d_i}(\Fs)$, the second index serving only to distinguish different copies of $\GL_{d_i}(\Fs)$.
Denote by $\widetilde{\GL_{d_i,i}(\Fs)}$ the inverse image in $\widetilde{G}$ of the subgroup of $\bp(\widetilde{M})$ with $i$-th factor equal to $\GL_{d_i,i}(\Fs)$ and all other factors equal to $\{1\}$.
Define similarly $\tH_d$.

We will say that  \it the standard Levi subgroup $\widetilde{M}$ of $\widetilde{G}$ decomposes \rm if the map
$$\widetilde{\GL_{d_1,1}(\Fs)} \times\cdots\times \widetilde{\GL_{d_r,r}(\Fs)}\times \tH_d\rightarrow \widetilde{M}$$
$$(\widetilde{g_1},\dots,\widetilde{g_r},\widetilde{h})\rightarrow \widetilde{g_1}\cdots\widetilde{g_r}\widetilde{h}$$
is a group homomorphism.

We will say that the \it Levi subgroups of $\widetilde{G}$ decompose \rm if all standard Levi subgroups of $\widetilde{G}$ decompose.

The definitions still make sense for the cover of a general linear group with $H_d$ omitted.
\end{defn}

\begin{rmk}
  By \cite[Cor.~4.6]{li2025decomposedlevisubgroupsbdcovers}, the standard Levi subgroups of the Brylinski--Deligne covers of many classical groups, namely those of symplectic groups,  special odd orthogonal groups or special even orthogonal groups with dimension of the defining space $\ge 6$, all decompose.
  
  For a cover $\widetilde{\rO_{n}}$ of a  full orthogonal group $\rO_{n}$, consider the decomposability of Levi subgroups $\tM$ as follows.
  Write $\widetilde{\SO_{n}}$ for the pull-back of $\widetilde{\rO_{n}}$  along $\SO_{n} \rightarrow \rO_{n}$.
  Assume first that $M \isom \GL_{d_{1}} (\Fs)\times\cdots \times\GL_{d_{r}} (\Fs)$.
  Then $M$ is also a Levi subgroup of $\SO_{n}$ so that the decomposability of $\tM$ is answered by the case of special orthogonal groups.
  More precisely, if  $\widetilde{\SO_{n}}$ is a Brylinski--Deligne cover and if $n\ge 6$ when $n$ is even and put no condition on $n$ when $n$ is odd, then $\tM$ decomposes.
  Assume next that $M \isom \GL_{d_{1}} (\Fs)\times\cdots \times\GL_{d_{r}}(\Fs) \times \rO_{d}$ where $d\ge 1$.
  Then   \cite[Prop.~7.1]{li2025decomposedlevisubgroupsbdcovers} reduces the decomposability of $\tM$
  to
  \begin{enumerate}
  \item the decomposability of $\tM \cap \widetilde{\SO_{n}}$,
  \item the commutation relations of  preimages in $\widetilde{\rO_{n}}$ of any element $\tau$ in $\rO_{d}\setminus \SO_{d}$ and preimages in $\widetilde{\rO_{n}}$ of the elements $t$ in the maximal torus of  the $\GL$-factors of $M$.
  \end{enumerate}
\end{rmk}

\begin{lemma}[{see also \cite[Section~2.6]{WWLi-trace-formula-I-MR3176600}}]\label{lemma:conj} 
Let $\tH$ be a subgroup of $\tG$ and denote by $\tN$ a subgroup of its normalizer.

Assume that $\bp(\tN)$ commutes with $\bp(\tH)$.

Then, for all $\tg\in\tN$, there is a  character $\chi_{\tg}$ of $\tH$  with values in $\bmu_\bm$ inflated from a character of $H$ such that, for all $\tilh\in\tH$, $\tg\tilh\tg^{-1}=\chi_{\tg}(\tilh)\tilh$. The map $\tg\mapsto \chi_{\tg}$ from $\tN$ to the group of $\bmu_\bm$-valued smooth characters of $\tH$ is a smooth group homomorphism.
\end{lemma}

\begin{proof}
Let $\tg$ be in $\tN$.
As $\bp(\tN)$ commutes with $\bp(\tH)$, $\bp(\tg^{-1}\tilh\tg)=\bp(\tilh)$ for all $\tilh\in\tH$.
So, there exists $\xi_{\tilh}\in\bmu_\bm$ such that $\tg\tilh\tg^{-1}=\xi_{\tilh}\tilh$.
As conjugation is a group homomorphism, for $\tilh, \tilh'\in\tH$, $\xi_{\tilh\tilh'}=\xi_{\tilh}\xi_{\tilh'}$, which proves that $\tilh\mapsto \xi_{\tilh}$ is a $\bmu_\bm$-valued character of $\tH$ which is smooth as conjugation by $\tg$ is smooth.

Write $\Inn(\tg)$ for the above automorphism induced by $\tg\in\tN$ on $\tH$ by conjugation.
The map $\tg\mapsto \Inn(\tg)$ is a homomorphism from $\tN$ to the group of automorphism of $\tH$, which shows that $\tg\mapsto \chi_{\tg}$ is a group homomorphism.
It is smooth because $\tG $ is a topological group.
\end{proof}

\begin{defn}\label{defn:mu-prod-of-covers}
For $\bmu_{\bm}$-covers of a general linear group and a classical group $H$,  we define
$$\widetilde{\GL_{d_{1},1}(\Fs)} \times_{\bmu_{\bm}} \cdots\times_{\bmu_{\bm}} \widetilde{\GL_{d_{r},r}(\Fs)} \times_{\bmu_{\bm}} \tH$$
to be the quotient of $\widetilde{\GL_{d_{1},1}(\Fs)}\times\cdots\times  \widetilde{\GL_{d_{r},r}(\Fs)}\times \tH$
by the subgroup
\begin{equation*}
 C = \{ (\zeta_{1},\dots, \zeta_{r},\zeta_{0}) \in \bmu_{\bm} \times\dots\times \bmu_{\bm} \times \bmu_{\bm} \spacedvert \zeta_{1}\cdots \zeta_{r}\zeta_{0} = 1 \}.
\end{equation*}
\end{defn}

\begin{rmk}\label{rmk:prod-cover-repns}
  If $\widetilde{\rho}_{1},\dots, \widetilde{\rho}_{r}$ are  smooth representations of $\widetilde{\GL_{d_{1},1}(\Fs)} ,\dots, \widetilde{\GL_{d_{r},r}(\Fs)}$ and $ \widetilde{\tau} $ is a  representation of   $\tH$ such that $\widetilde{\rho}_{1}|_{\bmu_{\bm}}, \ldots , \widetilde{\rho}_{r}|_{\bmu_{\bm}}$ and $\widetilde{\tau}|_{\bmu_{\bm}}$ act by the same character of $\bmu_{\bm}$, then the representation
\begin{align*}
        \widetilde{\rho}_{1} \otimes\cdots\otimes  \widetilde{\rho}_{r} \otimes \widetilde{\tau}
\end{align*}
of $\widetilde{\GL_{d_{1},1}(\Fs)} \times \cdots\times\widetilde{\GL_{d_{r},r}(\Fs)}\times\tH$
factors through a representation of $\widetilde{\GL_{d_{1},1}(\Fs)} \times_{\bmu_{\bm}} \cdots\times_{\bmu_{\bm}} \widetilde{\GL_{d_{r},r}(\Fs)} \times_{\bmu_{\bm}} \tH$ that we will denote by
\begin{align*}
 \widetilde{\rho}_{1} \otimes_{\bmu_{\bm}}\cdots\otimes_{\bmu_{\bm}}  \widetilde{\rho}_{r} \otimes_{\bmu_{\bm}} \widetilde{\tau} .
\end{align*}

 All irreducible smooth representations of $\widetilde{\GL_{d_{1},1}(\Fs)}\times_{\bmu_{\bm}} \cdots\times_{\bmu_{\bm}} \widetilde{\GL_{d_{r},r}(\Fs)}\times_{\bmu_{\bm}} \tH$ have this form.
\end{rmk}

\begin{prop}\label{prop:decomposed-Levis} 
Let $\tG$ be a cover of a symplectic, full or special orthogonal, or unitary group $G$. Let $\tM$ be a standard Levi subgroup of $\tG$ that decomposes, $$\bp(\widetilde{M})=\GL_{d_{1},1}(\Fs) \times \cdots\times\GL_{d_{r},r}(\Fs)\times H_d$$ with $d\ne 1$ if $G$ is of type $D$. The associated homomorphism 
\begin{align*}
\widetilde{\GL_{d_{1},1}(\Fs)}\times\cdots\times\widetilde{\GL_{d_{r},r}(\Fs)}\times\tH_d\rightarrow\tM
\end{align*}
induces an isomorphism 
$$\widetilde{\GL_{d_{1},1}(\Fs)}\times_{\bmu_{\bm}} \cdots\times_{\bmu_{\bm}}\widetilde{\GL_{d_{r},r}(\Fs)}\\ \times_{\bmu_{\bm}}\tH_d\rightarrow\tM.$$

Note that here $k_i=1$ for all $i$ in the notations of \ref{prop:str-Weyl-group}
. Denote simply by $s_i'$ and $s_i$ the reflections $s_{i,k_i}'$ and $s_{i,k_i}$.
Then the following hold.

\begin{enumerate}[leftmargin=*]
\item  For each $i=1,\dots ,r$, putting $M'=s_i'M$ and writing $\widetilde{\GL_{d_\cdot ,\cdot }(\Fs)'}$ for the factors of $\tM'$,
  conjugation by $s_i'$  induces isomorphisms $\widetilde{\GL_{d_{i},i}(\Fs)}\rightarrow\widetilde{\GL_{d_{i},i+1}(\Fs)'}$ and $\widetilde{\GL_{d_{i+1},i+1}(\Fs})\rightarrow\widetilde{\GL_{d_{i+1},i}(\Fs)'}$ that induce the identity on $\GL_{d_{i}}(\Fs)$ and $\GL_{d_{i+1}}(\Fs)$ respectively and are given via the commutative diagrams
  \begin{equation*}
    \begin{tikzcd}
      \widetilde{\GL_{d_{i},i}(\Fs)} \ar[r,"\sim"] \ar[d]  & \widetilde{\GL_{d_{i},i+1}(\Fs)}' \ar[d]\\
      \tM \ar[r,"s_i' "]                              & \tM'
    \end{tikzcd}
  \end{equation*}
  
  and
  
   \begin{equation*}
    \begin{tikzcd}
      \widetilde{\GL_{d_{i+1},i+1}(\Fs)} \ar[r,"\sim"] \ar[d]  & \widetilde{\GL_{d_{i+1},i}(\Fs)'} \ar[d]\\
      \tM \ar[r,"s_i' "]                              & \tM'
    \end{tikzcd}
  \end{equation*}

Conjugation by $s_i'$ leaves invariant all other $\widetilde{\GL_{d_{i'},i'}(\Fs)}$, $i'\ne i, i+1$ and also $\tH$.

\item Suppose now $d\ne 0$ or that $G$ is not the special orthogonal group or that $d_i$ is even.

Conjugation by $s_i$ induces an automorphism $^{\bT_i}$ on $\widetilde{\GL_{d_{i},i}(\Fs)}$ that induces the outer automorphism on $\GL_{d_{i},i}(\Fs)$ via the commutative diagram
  \begin{equation*}
    \begin{tikzcd}
      \widetilde{\GL_{d_{i},i}(\Fs)} \ar[r,"\sim"] \ar[d]  & \widetilde{\GL_{d_{i},i}(\Fs)} \ar[d]\\
      \tM \ar[r,"s_{i}"]                              & \tM
    \end{tikzcd}
  \end{equation*}

It leaves invariant all other $\widetilde{\GL_{d_{i'},i'}(\Fs)}$, $i'\ne i$.

If $G$ is not a special even orthogonal group or $d_i$ is even, conjugation induces also an action $^{\bT_i}$ on $\tH_d$ by a non-genuine character\footnote{This character is always trivial when $G$ is a symplectic group.} $\tH_d\rightarrow\bmu_{\bm}$.

If $G$ is a special even orthogonal group and $d_i$ is odd, then conjugation by $s_{i}$ induces an automorphism $^{\bT_i}$ on $\tH_d$ that induces the outer automorphism on $H_d$. 

\item Suppose $d=0$ and that $G$ is an even special orthogonal group.

\begin{enumerate}
\item Suppose $d_i$, $d_{i+1}$ odd. Write $M'= {s_{i}' }^+M$. Conjugation by ${s_{i}'} ^+$ induces isomorphisms $\widetilde{\GL_{d_{i},i}(\Fs)}\rightarrow\widetilde{\GL_{d_{i},i+1}(\Fs)'}$ and $\widetilde{\GL_{d_{i+1},i+1}(\Fs})\rightarrow\widetilde{\GL_{d_{i+1},i}(\Fs)'}$ given by commutative diagrams analogous to the ones in (1) that induce the outer automorphisms on $\GL_{d_{i},i}(\Fs)$ and $\GL_{d_{i+1},i+1}(\Fs)$ respectively.

\item Suppose $d_i, d_{i'}$ odd with $i<i'$. Put $t_{i,i'}=s_{i}s_{i'}$. Conjugation by $t_{i,i'}$ induces automorphisms denoted  $^{\bT_{i,i'}}$ on $\widetilde{\GL_{d_{i},i}(\Fs)}$ and $\widetilde{\GL_{d_{i'},i'}(\Fs)}$ given by commutative diagrams analogous to the ones in (2) that induce the outer automorphisms on $\GL_{d_{i},i}(\Fs)$ and $\GL_{d_{i'},i'}(\Fs)$. It leaves invariant  all other $\widetilde{\GL_{d_{i''},i''}(\Fs)}$, $i''\ne i, i'$. In addition, if $i''<i$ (resp. $i<i''$), then $^{\bT_{i,i'}}$ equals $^{\bT_{i'',i}}$ (resp. $^{\bT_{i,i''}}$) on $\widetilde{\GL_{d_{i},i}(\Fs)}$.
\end{enumerate}

\end{enumerate}
\end{prop}

\begin{rmk}\label{rmk:identites Weyl group-action}
Let us explain the effect on the elements of $\tM$ that follows from \ref{prop:decomposed-Levis}: write $\alpha_i'=\alpha_{i,k_i}'$, $\alpha_{i;i'}'=\alpha_i'+\cdots +\alpha_{i'}'$, ${\alpha_{i;i'}' }^+=\alpha_i'+\cdots +\alpha_{i'-1}'+{\alpha_{i'}'} ^+$ and denote here $s_{i,i'}'$ (resp.
${s_{i,i'}'} ^+$) the corresponding simple reflection in $M_{\alpha_{i;i'}'}$ (resp.
$M_{{\alpha_{i;i'}'} ^+}$).

Write $\varphi_{i,i+1}$ for the isomorphism induced by $s_i'$ in (1) in the first diagram, $\varphi_{i+1,i}$ for the isomorphism induced by $s_i'$ in (1) in the second diagram, and $\varphi_{i,i'}$ for the isomorphism obtained by composition taking possibly the conjugated Levi  subgroups.
Remark that $\varphi_{i+1,i}$ is the inverse of $\varphi_{i,i+1}$ when $s_i'M=M$. Otherwise, the inverse of $\varphi_{i,i+1}$ is $\varphi_{i+1,i}'$, where the $'$ means that it is taken relative to the conjugated Levi subgroup $M'=s_i'M$.

Let $m\in\tM$ and choose $x_i\in \widetilde{\GL}_{d_{i},i}(\Fs)$, $i=1,\dots r$, $h\in\tH_d$ such that $m=x_1\cdots x_r h$, $h=1$ if $d=0$ and $G$ is the special even orthogonal group.

Then, one has the following, if the objects are defined:  %

(1) $s_{i,i'}' m=x_1\cdots x_{i-1}\varphi_{i',i}(x_{i'})x_{i+1}\cdots x_{i'-1}\varphi_{i,i'}(x_i)x_{i'+1}\cdots x_rh$ as element of $s_{i,i'}'\tM$ using \ref{prop:decomposed-Levis} (1);

(2) $s_im=x_1\cdots x_{i-1}\varphi_{r,i}(^{\bT_r'}\varphi_{i,r}(x_{i}))x_{i+1}\cdots x_r\ ^{\bT_r'}h$ as elements of $\tM$ using \ref{prop:decomposed-Levis} (2) and $^{\bT_r'}$ taken relatively to the conjugate Levi subgroup;

(3) $t_{i,i'}m=x_1\cdots x_{i-1}\varphi_{r-1,i}(^{\bT_{r-1,r}'}\varphi_{i,r-1}(x_{i}))x_{i+1}\cdots x_{i'-1}\varphi_{r,i'}(^{\bT_{r-1,r}'}\varphi_{i',r} (x_{i'}))$ $x_{i'+1}\cdots x_r$ as elements of $\tM$ using \ref{prop:decomposed-Levis} (3) and $^{\bT_{r-1,r}'}$ taken relatively to the conjugate Levi subgroup;

(4) ${s_{i,i'}' }^+ m=x_1\cdots x_{i-1}\varphi_{i',i}(^{\bT_{i,i'}}x_{i'})x_{i+1}\cdots x_{i'-1}\varphi_{i,i'}(^{\bT_{i,i'}}x_i)x_{i'+1}\cdots x_r$ as element of ${s_{i,i'}'} ^+\tM$ using part (1) and the identity ${s_{i,i'}'} ^+=s_{i,i'}' t_{i,i'}$.
\end{rmk}

\begin{proof} By definition, $\tM$ is generated as a group by the $\widetilde{\GL_{d_{i},i}(\Fs)}$ along with $\tH$.
The fact that $\tM$ decomposes implies that the $\widetilde{\GL_{d_{i},i}(\Fs)}$ commute with each other and with $\tH$.
In particular, the homomorphism
\begin{align*}
\widetilde{\GL_{d_{1},1}(\Fs)}\times\cdots\times\widetilde{\GL_{d_{r},r}(\Fs)}\times\tH_d\rightarrow\tM
\end{align*} 
is surjective.
It is clear that it factors through $$\widetilde{\GL_{d_{1},1}(\Fs)}\times_{\bmu_{\bm}} \cdots\times_{\bmu_{\bm}}\widetilde{\GL_{d_{r},r}(\Fs)}\times_{\bmu_{\bm}}\tH_d\rightarrow\tM.$$ It remains to show that the latter homomorphism is injective.
By compatibility with the projection $\bp $, one sees that the composition with $\bp $ has kernel $\bmu_{\bm}$, which means that the homomorphism is injective as the restriction to $\bmu_{\bm}$ is just the embedding.

For (1), one can see $s_i'$ as an element of the standard Levi subgroup $\tM_{\alpha_{i}'}$ which is canonically isomorphic to
$$\widetilde{\GL_{d_{1},1}(\Fs)}\times_{\bmu_{\bm}} \cdots\times_{\bmu_{\bm}} \widetilde{\GL_{d_{i}+d_{i+1},i}(\Fs)}\times_{\bmu_{\bm}} \cdots\times_{\bmu_{\bm}} \widetilde{\GL_{d_{r},r}(\Fs)} \times_{\bmu_{\bm}} \tH_d$$ and contains $\tM$ with $ \widetilde{\GL_{d_{i}+d_{i}}(\Fs)}$ containing $\widetilde{\GL_{d_i,i}}\times_{\bmu_{\bm}}\widetilde{\GL_{d_i, i+1}}$.
This implies that conjugation by $s_{i}'$ leaves invariant all other $\widetilde{\GL_{d_{i'},i'}(\Fs)}$, $i'\ne i, i+1$ and also $\tH$.

\null For (2), $s_r$ can be seen as an element of the last factor of the standard Levi subgroup $\tM_{\alpha_{r}}$ that is given by $\widetilde{\GL_{d_{1},1}(\Fs)}\times_{\bmu_{\bm}} \cdots\times_{\bmu_{\bm}} \widetilde{\GL_{d_{r-1}, r-1}(\Fs)}\times_{\bmu_{\bm}} \tH_{d_r+d}$ and contains $\tM$.
So the elements of $\widetilde{\GL_{d_{i},i}(\Fs)}$, $i\ne r$, commute with 
$s_{r}$, while the commutative diagram shows that conjugation by $s_{r}$ induces on $\widetilde{\GL_{d_{r},r}(\Fs)}$ an automorphism that lifts the outer automorphism of $\GL_{d_{r},r}(\Fs)$.
As $s_{r}$ commutes with the elements of $H_d$ if $d_r$ is even or $G$ is not a special orthogonal group, conjugation by $s_{r}$ induces by \ref{lemma:conj} on $\tH_d$ a $\bmu_{\bm}$-valued character that factors through $H_d$ and is smooth.
If $G$ is a special orthogonal group and $d_r$ is odd, then $s_{r}$ induces on $\tH_d$ an automorphism which induces the outer automorphism on $H_d$.

The case of general $s_{i}$ can be reduced to the case $i=r$ by permuting the $i$th factor $\widetilde{\GL_{d_{i},k_i}(\Fs)}$ successively with the  $\widetilde{\GL_{d_{i'},i'}(\Fs)}$ which lie on the way to $\tH_d$ using part (1), conjugating by the analogue of $s_{r}$ and permuting back.
Each permutation on the way to $\tH_d$ corresponds to a conjugation by the Weyl group element corresponding to a maximal Levi subgroup $\GL_{d_i}(\Fs)\times \GL_{d_{i'}}(\Fs)$ inside $\GL_{d_i+d_{i'}}(\Fs)$ and on the way back to the Weyl group element of the conjugate Levi which is its inverse.
So this leaves invariant the  $\widetilde{\GL_{d_{i'},i'}(\Fs)}$ different from $\widetilde{\GL_{d_i,i}(\Fs)}$ while acting on $\widetilde{\GL_{d_i,i}(\Fs)}$ and $\tH_d$ as previously.
Remark that the  precise actions may depend on $i$.

\null For (3) (a), suppose now that $G$ is the even special orthogonal group, that $d=0$, and $d_i$, $d_{i+1}$ odd.
If $i=r-1$ then ${s_{r-1}' }^+$ can be realized in $H_{d_{r-1}+d_r}$ inducing the outer automorphisms on $\GL_{d_{r-1}}(\Fs)$ and $\GL_{d_r}(\Fs)$ as indicated and leaving the other $\widetilde{\GL_{d_{i'},i'}(\Fs)}$, $i'\ne r-1, r$ invariant.
For $i<r-1$, one has to conjugate as in case (2), but this still leaves invariant  the other $\widetilde{\GL_{d_{i'},i'}(\Fs)}$.

For (3) (b), assume first $i=r-1$ and $i'=r$.
Then, $t_{i,i'}$ can be realized in $\tH_{d_{r-1}+d_r}$ and the rest follows as previously.
Otherwise, analogous to (2) and (3), one has to conjugate and concludes as in (2), (3).

To show the last property of (3) (b), compose $^{\bT_{i'',i}}$ (resp.
$^{\bT_{i,i''}}$) with $^{\bT_{i,i'}}$, and $^{\bT_{i,i'}}$ with itself.
This is the same as to conjugate with the products $t_{i,i''}t_{i,i'}=s_{i''}s_{i'}s_i^2$ and $t_{i,i'}^2=s_i^2s_{i'}^2$.
But, $s_i^2$ is an element of $M$, so the action on $\widetilde{\GL_{d_i,i}}(\Fs)$ is in both cases the same inner automorphism, while $s_{i'}^2$ and $s_{i'}s_{i''}$ leave $\widetilde{\GL_{d_i,i}}(\Fs)$ invariant by what preceded.
\end{proof}

\begin{rmk}\label{rmk:identif decomposed-Levis}
  In the setting of \ref{prop:str-Weyl-group} with the assumption that the $d_i$ are distinct, we will now use double indexation $\GL_{d_{i},j}(\Fs)$, $i\in\{1,\dots ,r\}$, $j\in\{1,\dots ,k_i\}$ and consider on the set of pairs $(d_i,j)$ in the lexicographic order.
  
  It follows from \ref{prop:decomposed-Levis} that the covers $\widetilde{\GL_{d_{i},j}(\Fs)}$, $j=1,\dots ,k_i$, are all isomorphic.
Denote by $\varphi_{i;j,j+1}$ the isomorphism $\widetilde{\GL_{d_{i},j}(\Fs)}\rightarrow\widetilde{\GL_{d_{i},j+1}(\Fs)}$ induced by conjugation by $s_{i,j}$, by $\varphi_{i;j+1,j}$ its inverse and, in general, by $\varphi_{i;j,j'}$ the isomorphism $\widetilde{\GL_{d_{i},j}(\Fs)}\rightarrow\widetilde{\GL_{d_{i},j'}(\Fs)}$ obtained from the previous ones by composition.
  
Let us put $\widetilde{\GL_{d_{i}}(\Fs)}=\widetilde{\GL_{d_{i},1}(\Fs)}$ for $i=1,\dots ,r$.
  
Using the isomorphisms $\varphi_{i;j,1}:\widetilde{\GL_{d_{i},j}(\Fs)}\rightarrow\widetilde{\GL_{d_{i},1}(\Fs)}$ in the above notation, we will identify $\tM$, $M$ written as in \ref{prop:str-Weyl-group}, with
\begin{align*}
  \widetilde{\GL_{d_{1}}(\Fs)}\times_{\bmu_{\bm}} \cdots\times_{\bmu_{\bm}}\widetilde{\GL_{d_{1}}(\Fs)} \times_{\bmu_{\bm}}\cdots\times_{\bmu_{\bm}}\widetilde{\GL_{d_{r}}(\Fs)}\times_{\bmu_{\bm}}\cdots\times_{\bmu_{\bm}}\widetilde{\GL_{d_{r}}(\Fs)}\\ \times_{\bmu_{\bm}}\tH_d\rightarrow\tM
\end{align*}
and call this \it  a realization of $\tM$\rm .

A representation
$$\widetilde{\rho}_{1,1} \otimes_{\bmu_{\bm}} \cdots \otimes_{\bmu_{\bm}} \widetilde{\rho}_{1,k_{1}}\otimes_{\bmu_{\bm}}  \cdots \otimes_{\bmu_{\bm}} \widetilde{\rho}_{r,1} \otimes_{\bmu_{\bm}} \cdots \otimes_{\bmu_{\bm}} \widetilde{\rho}_{r,k_{r}}\otimes_{\bmu_{\bm}} \widetilde{\tau}$$
of
\begin{align*}
  \widetilde{\GL_{d_{1}}(\Fs)}\times_{\bmu_{\bm}} \cdots\times_{\bmu_{\bm}}\widetilde{\GL_{d_{1}}(\Fs)} \times_{\bmu_{\bm}}\cdots\times_{\bmu_{\bm}}\widetilde{\GL_{d_{r}}(\Fs)}\times_{\bmu_{\bm}}\cdots\times_{\bmu_{\bm}}\widetilde{\GL_{d_{r}}(\Fs)}\\ \times_{\bmu_{\bm}}\tH_d\rightarrow\tM.
\end{align*}
will then be identified with a representation $\widetilde{\sigma }$ of $\tM$ in the following way: write $\tm\in\tM$ as a product $x_{1,1}\cdots x_{1,k_1}\cdots x_{r,1}\cdots x_{r,k_r}h$ with $x_{i,j}\in\widetilde{\GL_{d_{i},j}(\Fs)}$ and $h\in\tH_d$.

Then, $\widetilde{\sigma }(m)$ is the tensor product of the $\widetilde{\rho}_{i,j}(\varphi_{i;j,1}(x_{i,j}))$ and $\widetilde{\tau}(h)$.

We will call it the \it decomposition of $\widetilde{\sigma }$ according to the given realization of $\tM$.
\rm \end{rmk}

\begin{prop} \label{prop:action-on-representations} Let $\widetilde{\sigma}$ be a supercuspidal representation of $\tM$ with decomposition $$\widetilde{\rho}_{1,1} \otimes_{\bmu_{\bm}} \cdots \otimes_{\bmu_{\bm}} \widetilde{\rho}_{1,k_{1}}\otimes_{\bmu_{\bm}}  \cdots \otimes_{\bmu_{\bm}} \widetilde{\rho}_{r,1} \otimes_{\bmu_{\bm}} \cdots \otimes_{\bmu_{\bm}} \widetilde{\rho}_{r,k_{r}}\otimes_{\bmu_{\bm}} \widetilde{\tau}$$
according to the realization
\begin{align*}
  \widetilde{\GL_{d_{1}}(\Fs)}\times_{\bmu_{\bm}} \cdots\times_{\bmu_{\bm}}\widetilde{\GL_{d_{1}}(\Fs)} \times_{\bmu_{\bm}}\cdots\times_{\bmu_{\bm}}\widetilde{\GL_{d_{r}}(\Fs)}\times_{\bmu_{\bm}}\cdots\times_{\bmu_{\bm}}\widetilde{\GL_{d_{r}}(\Fs)}\\ \times_{\bmu_{\bm}}\tH_d\rightarrow\tM ,
\end{align*}
where the $d_{i}$ are supposed to be distinct.

Then,

\begin{enumerate}
\item   For all $i=1,\ldots ,r$ and  $j=1,\dots ,k_i-1$, the decomposition of the representation $s_{i,j}\widetilde{\sigma}$ of $\tM$ is obtained from the one of $\widetilde{\sigma}$ by permuting $\widetilde{\rho}_{i,j}$ and $\widetilde{\rho}_{i,j+1}$. \label{item:8}

\item   Consider $^{\bT_{(i',j')}}$ in the notations of \ref{prop:decomposed-Levis}(2), replacing the indexation in \ref{prop:decomposed-Levis} by the double-indexation with the lexicographic order.
  For $i'=1,\dots ,r$, $j'=1,\dots ,k_{i'}$, let us denote the factors of the decomposition of $^{\bT_{(i',j')}}\widetilde{\sigma}$ via the above realization of $\tM$ by $^{\bT_{(i',j')}}\widetilde{\rho}_{i,j}$.
  One has $^{\bT_{(i',j')}}\widetilde{\rho}_{i,j}=\widetilde{\rho}_{i,j}$ if $(i,j)\ne (i',j')$ and, in addition, $^{\bT_{(i',j')}}$ does not depend on the choice of $j'$ in the sense that, if  $j_1'\in\{1,\dots , k_{i'}\}$ is such that $\widetilde{\rho}_{i',j'}=\widetilde{\rho}_{i',j_1'}$  then $^{\bT_{(i',j')}}\widetilde{\rho}_{i',j'}= ^{\bT_{(i',j'_1)}}\widetilde{\rho}_{i',j_1'}$. \label{item:9}

\item
    Assume $\tG$ to be a special even orthogonal group and $d=0$.
  Let $i', i'' = 1,\ldots ,r$, $j' = 1, \ldots ,k_{i'}$ and $j'' = 1, \ldots ,k_{i''}$ be such that $d_{i'}, d_{i''}$ are odd and $(i',j')<(i'',j'')$ in the lexicographic order.
  Consider $^{\bT_{(i',j'),(i'',j'')}}$ in the notations of \ref{prop:decomposed-Levis}(3), replacing the indexation in \ref{prop:decomposed-Levis} by the double-indexation with the lexicographic order.
  Denote the factors of the decomposition of $^{\bT_{(i'j'),(i'',j'')}} \widetilde{\sigma}$ via the above realization of $\tM $ by $^{\bT_{(i',j'),(i'',j'')}}\widetilde{\rho}_{i,j}$.

  One has $^{\bT_{(i',j'),(i'',j'')}}\widetilde{\rho}_{i,j}=\widetilde{\rho}_{i,j}$ if $(i,j)\ne \{(i',j'), (i'',j'')\}$.
  In addition, $^{\bT_{(i',j'),(i'',j'')}}$ does not depend on the choice of $j'$ and $j''$ in the sense that, if  $j_1' \in \{ 1, \ldots ,k_{i'}  \}$ and $j_1'' \in \{ 1, \ldots ,k_{i''}  \}$ satisfying $(i',j'_1)<(i'',j''_1)$ are such that $\widetilde{\rho}_{i',j'}=\widetilde{\rho}_{i',j_1'}$ (resp.
  $\widetilde{\rho}_{i'',j''}=\widetilde{\rho}_{i'',j_1''}$), then $^{\bT_{(i',j'),(i'',j'')}}\widetilde{\rho}_{i',j'}= ^{\bT_{(i',j_1'),(i'',j'')}}\widetilde{\rho}_{i',j_1'}$ (resp.
  $^{\bT_{(i',j'),(i'',j'')}}\widetilde{\rho}_{i'',j''}= ^{\bT_{(i',j'),(i'',j'')}}\widetilde{\rho}_{i'',j_1''}$). 
\label{item:10}
\end{enumerate}
\end{prop}

\begin{rmk} \label{rmk:definition T} Part \eqref{item:9} of the proposition allows us to write just $^{\bT_i}$ for the action on the components $\widetilde{\rho}_{i,j}$ with $j=1,\dots ,k_i$ induced by the actions of ${\bT_{(i,j)}}$ on $\widetilde{\sigma}$.
Similarly, for \eqref{item:10}, we can just write $^{\bT_{i,i'}}$ for the action - up to isomorphism - on the components $\widetilde{\rho}_{i,j}$ and $\widetilde{\rho}_{i',j'}$ induced by the action of $^{\bT_{(i,j),(i',j')}}$ on $\widetilde{\sigma }$, $j=1,\dots ,k_i$ and $j'=1,\dots ,k_{i'}$.
(We use here also \ref{prop:decomposed-Levis} (3) (b).)
This includes the case $i'=i$ for $k_i>1$.

Both, $^{\bT_i}$ and $^{\bT_{i,i'}}$, stand for automorphisms of $\tM$ which introduce outer automorphisms on $\GL_{d_i}(\Fs)$ and $\GL_{d_{i'}}(\Fs)$ respectively via the commutative diagram defined by the projections.
Additionally, we will write $^{\bT_i}$ also for $^{\bT_{i,i}}$, although the latter always involves simultaneous action on two different copies of $\widetilde{\GL_{d_i}(\Fs)}$.
\end{rmk}

\begin{proof} \eqref{item:8} is a direct computation from the definition.

  For \eqref{item:9}, the first part follows from \ref{prop:decomposed-Levis}(2). Regarding the independency of $j'$, consider $(i',j')=(i,j)$.
   Using the identities in \ref{rmk:identites Weyl group-action} with \ref{rmk:identif decomposed-Levis}, one sees that, for $x_{i,j}\in\widetilde{\GL_{d_{i},j}}(\Fs)\subseteq\tM$,  $$^{\bT_{(i,j)}}x_{i,j}=\varphi_{i;1,j}\ ( ^{\bT_{(i,1)}}(\varphi_{i;j,1}(x_{i,j}))).$$
It follows that the $(i,j)$-th component of the tensor product  of automorphisms which corresponds to $^{\bT_{(i,j)}}\widetilde{\sigma }(x_{i,j})$ (the other components are the identity) is
$${\rho_{i,j}}(\varphi_{i;j,1}(^{\bT_{(i,j)}}x_{i,j}))={\rho_{i,j}}( ^{\bT_{(i,1)}}(\varphi_{i;j,1}(x_{i,j}))).$$
This proves that the $(i,j)$-component of  the decomposition of $^{\bT_{(i,j)}}\widetilde{\sigma }$ is $^{\bT_{(i,1)}}\widetilde{\rho_{i,j}}$.
It does not depend on $j$.

For \eqref{item:10}, the first part follows from \ref{prop:decomposed-Levis}(3)b). For the second part, it is enough to consider $(j',j'')=(1,k_{i''})$.
The rest will follow from this.
Using the identities in  \ref{rmk:identites Weyl group-action} with \ref{rmk:identif decomposed-Levis}, one sees that, for $x_{i',j_1'}\in\widetilde{\GL_{d_{i'},j_1'}}(\Fs)$ (resp.
for $x_{i'',j_1''}\in\widetilde{\GL_{d_{i''},j_1''}}(\Fs)$), $$^{\bT_{(i',j_1'),(i'',j'')}}x_{i',j'_1}=\varphi_{i';1,j_1'}\ ( ^{\bT_{(i',1),(i'',j'')}}(\varphi_{i';j_1',1}(x_{i',j'_1})))$$ $$(resp.
^{\bT_{(i',j'),(i'',j_1'')}}x_{i'',j_1''}=\varphi_{i'';k_{i''},j_1''}\ ( ^{\bT_{(i',j'),(i'',k_{i''})}}(\varphi_{i'';j_1'',k_{i''}}(x_{i'',j_1''})))).$$ It follows then as for \eqref{item:9} that the $(i',j_1')$-th component (resp.
$(i'',j_1'')$-th component) of $^{\bT_{(i',j_1'),(i'',j'')}}\widetilde{\sigma}$ (resp.\ $^{\bT_{(i',j'),(i'',j''_1)}}\widetilde{\sigma}$) in its decomposition is $^{\bT_{(i',1),(i'',j'')}}\widetilde{\rho}_{i',j'_1}$ $$(resp.
\widetilde{\rho}_{i'',j_1''}\circ \varphi_{i'';k_{i''},1}\circ ^{\bT_{(i',j'),(i'',k_{i''})}}\circ \varphi_{i'';1,k_{i''}}).$$ It does not depend on $j_1'$ (resp.
$j_1''$).

\end{proof}

\begin{prop} \label{prop:Sigma-O-mu-cover-decomposed-Levis}
  Assume that the assumptions in \ref{prop:decomposed-Levis} are satisfied w.r.t.
$\tM $ and $\tG$.
Let $\widetilde{\sigma}$ be a supercuspidal representation of $\tM$.
\begin{enumerate}[leftmargin=*]
  \item Using the realization \ref{rmk:identif decomposed-Levis} and without changing the Bernstein component\footnote{For $G$ of type $D$ and $d=0$, the fact that the last simple root is not a root of $M$ may require an outer automorphism of $G$ or even an isomorphism with another cover of $G$ inducing an outer automorphism on $G$ (see also Proposition~\ref{prop:str-Weyl-group} (1)).}, we can assume that $\tM$ and $\widetilde{\sigma}$ have decompositions
    \begin{align*}
      \begin{array}{c}
        \widetilde{\GL_{d_{1}} (\Fs)}\times_{\bmu_{\bm}} \cdots \times_{\bmu_{\bm}} \widetilde{\GL_{d_{1}} (\Fs)}
        \times_{\bmu_{\bm}} \cdots \times_{\bmu_{\bm}} \widetilde{\GL_{d_{r}} (\Fs)}\times_{\bmu_{\bm}} \cdots\times_{\bmu_{\bm}} \widetilde{\GL_{d_{r}}(\Fs)}
        \times_{\bmu_{\bm}} \tH_{d},\\ 
        \text{and} \\      
        \widetilde{\rho}_{1,1} \otimes_{\bmu_{\bm}} \cdots \otimes_{\bmu_{\bm}} \widetilde{\rho}_{1,k_{1}}
        \otimes_{\bmu_{\bm}}  \cdots \otimes_{\bmu_{\bm}} \widetilde{\rho}_{r,1} \otimes_{\bmu_{\bm}} \cdots \otimes_{\bmu_{\bm}} \widetilde{\rho}_{r,k_{r}}
        \otimes_{\bmu_{\bm}} \widetilde{\tau},
      \end{array}
    \end{align*}
    allowing consecutive $d_i$ to be equal, where  $H_{d}$ is a reductive group of the same type as $G$ such that $d\neq 1$ if $G$ is of type $D$, and satisfy the following conditions:
    \begin{enumerate}
    \item For $i=1,\ldots ,r$,
     \begin{enumerate}
     \item if $G$ is not a special even orthogonal group, $d\ne 0$ or $d_i$ is even, the representations $\widetilde{\rho}_{i,j}$, $j=1,\ldots , k_{i}$, are all isomorphic;
     \item otherwise, the representations $\widetilde{\rho}_{i,j}$, $j=1,\ldots , k_{i}$, are all isomorphic except for at most one $i$.
In this case, one may assume that all representations $\widetilde{\rho}_{i,j}$ with $j\leq k_i-1$ are isomorphic while $\widetilde{\rho}_{i,k_i}\simeq\ ^{\bT_{i}} \widetilde{\rho}_{i,k_i-1}$.
     \end{enumerate}
    \item For $i=1,\ldots ,r$; $i'=1,\ldots ,r$, if defined, $^{\bT_i}\widetilde{\rho}_{i,k_{i}}$ (resp.
$^{\bT_{i,i'}}\widetilde{\rho}_{i,k_{i}}$) and $\widetilde{\rho}_{i,k_{i}} $ are isomorphic if in the same inertial orbit.
    \item  For  $i,i' = 1, \ldots , r$; $i'\ne i$,
    \begin{enumerate}
     \item if $G$ is not a special even orthogonal group, $d\ne 0$ or $d_i$ is even, $\widetilde{\rho}_{i',k_{i'}}$ (resp.
$^{\bT_i}\widetilde{\rho}_{i',k_{i'}}$ if defined) and $\widetilde{\rho}_{i,k_{i}}$ are never in the same inertial orbit if $i'\ne i$;
      \item otherwise, for $i,i',i''$  all distinct, $d_i$, $d_{i'}$ and $d_{i''}$ odd, $\widetilde{\rho}_{i',k_{i'}}$ (resp.
$^{\bT_{i',i''}}\widetilde{\rho}_{i',k_{i'}}$, $^{\bT_{i',i''}}\widetilde{\rho}_{i'',k_{i''}}$)  and $\widetilde{\rho}_{i,k_{i}}$ are never in the same inertial orbit and the same is true for $\widetilde{\rho}_{i,1}$ instead of  $\widetilde{\rho}_{i,k_{i}}$ including the exceptional case in (1)(a)(ii).
      \end{enumerate}
     \end{enumerate}

  \item  The components $\Sigma_{\cO,\mu, i}$'s of the root system $\Sigma_{\cO,\mu}$ and a part $R(\cO)_{i}$ of the $R$-group are determined as follows.

Suppose first that $\mu^{\tM_{\alpha _{i,1}}}$  is not regular or $k_i= 1$.

\begin{enumerate}
\item Assume that $d\neq 0$ or $G$ is of type $B$.
Then we have the following.
        \begin{enumerate}
        \item If $\mu^{\tM_{\alpha_{i,k_{i}}}}$ has a pole, then $^{\bT_i}\widetilde{\rho}_{i,k_{i}} \simeq\widetilde{\rho}_{i,k_{i}} $, $^{\bT_i}\widetilde{\tau}\simeq\widetilde{\tau}$ and a base for $\Sigma_{\cO,\mu,i}$ is
          \begin{equation}\label{eq:Sigma-O-mu-pole-1}
            \{ \alpha_{i,1}, \alpha_{i,2}, \ldots, \alpha_{i,k_{i}} \}
          \end{equation}
          This is   a root system of type $B_{k_{i}}$.
          
          We put $R(\cO)_i=\{1\}$.
        \item   Suppose that $\mu^{\tM_{\alpha_{i,k_{i}}}}$ is regular and 
          $^{\bT_i}\widetilde{\rho}_{i,k_{i}}\isom \widetilde{\rho}_{i,k_{i}}$.
          Then, if $k_{i}\neq 1$,
 a base for $\Sigma_{\cO,\mu,i}$ is
          \begin{equation}\label{eq:Sigma-O-mu-no-pole-dual-1}
            \{ \alpha_{i,1}, \alpha_{i,2}, \ldots, \alpha_{i,k_{i}-1},  \alpha_{i,k_{i}-1} + 2\alpha_{i,k_{i}}\}.
          \end{equation}
          This is a root system of type $D_{k_{i}}$.
          If $k_{i}= 1$, then $\Sigma_{\cO,\mu,i}$ is empty.
           
       We put $R(\cO)_i=\langle s_{i,k_i}\rangle$,   if $^{\bT_i}\widetilde{\tau}\simeq\widetilde{\tau}$.
       Otherwise, we put $R(\cO)_i=\{ 1 \}$.
           \item Otherwise, a base for $\Sigma_{\cO,\mu,i}$ is
          \begin{equation}\label{eq:Sigma-O-mu-no-pole-not-dual-1}
            \{ \alpha_{i,1}, \alpha_{i,2}, \ldots, \alpha_{i,k_{i}-1}\}.
          \end{equation}
          This is a root system of type $A_{k_{i}-1}$.
          
          We put $R(\cO)_i=\{1\}$.
    \end{enumerate}

    \item Next assume that $d=0$ and $G$ is not of type $B$.
    Then we have the following.
    \begin{enumerate}
        \item If $\mu^{\tM_{2\alpha_{i,k_{i}}}}$ has a pole, then a base for $\Sigma_{\cO,\mu,i}$ is
          \begin{equation}\label{eq:Sigma-O-mu-pole-2}
            \{ \alpha_{i,1}, \alpha_{i,2}, \ldots, 2\alpha_{i,k_{i}} \}
          \end{equation}
          This is   a root system of type $C_{k_{i}}$.
          
          We put $R(\cO)_i=\{1\}$
        \item   Suppose that  $\mu^{\tM_{2\alpha_{i,k_{i}}}}$ is regular and  $^{\bT _i} \widetilde{\rho}_{i,k_{i}}\isom \widetilde{\rho}_{i,k_{i}}$.
          Then, if $k_{i}\neq 1$,  a base for $\Sigma_{\cO,\mu,i}$ is
          \begin{equation}\label{eq:Sigma-O-mu-no-pole-dual-2}
            \{ \alpha_{i,1}, \alpha_{i,2}, \ldots, \alpha_{i,k_{i}-1},  \alpha_{i,k_{i}-1} + 2\alpha_{i,k_{i}}\}.
          \end{equation}
          This is a root system of type $D_{k_{i}}$.
          If $k_{i}= 1$, then $\Sigma_{\cO,\mu,i}$ is empty.
          
          We put $R(\cO)_i=\langle s_{i,k_i}\rangle$, if $G$ is of type $C$ or $d_i$ is even, otherwise $R(\cO)_i=\{ 1\}$.
        \item Otherwise,
        \begin{enumerate}
        \item if $d_i$ even or $G$ is not a special orthogonal group or $d_i$ odd and $i$ does not satisfy the last part of (1)(a)(ii), a base for $\Sigma_{\cO,\mu,i}$ is
          \begin{equation}\label{eq:Sigma-O-mu-no-pole-not-dual-2}
            \{ \alpha_{i,1}, \alpha_{i,2}, \ldots, \alpha_{i,k_{i}-1}\}.
          \end{equation}
          This is a root system of type $A_{k_{i}-1}$.
          
          We put $R(\cO)_i=\{1\}$.
          \item if $G$ is the special orthogonal group, $d_i$ odd and the last part of (1)(a)(ii) is satisfied for  $i$, a base for $\Sigma_{\cO,\mu,i}$ is
          \begin{equation}\label{eq:Sigma-O-mu-no-pole-not-dual-3}
            \{ \alpha_{i,1}, \alpha_{i,2}, \ldots, \alpha_{i,k_{i}-2}, \alpha_{i,k_{i}-1}^+\}.
          \end{equation}
          This is a root system of type $A_{k_{i}-1}$ and we let $R(\cO)_i=\{1\}$.
          \end{enumerate}
              \end{enumerate}
       \end{enumerate}

     \item If $\mu^{\tM_{\alpha_{i,1}}}$  is regular, then all roots $\ne\alpha_{i,k_i}, 2\alpha_{i,k_i}$ and their linear combinations have to be removed from  $\Sigma_{\cO,\mu,i}$ and the corresponding reflections added to $R(\cO)_i$.
In the (2)(a)(ii), (2)(a)(iii), (2)(b)(ii), (2)(b)(iii)  cases, in particular, this means removing all the roots in the base and thus $\Sigma_{\cO,\mu,i}$ is empty.
In the (2)(a)(i) (resp. (2)(b)(i)) case, $\Sigma_{\cO,\mu,i}$ becomes of type $B_1\times\cdots\times B_1$ (resp. $C_1\times\cdots\times C_1$).
This applies accordingly to the case $k_i=1$.
     \item For the $R$-groups $R(\cO)$, one has the following.
  \begin{enumerate}
   \item Put $I=\{1,\dots , r\}$ and let $R(\cO)_i$  be the groups defined above.
Then, $R(\cO)_I=\prod_{i\in I}R(\cO)_i$ is a normal subgroup of $R(\cO)$.
   \item Assume that $d\ne 0$.
Denote by $I'$ the set of $i$ with $d_i$ odd such that $^{\bT _i} \widetilde{\rho}_{i,k_{i}}\isom \widetilde{\rho}_{i,k_{i}}$ and $^{\bT _i} \widetilde{\tau }\not\isom\widetilde{\tau }$.
Let $R(\cO)^{I'}$ be the subgroup of $W(M,\cO)$ generated by products of $s_{i,k_i}$, $i\in I'$, such that the product $\bT$ of the corresponding $\bT_i$'s satisfies $^{\bT }\widetilde{\tau }\isom\widetilde{\tau }$.
Then, $R(\cO)$ is the semidirect product of $R(\cO)_I$ and $R(\cO)^{I'}$.
    \item Assume that $G$ is the special even orthogonal group and $d=0$.
Denote by $I'$ the set of ordered pairs $(i,i')$ with $d_i$ and $d_{i'}$ odd and $^{\bT _{i,i'}} \widetilde{\rho}_{i,k_{i}}\isom \widetilde{\rho}_{i,k_{i}}$, $^{\bT _{i,i'}}\widetilde{\rho}_{i',k_{i'}}\isom \widetilde{\rho}_{i',k_{i'}}$.
Let $R(\cO)^{I'}$ be the subgroup of $W(M,\cO)$ generated by the products $s_{i,k_i}s_{i',k_{i'}}$, $(i,i')\in I'$.
Then, $R(\cO)$ is the semidirect product of $R(\cO)_I$ and $R(\cO)^{I'}$.
      \end{enumerate}
\end{enumerate}
\end{prop}

\begin{proof}
We can assume that $\tM$ has the form determined in Proposition~\ref{prop:str-Weyl-group} with not necessarily distinct $d_{i}$ after conjugating $\widetilde{\sigma}$ which does not change the Bernstein component, allowing an isomorphism with another cover of $G$ inducing an outer automorphism on $G$ in the case $G$ of type $D$ and $d=0$ to assure that the last simple root is not a root of $M$.

We can permute the factors of $\widetilde{\sigma}$ by \ref{prop:action-on-representations}, putting together the $\widetilde{\rho}_{i,j}$ which are in the same inertial orbit, so that, after possibly twisting by an unramified character, the property (1)(a) is satisfied.
By conjugating and applying $^{\bT_i}$ if defined, we can also make sure that there is no $\widetilde{\rho}_{i',j'}$ with $i'\ne i$ such that $^{\bT_{i'}}\widetilde{\rho}_{i',j'}$ is in the same inertial orbit as $\widetilde{\rho}_{i,j}$.

Concerning (1)(b), assume now that $\widetilde{\rho}_{i,j}$ and $^{\bT_{i}}\widetilde{\rho}_{i,j}$ are in the same inertial orbit.
Remark that $s_{i,k_i}$ lies in a factor isomorphic to $H_{d_i+d}$ of the Levi  subgroup $\tM_{\alpha_{i,k_i}}$.
So, $s_{i,k_i}\widetilde{\sigma }$ has the same components $\widetilde{\rho}_{i',j'}$ as $\widetilde{\sigma }$ except that the $(i,k_i)$-component is now a twist of $\widetilde{\rho}_{i,k_i}$ by some unramified character $\chi_{\lambda\alpha_{i,k_i}}$ which must be unitary.
Replacing $\widetilde{\sigma }$ by $\widetilde{\sigma }_{{\lambda\over 2}\alpha_{i,k_i}}$, this component becomes invariant by $^{\bT_i}$.
Twisting the other components $\widetilde{\rho}_{i,j}$ by the similar character, we come back to property (1)(a).
This gives $\tM$ and $\widetilde{\sigma}$ with properties (1) (a)(i), (b) and (c)(i). In addition, $\widetilde{\rho}_{i,j}$ with distinct $i$ lie in distinct inertial orbits.

Suppose now that $G$ is a special even orthogonal group, $d=0$.
If there is no $i$ with $d_i$ odd and $i'>i$ , such that $^{\bT_{i,i'}}\widetilde{\rho}_{i,k_{i}}$ is in the inertial orbit of $\widetilde{\rho}_{i',k_{i'}}$, then we have finished.
Otherwise, let $i$ be minimal for this property.
By \ref{prop:decomposed-Levis} (3)b) and the fact  that $\widetilde{\rho}_{i,j}$ with distinct $i$ lie in distinct inertial orbits, such an $i'$ is unique.
If $k_{i'}$ (resp. $k_i$) is even, we can apply $^{\bT_{i,i'}}$ successively to all $\widetilde{\rho}_{i',j'}$ (resp.
$\widetilde{\rho}_{i,j}$) and to $\widetilde{\rho}_{i,k_i}$ (resp.
$\widetilde{\rho}_{i',k_{i'}}$) so that all lie in the same inertial orbit.
It remains then to twist either all $\widetilde{\rho}_{i',j'}$ or all $\widetilde{\rho}_{i,j}$ by an unramified character and we are done.

If $k_i$ and $k_{i'}$ are both odd, but there is $i''>i$ different from $i$ and $i'$ with $d_{i''}$ odd, then we apply $^{\bT_{i',i''}}$ successively to all $\widetilde{\rho}_{i',j'}$ and each time to $\widetilde{\rho}_{i'',k_{i''}}$.
After twisting by an unramified character and reordering we are done, remarking that $^{\bT_{i',i''}}\widetilde{\rho}_{i'',k_{i''}}$ may not be in the inertial orbit of $\widetilde{\rho}_{i'',k_{i''}}$. By choice of $i$ above, one remarks however that $^{\bT_{i',i''}}\widetilde{\rho}_{i'',k_{i''}}$ cannot be isomorphic to an $\widetilde{\rho}_{i''',k_{i'''}}$ with $i'''\leq i$. We then reorder and do the previous steps until we arrive again at $\tM$ and $\widetilde{\sigma}$ with properties (1) (a)(i), (b) and (c)(i) and so that in addition $\widetilde{\rho}_{i,j}$ with distinct $i$ lie in distinct inertial orbits..

We pursue until there is either no $i$ with $d_i$ odd and $i'>i$ , such that $^{\bT_{i,i'}}\widetilde{\rho}_{i,k_{i}}$ is in the inertial orbit of $\widetilde{\rho}_{i',k_{i'}}$, or that an $i''$ as above does not exist. The construction makes sure that this process will stop. In the latter case, we may end up in case (a) (ii).
This can be avoided, if $k_i$ or $k_{i'}$ is a sum of multiples of other $k_{i''}$ with $d_{i''}$ odd or if $^{\bT_{i',i''}}\widetilde{\rho}_{i'',k_{i''}}$ is in the inertial orbit of $\widetilde{\rho}_{i'',k_{i''}}$ for some $i''$, but not otherwise.

Finally, if for some $i', i''$ with $d_{i'}, d_{i''}$ odd, one has that $^{\bT_{i',i''}}\widetilde{\rho}_{i'',k_{i''}}$ lies in the inertial orbit of $\widetilde{\rho}_{i'',k_{i''}}$ but both are not isomorphic - remark that it does not depend on $i''$ by the last part of \ref{prop:decomposed-Levis}) -, we can in the same way as for (1)(c)(i) twist  $\widetilde{\rho}_{i'',k_{i''}}$ and the other $\widetilde{\rho}_{i'',j''}$ by an unramified character so that all become isomorphic.

Concerning (2)(a), it is clear that in case (i) all roots $\alpha_{i,j}$ lie in $\Sigma_{\cO ,\mu}$.
In case (ii), only the short root $\alpha_{i,k_i}$ and its conjugates do not lie in $\Sigma_{\cO ,\mu}$, but $s_{i,k_i}$ is in the $R$-group if and only if $^{\bT_i}\tau\simeq\tau$.
\footnote{Remark that by \ref{thm:mu-zero-irreducibility}  the simple reflection $s_{i,k_i-1}^+$ associated to  $\alpha_{i,k_{i}-1} + 2\alpha_{i,k_{i}}$ must consequently always leave invariant $\tau $, even if $^{\bT_i}\tau\not\simeq\ \tau$.
This can be seen independently from the fact that $s_{i,k_i-1}^+=s_{i,k_i}s_{i,k_i-1}s_{i,k_i}$ acts in the same way on $\tH_d$ as $s_{i,k_i}^2$ (see \ref{rmk:identites Weyl group-action} (2)), the latter action being by inner automorphism.}

Finally, in case (iii), only the roots $\alpha_{i,j}$ with $j\ne k_i$ and its linear combinations lie in $\Sigma_{\cO,\mu}$ while the reflections w.r.t.
the other simple roots do not leave $\cO $ invariant and are consequently not in the $R$-group.

Concerning (2)(b), the same holds with some additional remarks in the cases of $d_i$ odd: in the situation (2)(b)(iii)(B), $s_{i,k_i}$ does not leave $M$ invariant by \ref{prop:str-Weyl-group}(5).
As $^{\bT_{i,i}}\widetilde{\rho}_{i,k_{i,i}}$ is isomorphic to $\widetilde{\rho}_{i,k_{i,i}+1}$, $\alpha_{i,k_{i,i}}^+$ is well an element of $\Sigma_{\cO ,i}$.

The point (3) is clear and also (4)(a).

The additional elements of the $R$-group appear when examining the elements in $W(M)$ made explicit in \ref{prop:decomposed-Levis} and examining which do not lie in $W_{\cO}$ but leave $\cO $ and $\Sigma_{\cO,\mu }$ invariant.

To see that $R(\cO)_I$ is a normal subgroup of $R(\cO)$ it is enough to show that it is normalized by $R(\cO)^{I'}$ in the cases (4)(b) and (4)(c) respectively.
In case (4)(b), for $i\in I'$, one is either in situation (2)(a)(ii) or in  situation (3).
In the situation (2)(a)(ii) one has $R(\cO)_i\ne 1$ under the assumptions of (4)(b), and consequently the factor $s_{i,k_i}$ commutes with $R(\cO)_I$ as do the other factors.
In the situation (3) $R(\cO)_i$ is the Weyl group of a root system of type $D_{k_i}$ and this is normalized by $s_{i,k_i}$, while commuting with the $s_{i',k_{i'}}$, $i\ne i'$.
In case (4)(c), one is either in situation (2)(b)(ii) or in (3) for both $i$ and $i'$ and the similar as above applies.
\end{proof}

\begin{rmk}
  Theorem~6.2 of \cite{Kaplan-Lapid-Zou-MR4663361} shows that when the cover is metaplectic in the sense of their paper, the point of reducibility of $\rho_{i,j} |\det |_{\Fz}^{s}\otimes_{\bmu_{\bm}} \rho_{i,j+1}$ is not at $0$ and thus in that case, our assumption that $\mu^{\tM_{\alpha^{-}_{i,1}}}$  is not regular always holds.
\end{rmk}

\begin{rmk}
The statements of Proposition~\ref{prop:Sigma-O-mu-cover-decomposed-Levis} for $k_{i}=1$ suggest that the root systems of type $D_{1}$ and $A_{0}$ should be regarded as the empty set, though, for uniformity of statements, it still counts as a component of $\Sigma_{\cO,\mu}$.
   In addition, root systems of type $B_{1}$, $C_{1}$ and $D_{2}$ degenerate to type $A_{1}$, $A_{1}$, $A_{1}\times A_{1}$ respectively.
   However, for later results, especially those in Section~\ref{sec:hecke-algebras}, we will keep the labels $B_{1}$, $C_{1}$, $D_{1}$ and $D_{2}$ respectively.
\end{rmk}

To finish the section, let us consider the particular case of a general linear group over $F$.
The definition~\ref{defn:prod-of-covers} makes sense for a cover of a general linear group as already remarked there and \ref{defn:mu-prod-of-covers} stands as it is with $H_d$ omitted.
The analogue of \ref{prop:str-Weyl-group} remains valid when keeping only the $\alpha_{i,j}$ and $s_{i,j}$ with $j=1,\dots ,k_i-1$ and one has the analogue of \ref{prop:decomposed-Levis} (1) so that the terminology of \it realization \rm of a Levi subgroup and \it decomposition \rm of a representation of such a Levi subgroup following remark~\ref{rmk:identif decomposed-Levis} applies.

With this, one gets as an analogue of \ref{prop:Sigma-O-mu-cover-decomposed-Levis}:

\begin{prop}\label{prop:Sigma-O-mu-cover-GL-decomposed-Levis}
     Suppose that $G$ is a general linear group and that the standard Levi subgroups of $\tG$ decompose. Let $\widetilde{\sigma}$ be a supercuspidal representation of a standard Levi subgroup $\tM$.
    \begin{enumerate}[leftmargin=*]
    \item Using the realization \ref{rmk:identif decomposed-Levis} and without changing the Bernstein component we can assume that $\tM$ and $\widetilde{\sigma}$ have decompositions
      \begin{align*}
        \begin{array}{c}
          \widetilde{\GL_{d_{1}}(\Fs)} \times_{\bmu_{\bm}} \cdots \times_{\bmu_{\bm}}\widetilde{\GL_{d_{1}}(\Fs)}
          \times_{\bmu_{\bm}} \cdots \times_{\bmu_{\bm}} \widetilde{\GL_{d_{r}}(\Fs)} \times_{\bmu_{\bm}} \cdots \times_{\bmu_{\bm}} \widetilde{\GL_{d_{r}}(\Fs)}\\ \text{and}\\
          \widetilde{\rho}_{1,1} \otimes_{\bmu_{\bm}} \cdots \otimes_{\bmu_{\bm}} \widetilde{\rho}_{1,k_{1}}
          \otimes_{\bmu_{\bm}}  \cdots \otimes_{\bmu_{\bm}} \widetilde{\rho}_{r,1} \otimes_{\bmu_{\bm}} \cdots \otimes_{\bmu_{\bm}} \widetilde{\rho}_{r,k_{r}}
        \end{array}
      \end{align*}
      and can be assumed to satisfy the following conditions:
      \begin{enumerate}
      \item For $i=1,\ldots ,r$ and $j,j'=1,\ldots , k_{i}$, the representations $\widetilde{\rho}_{i,j}$ and $\widetilde{\rho}_{i,j'}$ for $j,j'=1,\ldots , k_{i}$ are all isomorphic.
            \item  For  $i,i' = 1, \ldots , r$ and $i\neq i'$,    $\widetilde{\rho}_{i,k_{i}}$ and $\widetilde{\rho}_{i',k_{i'}}$ are not in the same inertial orbit.
      \end{enumerate}
      
    \item  The irreducible components $\Sigma_{\cO,\mu, i}$ of the root system $\Sigma_{\cO,\mu}$ are indexed by the set $\{1,\dots ,r\}$ with $\Sigma_{\cO,\mu, i}$ of type $A_{k_{i}-1}$  with base $\{ \alpha_{i,1}, \alpha_{i,2}, \ldots, \alpha_{i,k_{i}-1}\}$ if $k_i\neq 1$ and  $\mu^{\tM_{\alpha^{-}_{i,1}}}$  is not regular.
Otherwise, $\Sigma_{\cO,\mu, i}=\emptyset$.
    \item The group $R(\cO)$ is a product of groups $R(\cO)_i$ indexed by the set $\{1,\dots ,r\}$ with $R(\cO)_i=1$ if $k_i=1$ or if $\mu^{\tM_{\alpha^{-}_{i,1}}}$ is not regular.
Otherwise, $R(\cO)_i$ is equal to the symmetric group generated by the simple reflections $\{s_{i,1},\dots , s_{i,k_i-1}\}$.
    \end{enumerate}
\end{prop}

\begin{proof} The proof is analogous to that of Proposition~\ref{prop:Sigma-O-mu-cover-decomposed-Levis} using the remarks above before the statement.
\end{proof}

\section{Restriction of a supercuspidal representation of $\tM$ to $\tM^1$ --- decomposed Levis }
\label{sec:restr-cusp-repr}

  Next we consider the restriction of $\widetilde{\sigma}$  to $\tM^{1}$.
Until now the computation of the full intertwining algebra (for reductive groups) has only been worked out under a certain multiplicity $1$ condition,
which, as is shown in Proposition~\ref{prop:M1-restriction-mult-one}, is satisfied  for covers of classical groups whose Levi subgroups decompose (Definition~\ref{defn:prod-of-covers}).

  \begin{prop}\label{prop:M1-restriction-mult-one} Assume that the Levi subgroups decompose or that $A_M$ is of rank $1$.
    Let $(\widetilde{\sigma}_{1}, \tE_{1})$ be an irreducible component of $\widetilde{\sigma}|_{\tM^{1}}$.
    Let $\tM^{\widetilde{\sigma}} = \{ \tm\in \tM \spacedvert  \lsup{\tm}{{}\widetilde{\sigma}_{1}} \isom \widetilde{\sigma}_{1}\}$, which is a normal subgroup of $\tM$.
    Then $\widetilde{\sigma}_{1}$ can be extended to a representation $\widetilde{\sigma}_{2}$ of $\tM^{\widetilde{\sigma}}$ such that
\begin{equation*}
  \widetilde{\sigma}|_{\tM^{\widetilde{\sigma}}} \isom \bigoplus_{\tm\in \tM/ \tM^{\widetilde{\sigma}}}\lsup{\tm}{{}\widetilde{\sigma}_{2}},
\end{equation*}
where the components $\lsup{\tm}{{}\widetilde{\sigma}_{2}}$'s are pairwise non-isomorphic irreducible representations.
Moreover $\widetilde{\sigma} \isom \Ind_{\tM^{\widetilde{\sigma}}}^{\tM} \widetilde{\sigma}_{2}$ and
\begin{equation*}
  \widetilde{\sigma}|_{\tM^{1}} \isom  \bigoplus_{\tm \in \tM/ \tM^{\widetilde{\sigma}}} \lsup{\tm}{{}\widetilde{\sigma}_{1}},
\end{equation*}
where $\lsup{\tm}{{}\widetilde{\sigma}_{1}}$ are pairwise non-isomorphic irreducible representations.
\end{prop}

\begin{rmk}
  The proof below remains valid for other cases of $G$ whose Levi subgroups can be ``reduced'' to the rank $1$ case.
Proposition~\ref{prop:M1-restriction-mult-one} also implies that the equivalence class of $\ind_{\tM^{1}}^{\tM} \tE_{1}$ does not depend on the choice of the irreducible component $(\widetilde{\sigma}_{1},\tE_{1})$.
\end{rmk}

  \begin{proof}
When $A_{M}$  is of rank $1$,  the proof   in \cite[Proposition~1.16]{Heiermann-MR2827179} carries over  after replacing $A_{M}$ by $\tA_{M}^{\dagger}$, where we recall that $A_{M}^{\dagger} = \{ a^{\bm} \spacedvert a\in A_{M} \}$ so that   $\tA_{M}^{\dagger}$ is central in $\tM$, noting that the result \cite[Lemma 3.29]{Bernstein-Zelevinsky-MR425030} used in the proof of \cite[Proposition~1.16]{Heiermann-MR2827179} generalizes easily to covering groups.

    Next we remove the assumption that $A_{M}$ is of rank $1$.
    Since $\tG$ is a $\bmu_{\bm}$-central extension of a symplectic group, an orthogonal, a special orthogonal, a unitary or a general linear group, under the assumption that the Levi subgroups decompose, we have
    \begin{align*}
      \tM &\isom \widetilde{\GL_{d_{1}}(\Fs)} \times_{\bmu_{\bm}} \cdots\times_{\bmu_{\bm}} \widetilde{\GL_{d_{r}}(\Fs)} \times_{\bmu_{\bm}} \tH, \\
      \widetilde{\sigma} &\isom \widetilde{\rho}_{1} \otimes_{\bmu_{\bm}} \widetilde{\rho}_{2} \otimes_{\bmu_{\bm}} \cdots \otimes_{\bmu_{\bm}} \widetilde{\rho}_{r} \otimes_{\bmu_{\bm}} \widetilde{\tau},\\
      \tM^{1} &\isom \widetilde{\GL_{d_{1}}(\Fs)^{1}} \times_{\bmu_{\bm}} \widetilde{\GL_{d_{2}}(\Fs)^{1}} \times_{\bmu_{\bm}} \cdots \times_{\bmu_{\bm}} \widetilde{\GL_{d_{r}}(\Fs)^{1}}\times_{\bmu_{\bm}} \tH ,\\
      \tM^{\widetilde{\sigma}} &\isom \widetilde{\GL_{d_{1}}(\Fs)^{\widetilde{\rho}_{1}}} \times_{\bmu_{\bm}} \widetilde{\GL_{d_{2}}(\Fs)^{\widetilde{\rho}_{2}}} \times_{\bmu_{\bm}} \cdots \times_{\bmu_{\bm}} {\widetilde{\GL_{d_{r}}(\Fs)^{\widetilde{\rho}_{r}}}} \times_{\bmu_{\bm}} \tH ,
    \end{align*}
    where $H$ is a possibly trivial reductive group of the same type as $G$, $\widetilde{\rho}_{i}$'s are irreducible supercuspidal unitary representations of $\widetilde{\GL_{d_{i}}(\Fs)}$ and $\widetilde{\tau}$ is an irreducible supercuspidal unitary representation of $\tH$ such that  $\bmu_{\bm}$ acts by the same character for all $\widetilde{\rho}_{i}$ and $\widetilde{\tau}$.
    By Remark~\ref{rmk:prod-cover-repns}, $\widetilde{\rho}_{1} \otimes_{\bmu_{\bm}} \widetilde{\rho}_{2} \otimes_{\bmu_{\bm}} \cdots \otimes_{\bmu_{\bm}} \widetilde{\rho}_{r} \otimes_{\bmu_{\bm}} \widetilde{\tau}$
    comes from the representation $\widetilde{\rho}_{1} \otimes \widetilde{\rho}_{2} \otimes \cdots \otimes \widetilde{\rho}_{r} \otimes \widetilde{\tau}$ of the direct product $\widetilde{\GL_{d_{1}}(\Fs)} \times \widetilde{\GL_{d_{2}}(\Fs)}  \times \cdots \times \widetilde{\GL_{d_{r}}(\Fs)} \times \tH$ on the same representation space.
    Thus we can treat each factor separately.
    The centre of each $\GL_{d_{i}}(\Fs)$-factor is of rank at most $1$ and the centre of $H$ is of rank $0$.
    The statements for the rank $0$ hold trivially and we have just treated the rank $1$ case.
  \end{proof}

  Let $X(\tM/\tM^{\widetilde{\sigma}})$ denote the set of unramified characters of $\tM$ that are trivial on $\tM^{\widetilde{\sigma}}$.
  For $\chi\in \Xnr(\tM)$, denote by $\tE_{\chi}$ the space $\tE$ equipped with the representation $\widetilde{\sigma}\otimes\chi$.

  \begin{cor}\label{cor:stab-O-defn-phi-sigma-chi}
    With the notation of Proposition~\ref{prop:M1-restriction-mult-one}, we have $\Stab \cO = X(\tM/\tM^{\widetilde{\sigma}})$ and for every $\chi\in X(\tM/\tM^{\widetilde{\sigma}})$, the map $\varphi_{\widetilde{\sigma},\chi} : \tE \rightarrow \tE_{\chi}$ which sends $e \in \widetilde{\sigma}(\tm)\tE_{1}$ to $\chi(\tm)e$ for $\tm\in \tM$ is an isomorphism intertwining $\widetilde{\sigma}$ and $\widetilde{\sigma}\otimes\chi$.
  \end{cor}
  \begin{proof}
    The proof of \cite[Corollary~1.17]{Heiermann-MR2827179} applies verbatim, once we  change $M$ to $\tM$ and $M^{\sigma}$ to $\tM^{\widetilde{\sigma}}$ in the notation there.
  \end{proof}

  \section{The Intertwining operators $A_{w}$ --- decomposed Levi  subgroups}
  \label{sec:intertw-oper-iii}

  By the results in Section \ref{sec:decom-levis}, we can always choose a base point  $\sigma\in \cO$ that is not only fixed by $W_{\cO}$ but also by $R(\cO)$.

  Recall that $B=B_{\tM}=\CC [\tM/\tM^1]$ is the ring of polynomials on the complex affine variety $\Xnr(\tM) = \Xnr(M)$.
    To $\tm\in \tM$, we associate the element $b_{\tm}\in B$ given by
  \begin{align*}
    b_{\tm} : \Xnr(\tM) &\rightarrow \CC\\
    \chi &\mapsto \chi(\tm) ,
  \end{align*}
  which is determined by the class modulo $\tM^1$ of $\tm$.
  We will also write it as  $b_{m}$ where $m\in M$  is the image of $\tm$. This is an integral ring and we write $K(B)$ for its fraction field.

  Let $\tE_{B} = \tE\otimes_{\CC} B$ and $\tE_{K(B)} = \tE\otimes_{\CC} K(B) = \tE_{B}\otimes_{B} K(B)$.
  Let $r_{B}$ be the representation of $\tM$ on $B$ given by $r_{B}(\tm) b = b b_{\tm}$.
  Denote by $\widetilde{\sigma}_{B}$ the representation $ \widetilde{\sigma}\otimes r_{B}$ of $\tM$ on $\tE_{B}$.

  As the representation $r_{B}$ can be identified by the above with $\ind_{\tM^{1}}^{\tM}(\CC)$, we find that $(\widetilde{\sigma} _B, \tE_B)$ is isomorphic to $(\ind_{\tM^{1}}^{\tM} (\widetilde{\sigma} |_{\tM^{1}}), \ind_{\tM^{1}}^{\tM} (\tE|_{\tM^{1}}))$.

  By restricting to $\tK\subseteq \tG$, we have an isomorphism between $i_{\tP}^{\tG} (\tE_{B})$ and $i_{\tK \cap \tP}^{\tK}(\tE) \otimes_{\CC} B$.
  We identify $i_{\tK \cap \tP}^{\tK}(\tE)$ with its image in $i_{\tP}^{\tG} (\tE_{B})$ via the canonical homomorphism $v\mapsto v\otimes 1$.
  We have analogous constructions for $\tE_{K(B)}$ by replacing $\tE_{B}$ by $\tE_{K(B)}$ above.

  Let $\chi\in\Xnr(\tM)$ and recall that $\tE_{\chi}$ denotes the representation space of $\widetilde{\sigma}\otimes\chi$.
  Denote by $B_{\chi}$ the ideal of $B$ which consists of elements $b\in B$ such that $b(\chi)=0$.
  Let $\spe_{\chi}$ be the specialisation map $B\rightarrow \CC, b\mapsto b(\chi)$.
  Denote also by $\spe_{\chi}$ the specialisation maps $\tE\otimes B \rightarrow \tE_{\chi}$ and $\tE\otimes K(B) \rightarrow \tE_{\chi}$ given by $e\otimes b \mapsto b(\chi) e$ where for the latter it is defined only for functions regular at $\chi$.
  By functoriality, we have also the induced maps $\spe_{\chi} : i_{\tP}^{\tG} \tE_{B} \rightarrow i_{\tP}^{\tG} \tE_{\chi}$ and $\spe_{\chi} : i_{\tP}^{\tG} \tE_{K(B)} \rightarrow i_{\tP}^{\tG} \tE_{\chi}$, which are $\tG$-equivariant.

  Let $w\in W(M)$.
Replacing $A_M$ by $\tA_M^\dagger$ in the proofs of \cite[IV.1.1]{Waldspurger-plancherel-MR1989693} (see also \cite{Heiermann-notes} and \cite{WWLi-trace-formula-II-MR3053009}), one gets, as in the reductive group case, the following description for the standard intertwining operators defined in Section~\ref{sec:intertw-oper}: there exists $b_{\tw}\in B$ and a $B$-linear intertwining map
  \begin{align*}
    J_{B,\tw} : i_{\tP}^{\tG} (\tE_{B}) \rightarrow i_{\tP}^{\tG} (\tw \tE_{B})
  \end{align*}
  such that for all $\chi\in\Xnr(\tM)$,
  \begin{equation}\label{eq:JBw-bw}
    b_{\tw}(\chi) (\lambda(\tw) J_{w^{-1}\tP | \tP}(\widetilde{\sigma}\otimes\chi) \spe_{\chi}) = \spe_{\chi} J_{B,\tw}.
  \end{equation}
  We note that $\tw \tE_{B}$ means $\tw (\tE_{B})$.
  The homomorphism
  \begin{equation*}
    J_{K(B),\tw} := b_{\tw}^{-1} J_{B,\tw}
  \end{equation*}
  is an element of $\Hom_{\tG}(i_{\tP}^{\tG} (\tE_{B}), i_{\tP}^{\tG} (\tw \tE_{K(B)}))$.

 For $w\in W(M,\cO)$, we will now define additional intertwining operators $\rho_{\tw}$ and $\tau_{\tw}$ that do not depend on  $\chi $, so that the operator $J_{K(B),\tw}$ composed with $\rho_{\tw}$ and $\tau_{\tw}$ has image in $i_{\tP}^{\tG} (\tE_{K(B)})$.

 Let $w\in W(M,\cO)$ and let $\tw\in \tK$  be a lift of $w$.
  We have the relation $\lsup{\tw}{b_{\tm}} = b_{\tw\tm \tw^{-1}}$ for $\tm\in\tM$.
  There exists an $\tM$-equivariant isomorphism $\tau_{\tw}$  defined by
  \begin{align*}
    \tau_{\tw}: B &\rightarrow B\\
     b &\mapsto  \lsup{\tw}{b},
  \end{align*}
  which intertwines $r_{B}$ and $\tw^{-1}r_{B}$.
  Via functoriality, $\tau_{\tw}$ induces an $\tM$-intertwining isomorphism between $\tw(\sigma_{B})$ and $(\tw\sigma)_{B}$ and a $\tG$-intertwining isomorphism between $i_{\tP}^{\tG} \tw(\sigma_{B})$ and $i_{\tP}^{\tG} (\tw\sigma)_{B}$ and these will also be denoted by $\tau_{\tw}$.

Now let $w\in W_{\cO}$ and let $\tw\in\tK$ be a lift of $w$.
  The operator
  \begin{equation*}
    \left( \left(\prod_{\alpha} (Y_{\alpha}(\chi) - 1) \right) \lambda(\tw) J_{w^{-1}\tP|\tP}(\widetilde{\sigma}\otimes\chi) \right)\Big|_{\chi=\triv} ,
  \end{equation*}
  where the product is over $\Sigma_{\cO,\mu}(P) \cap w^{-1}\Sigma_{\cO,\mu}(\overline{P})$,
  is well-defined and bijective by \cite[Proposition~2.4]{Heiermann-MR2827179}. Denote by $\rho_{\widetilde{\sigma},\tw}$ or $\rho_{\tw}$ its inverse and \cite[Proposition~2.4]{Heiermann-MR2827179} also gives the following.

  \begin{prop}\label{prop:rho-w-for-WO}
    Let $w\in W_{\cO}$ and $\tw\in\tK$ be a lift.
    The map $\rho_{\tw}$ is well-defined and is an isomorphism $i_{\tP}^{\tG}(\tw \tE) \rightarrow i_{\tP}^{\tG} (\tE)$.
    It comes from an isomorphism $\tw \tE\rightarrow \tE$ via functoriality, which we will also denote by $\rho_{\tw}$.
  \end{prop}

  \begin{defn}
    For each $r\in R(\cO)$, by Section \ref{sec:decom-levis}, we have $r \widetilde{\sigma} \isom \widetilde{\sigma}$.
    Fix a lift $\tilr \in \tK$ for $r$.
Take $\rho_{\tilr}: \tE \rightarrow \tE$ to be any map that intertwines $\tilr\widetilde{\sigma}$ and $\widetilde{\sigma}$.
It induces a map $\tE_{B} \rightarrow \tE_{B}$ that intertwines $\tilr\widetilde{\sigma}\otimes r_{B}$ and $\widetilde{\sigma}\otimes r_{B}$ and also a map between the induced representations, both of which we will still denote by $\rho_{\tilr}$.
 For other lifts of $r$, which are always of the form $\tilr\tm$ for some $\tm \in \tM$, define $\rho_{\tilr\tm} = \rho_{\tilr}\circ \widetilde{\sigma}(\tm) : \tE_{B}\rightarrow \tE_{B}$ that intertwines $\tilr\tm \widetilde{\sigma}\otimes r_{B}$ and $\widetilde{\sigma}\otimes r_{B}$.

With $\rho_{\tw}$ and $\rho_{\tilr}$ defined for $w \in W_{\cO}$ and $r \in R(\cO)$,  finally let $w\in W(M,\cO)$.
  Write $w=r w_{\cO}$ according to the decomposition $W(M,\cO) = R(\cO) \ltimes W_{\cO}$ given by Proposition~\ref{prop:WMO=RO-WO}, set $\rho_{\tw} = \rho_{\tilr} \rho_{\tw_{\cO}}$ where we use $\tw = \tilr \tw_{\cO}$ for the lift of $w$.
  The same symbol $\rho_{\tw}$ will also denote the isomorphisms $(\tw \tE)_{B} \rightarrow  \tE_{B}$ and $i_{\tP}^{\tG} (\tw \tE)_{B} \rightarrow i_{\tP}^{\tG} \tE_{B}$.
  We note that if we change the lift $\tw$ to $\tw\tm$ for some $\tm\in\tM$, then
  \begin{align*}
    \rho_{\tw\tm} = \rho_{\tw}\circ \widetilde{\sigma}(\tm).
  \end{align*}
\end{defn}
  
\begin{defn}
 For $w\in W(M,\cO)$ with a lift $\tw$, we obtain an element
  \begin{equation}\label{eq:defn-Aw}
    A_{\tw} := \rho_{\tw} \circ \tau_{\tw} \circ J_{K(B),\tw}
  \end{equation}
  in $\Hom_{\tG}(i_{\tP}^{\tG} \tE_{B} , i_{\tP}^{\tG} \tE_{K(B)})$ and when the lift $\tw$ lies in $\tK$, we will
write $A_{w}$ for $A_{\tw}$, as we show in Proposition~\ref{prop:Aw-Awm-change} that it does not depend on the choice of the lift in $\tK$.
  The intertwining operator $A_{w}$ extends canonically to an element in $\End_{\tG}(i_{\tP}^{\tG} \tE_{K(B)})$.
\end{defn}

\begin{rmk}\label{rmk:2-cocycle}
When $r = s_{\alpha_i,k_i}\in R(\cO)_i$, one can normalise $\rho_{\tilr}$ such that
  $\rho_{\tilr}^{2}\widetilde{\sigma}(\tilr^{-2}) = \Id$ as is done in \cite[Section~2.5]{Heiermann-MR2827179}.
However, in general, there is no reason why  $\rho_{\tilr'}\rho_{\tilr}$ would be equal to $\rho_{\tilr'\tilr}$.
One has only an equality up to a nonzero complex  constant.
As the composition is clearly associative, there is a $2$-cocycle $\eta$ on $R(\cO)$  valued in $\Bbb C^\times$ such that $\rho_{\tilr'}\rho_{\tilr}=\eta(\tilr',\tilr)\rho_{\tilr'\tilr}$ for all $r, r'\in R(\cO)$.
 \end{rmk}

  For $b\in K(B)$ and $\chi\in\Xnr(\tM)$, let $b_{\chi}\in K(B)$ be the right translation of $b$ by $\chi$, i.e.,
  \begin{equation*}
    b_{\chi}(\cdot) = b(\cdot\ \chi).
  \end{equation*}
  Then $(b_{\tm})_{\chi}  = \chi(\tm) b_{\tm}$.

  \begin{prop}\label{prop:Aw-Awm-change}
    Let $w\in W(M,\cO)$ such that $w=r w_{\cO}$ according to the decomposition $W(M,\cO) = R(\cO) \ltimes W_{\cO}$.
Then $A_{\tw}$ does not depend on the choice of lift of $w$ in $\tK$.
  \end{prop}

  \begin{rmk}
    We can define $A_{\tw \tm}$  with $\tm \in \tM$.
    Then we have the relation
    \begin{align*}
      A_{\tw \tm} =\ \lsup{w}{b_{\tm}} A_{\tw} .
    \end{align*}
  \end{rmk}

  \begin{proof}
  The proof of \cite[Proposition~3.1]{Heiermann-MR2827179} generalises.
  \end{proof}

  To state  some properties of $A_{w}$,     we need some more auxiliary intertwining operators.
  By unravelling the definitions,  we have the following
  \begin{lemma}\label{lemma:rho_chi}    
  Let $\chi\in\Xnr(\tM)$.
    Then the map
    \begin{align*}
      \rho_{\chi}:  B &\rightarrow  B\\
       b &\mapsto  b_{\chi}
    \end{align*}
    is an isomorphism between  the representations $r_{B}$ and $r_{B}\otimes\chi$ of $\tM$.
    In addition, for $w\in W^{G}(M)$ and a lift $\tw\in\tK$, we have $\rho_{\chi}^{-1} = \rho_{\chi^{-1}}$ and $\tau_{\tw}\circ \rho_{\chi} = \rho_{{\tw}{\chi}}\circ \tau_{\tw}$.
    The map $\rho_{\chi}$ induces an $\tM$-equivariant isomorphism, also denoted by $\rho_{\chi}$
    \begin{align*}
      \rho_{\chi}: \tE\otimes_{\CC} B &\rightarrow \tE_{\chi}\otimes_{\CC} B\\
      e\otimes b &\mapsto e\otimes b_{\chi}.
    \end{align*}
  \end{lemma}
\begin{defn}\label{defn:phi_chi}
    Let $\chi\in\Stab(\cO)$.
    Recall the intertwining map $\varphi_{\widetilde{\sigma},\chi} : \tE \rightarrow \tE_{\chi}$ defined in Corollary~\ref{cor:stab-O-defn-phi-sigma-chi}.
    Define the automorphism
    \begin{equation}
      \label{eq:phi_chi}
      \begin{split}
        \phi_{\chi} : \tE_{B} &\rightarrow \tE_{B}\\
        e\otimes b &\mapsto (\varphi_{\widetilde{\sigma},\chi}e) \otimes b_{\chi^{-1}}.
      \end{split}
    \end{equation}
    Denote also by $\phi_{\chi}$, the automorphism of $i_{\tP}^{\tG}(\tE_{B})$ (resp. $i_{\tP}^{\tG}(\tE_{K(B)})$) induced by functoriality.
  \end{defn}

  We now record some properties of $A_{w}$, but   we omit the proofs,  since they are similar to \cite[Section~3]{Heiermann-MR2827179}. (The identity (4) in Proposition~\ref{prop:relations-Aw-Ar-B} is weaker than in \cite{Heiermann-MR2827179}  as one has here, in general, a $2$-cocycle because of Remark~\ref{rmk:2-cocycle}.)
  From now on, lifts of  elements in $W(M,\cO)$ are required to lie in $\tK$.

  \begin{lemma}\label{lemma:sp-Aw}
    Let $w\in W(M,\cO)$, $v\in i_{\tP}^{\tG} \tE_{B}$ and $\chi\in \Xnr(\tM)$ such that $\widetilde{\sigma}\otimes w^{-1}\chi$ is a regular point for $J_{\tw^{-1}\tP|\tP}$.
    Then
    \begin{align*}
      \spe_{\chi} A_{w} v = \rho_{\tw} \lambda(\tw) J_{w^{-1}\tP|\tP}(\widetilde{\sigma}\otimes w^{-1}\chi) \spe_{w^{-1}\chi} v.
    \end{align*}
  \end{lemma}

  \begin{prop}\label{prop:relations-Aw-Ar-B}
    \begin{enumerate}
    \item For $w,w'\in W_{\cO}$ such that $\ell_{\cO}(ww') = \ell_{\cO}(w) + \ell_{\cO}(w')$, we have
    \begin{align}
    A_{w}A_{w'} = A_{ww'}.
    \end{align}
    \item If $s=s_{\alpha}$ is a simple reflection in $W_{\cO}$, then there exists a complex number $c_{s}''\neq 0$ such that
    \begin{align} \label{eq:cs-double-prime}
      A_{s}^{2} = c_{s}'' (\mu^{\tM_{\alpha}})^{-1}.
    \end{align}
    \item \label{item:rel-Ar-Aw} For $w\in W_{\cO}$ and $r\in R(\cO)$, we have
      \begin{align*}
        A_{r} A_{w} &= A_{rw},\\
        A_{w} A_{r} &= A_{r} A_{r^{-1}wr}.
      \end{align*}
    \item \label{item:rel-Ar-Ar}  For $r,r'\in R(\cO)$, we have
      \begin{align*}
        A_{r} A_{r'} = \eta(r,r') A_{rr'} ,
      \end{align*}
      for some  $2$-cocycle $\eta$ on $R(\cO)$ valued in $\CC^{\times}$.
    \item \label{item:rel-Aw-b} For $b\in B$ and $w\in W(M,\cO)$, we have
      \begin{align*}
        A_{w} b =\ \lsup{w}{b} A_{w}.
      \end{align*}
    \end{enumerate}
  \end{prop}

We can now describe certain spaces of intertwining operators explicitly.

  \begin{prop}\label{prop:phi_chi-basis}    As $K(B)$-vector spaces,
    \begin{align*}
      \Hom_{\tM}(\tE_{B}, \tE_{K(B)}) \isom \bigoplus_{\chi \in \Stab(\cO)} K(B) \phi_{\chi},
    \end{align*}
    where $\phi_{\chi}$ is defined in \eqref{eq:phi_chi}.
  \end{prop}

  \begin{proof}
    This is a generalisation of Proposition~3.6 of \cite{Heiermann-MR2827179} to the case of finite central extensions when the multiplicity one result of Proposition~\ref{prop:M1-restriction-mult-one} holds.
    Using $\tE_{B} \isom \ind_{\tM^{1}}^{\tM} \tE|_{\tM^{1}}$, Frobenius reciprocity for compact induction \cite[Proposition~2.29]{Bernstein-Zelevinsky-MR425030} and Proposition~\ref{prop:M1-restriction-mult-one}, we have

    \begin{equation}\label{eq:isom-beta}
      \begin{aligned}
        \Hom_{\tM} (\tE_B, \tE_{K(B)})&=\Hom_{\tM} (\ind_{\tM^{1}}^{\tM} \tE|_{\tM^{1}}, \tE_{K(B)})\\
        =\Hom_{\tM^{1}} ( \tE|_{\tM^{1}}, \tE_{K(B)}|_{\tM^{1}}) &= \bigoplus_{\tm \in \tM/\tM^{\widetilde{\sigma}}} \Hom_{\tM^{1}} (  \widetilde{\sigma}(\tm) \tE_{1}, \tE_{K(B)}|_{\tM^{1}}).
      \end{aligned}
    \end{equation}
    As the components  are pairwise non-isomorphic, we find that  $\Hom_{\tM^{1}} ( \tE|_{\tM^{1}}, \tE_{K(B)}|_{\tM^{1}}) \isom K(B)^{\tM / \tM^{\widetilde{\sigma}}}$.

    Write $\beta$ for the second isomorphism in \eqref{eq:isom-beta}.
    For $e\in \tE$, letting $v_{e}$ denote the element in $\ind_{\tM^{1}}^{\tM} \tE|_{\tM^{1}}$ which has value $e$ at $1$ and vanishes outside $\tM^{1}$, we find that $\beta $ sends $\varphi\in \Hom_{\tM} (\ind_{\tM^{1}}^{\tM} \tE|_{\tM^{1}}, \tE_{K(B)})$ to $e\mapsto \varphi(v_{e})$.

It remains to check that the $\beta(\phi_{\chi})$'s for $\chi\in\Stab(\cO)$ form a basis of $\Hom_{\tM^{1}} ( \tE|_{\tM^{1}}, \tE_{K(B)}|_{\tM^{1}})$.
    By unravelling the definitions, for $e\in \tE$,
    \begin{align*}
      \beta(\phi_{\chi})(e) = \phi_{\chi} (v_{e}) = \varphi_{\widetilde{\sigma},\chi} (e) \otimes 1.
    \end{align*}
    Thus on the component $\widetilde{\sigma}(\tm) \tE_{1}$ of $\tE$, $\beta(\phi_{\chi})$ acts by multiplication by $\chi(\tm)$.
    By the duality between $\tM/\tM^{\widetilde{\sigma}}$ and $\Stab(\cO)$ and linear independence of characters, we see that the  $\beta(\phi_{\chi})$'s form a basis of $\Hom_{\tM^{1}} ( \tE|_{\tM^{1}}, \tE_{K(B)}|_{\tM^{1}})$.
  \end{proof}

  \begin{thm}\label{thm:Hom-i-EB-gen-by-phiA}
    The intertwining operators $\phi_{\chi} A_{w}$, for $\chi\in\Stab(\cO)$ and $w\in W(M,\cO)$ form a basis of $\Hom_{\tG}(i_{\tP}^{\tG}(\tE_{B}), i_{\tP}^{\tG}(\tE_{K(B)}))$ as a $K(B)$-vector space,
    \begin{align*}
      \Hom_{\tG}(i_{\tP}^{\tG}(\tE_{B}), i_{\tP}^{\tG}(\tE_{K(B)})) = \bigoplus_{\substack{w\in W(M,\cO)\\ \chi \in \Stab(\cO)}} K(B) \phi_{\chi}A_{w}.
    \end{align*}
  \end{thm}

  \begin{proof}
    For the linear independency the proof of Proposition~3.7 of \cite{Heiermann-MR2827179} carries over.
    It relies crucially on Proposition~3.6 of \cite{Heiermann-MR2827179}, which we have generalised in Proposition~\ref{prop:phi_chi-basis}.

    Then we follow the arguments in Theorem~3.8 of \cite{Heiermann-MR2827179} to show that the operators span $\Hom_{\tG}(i_{\tP}^{\tG}(\tE_{B}), i_{\tP}^{\tG}(\tE_{K(B)}))$.
    The proof of Theorem~3.8 of \cite{Heiermann-MR2827179} makes use of the geometric lemma, which is valid for finite central extensions by the very general \cite[Theorem~5.2]{BZ1977-MR579172}.
        It shows that the (non-trivial) subquotients are $\Hom_{\tM}(\tw \tE_{B}, \tE_{K(B)})$ for $w\in W(M,\cO)$, each of which, by Proposition~\ref{prop:phi_chi-basis}, is of dimension $|\Stab(\cO)|$.
    Comparing the dimensions, we find that the intertwining operators $\phi_{\chi} A_{w}$ must span $\Hom_{\tG}(i_{\tP}^{\tG}(\tE_{B}), i_{\tP}^{\tG}(\tE_{K(B)}))$.
  \end{proof}

  \begin{lemma}[{\cite[Lemme 2.7]{Heiermann-MR2827179}}]\label{lemma:braid-b}
    Let $\alpha,\alpha'\in \Delta_{\cO}$.
    Set $s=s_{\alpha}$ and $s'=s_{\alpha'}$.
    Let $m \in M_{\alpha}^{1} \cap M$ and $m' \in M_{\alpha'}^{1} \cap M$.
    Then
    \begin{align*}
      b_{m} \lsup{s}{b_{m'}} \lsup{ss'}{b_{m}} \lsup{ss's}{b_{m'}} \cdots =
      b_{m'} \lsup{s'}{b_{m}} \lsup{s's}{b_{m'}} \lsup{s's s'}{b_{m}} \cdots
    \end{align*}
    where the number of factors of each side is equal to the order of $ss'$.
  \end{lemma}
  
  \begin{rmk}\label{rmk:braid-bnoncmt}
  One observes that the statement of the lemma is only about $M/M^1$ and a system of symmetries coming from the Weyl group of a root system related to $M/M^1$. 
  \end{rmk}

  \section{The  algebra of intertwining operators over the rational function field --- decomposed Levi subgroups}
  \label{sec:preliminaries-intertwining-algebra}

   The aim of this section is to treat $\ind_{\tM^{1}}^{\tM} \tE_{1}$ and to write down analogues of the results from  Section~\ref{sec:intertw-oper-iii} for $\tE_B = \tE\otimes_{\CC} \ind_{\tM^{1}}^{\tM} (\CC) \isom \ind_{\tM^{1}}^{\tM} \tE\vert_{\tM^{1}}$.
  The operators $\phi _{\chi }$ will disappear when describing the analogous intertwining algebra.

  As in Proposition~\ref{prop:M1-restriction-mult-one}, we fix an irreducible component $(\widetilde{\sigma}_{1},\tE_{1})$ of $\widetilde{\sigma}|_{\tM^{1}}$.
  Recall that $B$ is the ring of regular functions on $\Xnr(\tM) = X(\tM/\tM^{1})$.
  Let $\BO$ be the subset of $B$ consisting of functions that are invariant under translation by elements of $\Stab(\cO) = X(\tM/\tM^{\widetilde{\sigma}})$ (Corollary~\ref{cor:stab-O-defn-phi-sigma-chi}).
  In other words, it is the ring of regular functions on the quotient affine variety $X(\tM/\tM^{1})/\Stab(\cO) \isom X(\tM/\tM^{1})/X(\tM/\tM^{\widetilde{\sigma}})$ and we  identify it with the group algebra $\CC[\tM^{\widetilde{\sigma}}/\tM^{1}]$.
  Let $\KBO$ be the fraction field of $\BO$.
  We note that $\BO$ is a factorial ring, since $\tM^{\widetilde{\sigma}}/\tM^{1}$ is a free $\ZZ$-module of the same rank as $\tM/\tM^{1}$, which can be seen as follows.
  It is clear that  $\tA_{M}^{\dagger}\tM^{1}\subseteq \tM^{\widetilde{\sigma}} \subseteq \tM$.
  We have a surjective homomorphism $\tM \rightarrow \ZZ^{r}$ whose kernel is $\tM^{1}$, where $r$ is the rank of $A_{M}$ and  the image of the restriction $\tA_{M}^{\dagger}\tM^{1} \rightarrow \ZZ^{r}$ is of finite index in $\ZZ^{r}$.

Write $\cR(\tM/\tM^{1})$ for a set of representatives of $\tM/\tM^{1}$ and $\cR(\tM/\tM^{\widetilde{\sigma}})$ for that of $\tM/\tM^{\widetilde{\sigma}}$.

Lemmes 4.1 and 4.2 of \cite{Heiermann-MR2827179} readily generalise to the present case.
  \begin{lemma}\label{lemma:EBO-identification}
    The image of $\ind_{\tM^{1}}^{\tM} \tE_{1}$ under the isomorphism $\ind_{\tM^{1}}^{\tM} \tE|_{\tM^{1}} \rightarrow \tE\otimes_{\CC} B$ is given by
    \begin{align}\label{eq:ind-E1-image}
      \sum_{\tm \in \cR(\tM/\tM^{1})} \widetilde{\sigma}(\tm)\tE_{1} \otimes b_{\tm}
    \end{align}
    and it is a  $\BO$-submodule of $\tE\otimes_{\CC} B$.
  \end{lemma}

  Write $\tEBO$ for the $\BO$-module $\ind_{\tM^{1}}^{\tM} \tE_{1}$ which is identified with \eqref{eq:ind-E1-image}.
  The notation does not imply that $\ind_{\tM^{1}}^{\tM} \tE_{1}$ is isomorphic to $\tE\otimes_{\CC} \BO$.
  Set $\tEKBO = \tEBO \otimes_{\BO} \KBO$.

  \begin{lemma}\label{lemma:hom-EBO}
    We have
    \begin{align*}
      \Hom_{\tM}(\tEBO, \tEKBO) \isom \KBO.
    \end{align*}
  \end{lemma}

  \begin{lemma}\label{lemma:bAw-preserves-iEBO}
    For every $w\in W(M,\cO)$, there exists $\tm_{w}\in \tM$ such that $b_{\tm_{w}} A_{w}$ leaves the  space $(i_{\tP}^{\tG} \tEBO)_{\KBO}$ invariant.
    
    Furthermore, if $w=s_{\alpha}$ for $\alpha\in\Delta_{\cO,\mu}$, then $\tm_{w}$ can be chosen to lie in $\tM \cap \tM_{\alpha}^{1}$; the same is true for $\tm_{r}$
     if $r =s_{\alpha} \in R(\cO)$, for some $\alpha\in\Sigma(A_M)$. If $r$ is a product of $s_{i,k_i}$ in the context of  Prop~\ref{prop:Sigma-O-mu-cover-decomposed-Levis} (4) (b) or (c), then the image $m_{r}$ of $\tm_{r}$ in $M$ can be written as a product of elements in the various  $M \cap M_{\alpha_{i,k_i}}^{1}$\footnote{In the case of Prop~\ref{prop:Sigma-O-mu-cover-decomposed-Levis} (4) (c) $\alpha_{i,k_i}$ has to be interpreted as a root for $A_M$ relative to a split odd orthogonal group containing $G$. By  Remark~\ref{rmk:braid-bnoncmt} this will allow $b_{\tm_r}$ to enter in the context of Lemma~\ref{lemma:braid-b}.}.
    \end{lemma}
    
    \begin{proof}
      This proof generalizes readily from \cite[Lemme~4.5]{Heiermann-MR2827179}. We will give additional details for $r\in R(\cO )$. 
      
      If, in the notations of Prop~\ref{prop:Sigma-O-mu-cover-decomposed-Levis}, $s=s_{i,j}^\pm$ with $j=1,\dots ,k_i-1$, then it is quite easy to give directly an element $\tm_{r}\in\tM \cap \tM_{\alpha}^{1}$ such that $r\widetilde{\sigma}_1\simeq \tm_{r}\widetilde{\sigma}_1$ and the remaining part of the proof of \cite[Lemme~4.5]{Heiermann-MR2827179} remains valid. 
      
      If $r=s_{i,k_i}$ and $\alpha_{i,k_i}$ is a root for $A_M$, then \cite[Lemme~4.5]{Heiermann-MR2827179} generalizes directly.
      
      It remains the case of $r$ as in Prop~\ref{prop:Sigma-O-mu-cover-decomposed-Levis} (4) (b) and (c). In the case of Prop~\ref{prop:Sigma-O-mu-cover-decomposed-Levis} (4) (b), each $\alpha_{i,k_i}$ is still a root for $A_M$ and we can choose $\tm_{i,k_i}$ in $\tM\cap \tM_{\alpha_{i,k_i}}$ for each $i$ and it is immediate that the product of the $\tm_{i,k_i}$ gives $\tm_r$. In the situation Prop~\ref{prop:Sigma-O-mu-cover-decomposed-Levis} (4) (c) $\alpha_{i,k_i}$ is no more a root for $A_M$. However, one can consider $M$ as a Levi subgroup of some odd special orthogonal group $G'$ and define $\alpha_{i,k_i}$ relative to $G'$. Then $\tm_r$ will be a product of elements $\tm_{i,k_i}$ with image in $M$ with non trivial entry only in the last copy of the corresponding $\GL_{d_i}$. Consequently the image in $M$ of each $\tm_{i,k_i}$ lies in $M\cap M_{\alpha_{i,k_i}}^1$ with $M_{\alpha_{i,k_i}}^1$ defined relative to $G'$.

\end{proof}

  \begin{defn}\label{defn:Jw}
    \begin{enumerate}
    \item For   a simple reflection  $s = s_{\alpha}\in W_{\cO}$, fix $\tm_{s}\in \tM\cap\tM_{\alpha}^{1}$ as in Lemma~\ref{lemma:bAw-preserves-iEBO} and set $J_{s} = b_{\tm_{s}}A_{s}$.
    \item     For $w\in W_{\cO}$, write $w=s_{1}\cdots s_{\ell}$ as a reduced product of simple reflections  and set $J_{w} = J_{s_{1}} \cdots J_{s_{\ell}}$.   This is well-defined by Lemma~\ref{lemma:Jw-properties} below.
    \item     %
      For $r\in R(\cO)$ which is either $s_{\alpha}$ for some $\alpha\in\Delta(A_M)$ or a product of $s_{i,k_i}$ in the context of  Proposition~\ref{prop:Sigma-O-mu-cover-decomposed-Levis} (4) (b) or (c), fix $\tm_{r}\in \tM$ as in Lemma~\ref{lemma:bAw-preserves-iEBO} and set $J_{r} = b_{\tm_{r}}A_{r}$.
    \item For general $r\in R(\cO)$, fix a decomposition of $r = r_{1} \cdots r_{t}$ for some $t\in\ZZ_{\ge 0}$ with  $r_{i} = s_{\alpha}$ or a product of $s_{i,k_{i}}$ which are in (3) and set $J_{r} = J_{r_{1}} \cdots J_{r_{t}}$.
    \item     For $w\in W(M,\cO)$, write $w = r w_{\cO}$ with $r\in R(\cO)$ and $w_{\cO} \in W_{\cO}$ and set $J_{w} = J_{r}J_{w_{\cO}}$.
    \end{enumerate}
  \end{defn}

  The following is the analogue of \cite[Lemme~4.7]{Heiermann-MR2827179}  except for (4) which was forgotten in \cite{Heiermann-MR2827179} . 
  \begin{lemma}\label{lemma:Jw-properties}
    \begin{enumerate}
    \item \label{item:1}Let $w\in W_{\cO}$.
      Write $w=s_{1}\cdots s_{\ell}$ as a reduced product of simple reflections.
      Then the operator $J_{s_{1}} \cdots J_{s_{\ell}}$ does not depend on the choice of the decomposition of $w$.
    \item \label{item:2}If $s=s_{\alpha}$ is a simple reflection in $W_{\cO}$, then $J_{s}^{2} = c_{s}'' (\mu^{\tM_{\alpha}})^{-1}$ where $c_{s}''\in\CC^{\times}$ is defined in \eqref{eq:cs-double-prime}.
    \item \label{item:3}For every $r,r'\in R(\cO)$, we have $J_{r}J_{r'}= \eta(r,r') J_{rr'}$.
      In particular, if $r^{2} = 1$, $J_{r}^{2} = \eta(r,r)  \Id$. 
    \item \label{item:4}For all $r\in R(\cO)$, $w\in W_{\cO}$, we have $J_wJ_r= J_rJ_{r^{-1}wr}$. 
    \item \label{item:5}For $b\in \KBO$ and $w\in W(M,\cO)$, $J_{w} b =\ \lsup{w}{b} J_{w}$.
    \end{enumerate}
  \end{lemma}
  
  \begin{proof} The proofs of (1), (2) and (5) are analogous to  \cite[Lemme~4.7]{Heiermann-MR2827179}.
    
    For part (3), one has to show that the product of the $b_{m_\bullet}$ matches. As in \cite{Heiermann-MR2827179}, this follows from Lemma~\ref{lemma:braid-b} and the definitions. In the situation of Proposition~\ref{prop:Sigma-O-mu-cover-decomposed-Levis} (4) (b) or (c), where the elements in $R(\cO)^{I'}$ are products of reflections $s_{i,k_i}$, one has to use in addition Remark~\ref{rmk:braid-bnoncmt} with Lemma~\ref{lemma:bAw-preserves-iEBO}.

By Proposition~\ref{prop:relations-Aw-Ar-B} (3),  part (4)  is reduced to the equality $b_{m_r} {^rb_{m_w}}=b_{m_{rwr^{-1}}} {^{rwr^{-1}}b_{m_r}}$. To check this equality one can reduce to a $B_2$ situation by the nature of the elements of $R(\cO)$ given in Proposition~\ref{prop:Sigma-O-mu-cover-decomposed-Levis} and arguments similar to those used for the proof of part (3). One has to work with representatives of the elements $m_\bullet$ as given in \cite[3.3]{Heiermann-MR2643577}. The equality is then an elementary straightforward computation.
\end{proof}

  Next we proceed to prove a $\KBO$-analogue of Theorem~\ref{thm:Hom-i-EB-gen-by-phiA}.

  \begin{thm}\label{thm:Hom-i-EBO-gen-by-J}
    As $\KBO$-vector spaces,
    \begin{align*}
      \Hom_{\tG}(i_{\tP}^{\tG}(\tEBO), i_{\tP}^{\tG}(\tEKBO)) = \bigoplus_{w\in W(M,\cO)} \KBO J_{w}.
    \end{align*}
  \end{thm}

  \begin{proof}
    The proof follows a similar line to Theorem~\ref{thm:Hom-i-EB-gen-by-phiA}, as is done in the proof of Th\'eor\`eme~4.9 of \cite{Heiermann-MR2827179}.
    The intertwining operators $J_{w}$ for  ${w\in W(M,\cO)}$ are $\KBO$-linearly independent.
    Using the geometric lemma \cite[Theorem~5.2]{BZ1977-MR579172} which is valid  for  finite central extensions, we find that the  (non-trivial) subquotients are $\Hom_{\tM}(w \tEBO, \tEKBO)$ for $w\in W(M,\cO)$, which by Lemma~\ref{lemma:hom-EBO} is of dimension $1$.
    Comparing the dimensions, we find that the intertwining operators  $J_{w}$ for  ${w\in W(M,\cO)}$ must generate $\Hom_{\tG}(i_{\tP}^{\tG}(\tEBO), i_{\tP}^{\tG}(\tEKBO))$.
     \end{proof}

  \section{  The  algebra of intertwining operators over the ring of regular functions --- decomposed Levi  subgroups}
  \label{sec:intertwining-algebra}

  We continue to assume that the standard Levi subgroups decompose.
  We now define operators $T_{w}$ for $w\in W_{\cO}$ and prove in Theorem~\ref{thm:Jr-Tw-basis} an integral version  of Theorem~\ref{thm:Hom-i-EBO-gen-by-J}, namely,  that the intertwining operators $J_{r}T_{w}$  for $r\in R(\cO)$ and $w\in W_{\cO}$ form a basis of the free $B_{\cO}$-module $\End_{\tG}(i_{\tP}^{\tG}\tEBO)$.

  Let $s=s_{\alpha}$ for $\alpha\in\Delta_{\cO}$.
  By Lemma~\ref{lemma:Jw-properties} and Proposition~\ref{prop:expression-mu-function}, $J_{s}$ can possibly have poles only at $\chi\in\Xnr(\tM)$ such that $X_{\alpha}(\chi)=\pm 1$.
  We analyse these potential poles and  produce a regular operator $T_{s}$ by adjusting $J_{s}$.
  Similar to \cite[Lemma~5.1]{Heiermann-MR2827179}, we have,
  \begin{lemma}\label{lemma:Js-poles}
    Let $\alpha\in\Delta_{\cO}$.
    Set $s=s_{\alpha}$ and $c_{s} = c_{s}'^{-1} c_{s}''$ where $c_{s}'$ and $c_{s}''$ are defined in Proposition~\ref{prop:expression-mu-function} and Proposition~\ref{prop:relations-Aw-Ar-B} respectively.
    Then the following hold.
    \begin{enumerate}
    \item The operator $\spe_{\triv} (X_{\alpha} - 1) J_{s}$ is scalar on $i_{\tP\cap\tK}^{\tK} \tE$ and
      \begin{align*}
        (\spe_{\triv} (X_{\alpha} - 1) J_{s})^{2} = \frac{1}{4} c_{s} (1-q_{\Fz}^{-a_{s}})^{2} (1+q_{\Fz}^{-a_{s,-}})^{2}.
      \end{align*}
      Moreover, if $\chi\in\Xnr(\tM)$ satisfies $X_{\alpha}(\chi) = 1$, %
      then for $v\in i_{\tP\cap\tK}^{\tK} \tEBO$, we have
      \begin{align*}
        \spe_{\chi} (X_{\alpha} - 1) J_{s} v = (\spe_{\triv} (X_{\alpha} - 1) J_{s} )\spe_{\chi} v .
      \end{align*}

    \item Assume that $\mu^{\tM_{\alpha}}$ vanishes when $X_{\alpha} = -1$.
      Fix $\chi_{-1} \in \Xnr(\tM)$ such that $X_{\alpha}(\chi_{-1}) = -1$.
      Then the operator $\spe_{\chi_{-1}} (X_{\alpha} + 1) J_{s}$ is scalar on $i_{\tP\cap\tK}^{\tK} \tE$ and
      \begin{align*}
        (\spe_{\chi_{-1}} (X_{\alpha} + 1) J_{s})^{2} = \frac{1}{4} c_{s} (1+q_{\Fz}^{-a_{s}})^{2} (1-q_{\Fz}^{-a_{s,-}})^{2}.
      \end{align*}
            Moreover, if $\chi\in\Xnr(\tM)$ satisfies $X_{\alpha}(\chi) = -1$,
      then for $v\in i_{\tP\cap\tK}^{\tK} \tEBO$, we have
      \begin{align*}
        \spe_{\chi} (X_{\alpha} + 1) J_{s} v = (\spe_{\chi_{-1}} (X_{\alpha} + 1) J_{s}) \spe_{\chi} v.
      \end{align*}
    \end{enumerate}
  \end{lemma}

  \begin{defn}\label{defn:T-from-J}
    Let $\alpha\in\Delta_{\cO}$ and $s=s_{\alpha}$.
    Fix a square root $c_{s}^{\half}$ of $c_{s}$.
    Fix $\chi_{-1} \in \Xnr(\tM)$ such that $X_{\alpha}(\chi_{-1}) = -1$.
    By Lemma~\ref{lemma:Js-poles}, there exists $\epsilon_{1}, \epsilon_{-1} \in \{ \pm 1  \}$ such that
    \begin{align*}
      \spe_{\triv} (X_{\alpha}-1) J_{s} &= \half \epsilon_{1} c_{s}^{\half} (1-q_{\Fz}^{-a_{s}})(1+q_{\Fz}^{-a_{s,-}}) \spe_{\triv}\\
      \spe_{\chi_{-1}} (X_{\alpha}+1) J_{s} &= \half \epsilon_{-1} c_{s}^{\half} (1+q_{\Fz}^{-a_{s}})(1-q_{\Fz}^{-a_{s,-}}) \spe_{\chi_{-1}}.
    \end{align*}
  \end{defn}

  Set
  \begin{align*}
    R_{s} &=
    \begin{cases}
      - \epsilon_{1} q_{\Fz}^{a_{s}+a_{s,-}} c_{s}^{-\half} J_{s},  &\quad\text{if $\epsilon_{1}a_{s,-} \neq \epsilon_{-1}a_{s,-}$,} \\
      -X_{\alpha}\epsilon_{1} q_{\Fz}^{a_{s}+a_{s,-}} c_{s}^{-\half} J_{s},  &\quad\text{if $\epsilon_{1}a_{s,-} = \epsilon_{-1}a_{s,-}$.}
    \end{cases}
  \end{align*}
  Finally let
  \begin{align*}
    T_{s} = R_{s} + (q_{\Fz}^{a_{s} + a_{s,-}}) \frac{X_{\alpha} (X_{\alpha} - \frac{q_{\Fz}^{a_{s,-}} - q_{\Fz}^{a_{s}}}{q_{\Fz}^{a_{s}+a_{s,-}}-1})}{X_{\alpha}^{2}-1}.
  \end{align*}

  Analogous to \cite[Proposition~5.4]{Heiermann-MR2827179}, we have,
  \begin{prop}\label{prop:T-integral}
    The operators $T_{s_{\alpha}}$ for $\alpha\in\Delta_{\cO}$ belong to $\End_{\tG}(i_{\tP}^{\tG}(\tEBO))$.
  \end{prop}

Some basic computation as in \cite[Proposition~5.5]{Heiermann-MR2827179} shows the following.
  \begin{prop}\label{prop:hecke-relation}
    Let $s=s_{\alpha}$ for $\alpha\in\Delta_{\cO}$.
    Then $(T_{s}+1)(T_{s}-q_{\Fz}^{a_{s}+a_{s,-}}) = 0$.
  \end{prop}

  \begin{defn}
    Let $w\in W_{\cO}$.
    Fix a  decomposition  of  $w$ into a reduced product of simple reflections: $w= s_{1} \cdots s_{\ell}$.
    Define $T_{w} = T_{s_{1}} \cdots T_{s_{\ell}}$.
    In particular, $T_{1} = \Id$.
  \end{defn}
  After more setup, we will eventually show (Lemma~\ref{lemma:Ts-braid}) that $T_{w}$ does not depend on the choice of the reduced decomposition of $w$.
  The results before Lemma~\ref{lemma:Ts-braid} do not depend on this fact.

\begin{thm}\label{thm:Jr-Tw-basis}
  We have
  \begin{align*}
    \End_{\tG}(i_{\tP}^{\tG}(\tEBO)) = \bigoplus_{\substack{r\in R(\cO)\\w\in W_{\cO}}} \BO J_{r}T_{w}.
  \end{align*}
\end{thm}

\begin{proof}
  The arguments in the proof of Th\'eor\`eme~5.10 of \cite{Heiermann-MR2827179} generalise to the present case.
  One makes use of  Theorem~\ref{thm:Hom-i-EBO-gen-by-J} to show that  $J_{r}T_{w}$ form a $\KBO$-basis of $\Hom_{\tG}(i_{\tP}^{\tG}(\tEBO), i_{\tP}^{\tG}(\tEKBO))$.
  Then highly non-trivial arguments  in the line of \cite[5.6--5.9]{Heiermann-MR2827179} are needed to prove that these operators generate the intertwining algebra over $\BO$.
\end{proof}

We have analogues of Corollaries~5.11 and 5.12 of \cite{Heiermann-MR2827179}.
\begin{cor}
  The centre $\cZ_{\cO}$ of $\End_{\tG}(i_{\tP}^{\tG} \tEBO)$ is exactly the $W(M,\cO)$-invariant elements in $\BO$.
\end{cor}

Denote the fraction field of $\cZ_{\cO}$ by $K(\cZ_{\cO})$.
\begin{cor}
  The algebra $\End_{\tG}(i_{\tP}^{\tG} \tEBO)\otimes_{\cZ_{\cO}} K(\cZ_{\cO})$ is a $\KBO$-module which is canonically isomorphic to $\Hom_{\tG}(i_{\tP}^{\tG} \tEBO, i_{\tP}^{\tG} \tEKBO)$.
\end{cor}

\section{Hecke algebras with parameters --- decomposed Levi subgroups}
\label{sec:hecke-algebras}

We will now relate the algebra of intertwining operators $\End_{\tG}(i_{\tP}^{\tG} \tEBO)$ to certain affine Hecke algebras. The notion of a Hecke algebra with parameters was defined in \cite{Lusztig-MR991016}.

\begin{defn} Let $(\Lambda, \Lambda^{\vee}, \Sigma, \Sigma^{\vee},\Delta)$ be a based root datum, where $\Lambda$ and $\Lambda^{\vee}$ are free abelian groups of finite type in duality with respect to  a perfect $\ZZ$-bilinear map $\Lambda\times\Lambda^{\vee}\rightarrow\ZZ$, $\Sigma\subseteq\Lambda$ is a root system, $\Sigma^{\vee}\subseteq\Lambda^{\vee}$ is its dual root system and $\Delta$ is a base for $\Sigma$.
Denote by $W(\Sigma)$ the Weyl group of $\Sigma$.
Let $\Sigma_{i}$ be an irreducible component of $\Sigma$ and let $\Delta_{i}=\Delta\cap\Sigma_{i}$.
For each $\alpha\in\Delta_{i}$, take a real number $q_{\alpha}>1$ such that $q_{\alpha} = q_{\beta}$ if $\alpha$ and $\beta$ are conjugate by $W(\Sigma_{i})$.
In addition, if $\Sigma_{i}$ is of type $B$, take an extra real number $q_{i} \geq 1$.

Let $B_{0}(\Sigma)$ be the $\CC$-algebra with generators given by the symbols $U_{s_{\alpha}}$ for $\alpha\in\Delta$ which satisfy the braid relations
\begin{align*}
  U_{s_{\alpha}}U_{s_{\beta}}U_{s_{\alpha}} \cdots =   U_{s_{\beta}}U_{s_{\alpha}}U_{s_{\beta}} \cdots
\end{align*}
where each side has $m(\alpha,\beta)$ factors where $m(\alpha,\beta)$ is the order of the element $s_{\alpha}s_{\beta}$ in $W(\Sigma)$.
We define $U_{w} = U_{s_{1}}\cdots U_{s_{r}}$ where $w=s_{1}\cdots s_{r}$ is a reduced decomposition of $w$ into simple reflections.
This is well-defined due to the braid relation.
Define $\cH_{0}(\Sigma, \{ q_{\alpha} \}_{\alpha\in\Delta})$ to be the quotient of $B_{0}(\Sigma)$ by the ideal generated by $(U_{s_{\alpha}}+1)(U_{s_{\alpha}}-q_{\alpha})$.

Let  $C$ denote the group algebra  $\CC[\Lambda]$ and we denote the element in $C$ that is associated  to $\lambda\in\Lambda$ by $Z_{\lambda}$.
Let $K(C)$ denote its field of fractions.
Define $\cH(\Sigma, \{ q_{\alpha} \}, \{ q_{i}  \})$  to be the $\CC$-algebra generated by $U_{s_{\alpha}} Z_{\lambda}$ on which multiplication is deduced from that of $\cH_{0}(\Sigma, \{ q_{\alpha} \})$, that of $\CC[\Lambda]$ and the commutation rule of $U_{s_{\alpha}}$ and $Z_{\lambda}$, where $\alpha\in\Delta$ and $\lambda\in\Lambda$:
\begin{align}\label{eq:comm-rel-Z-U}
  Z_{\lambda} U_{s_{\alpha}} - U_{s_{\alpha}} Z_{s_{\alpha} \lambda} =
  \begin{cases}
    (q_{\alpha} - 1)\frac{Z_{\lambda} - Z_{s_{\alpha} \lambda}}{1 - Z_{-\alpha}}, &\quad\text{if $\alpha^{\vee} \not\in 2 \Lambda^{\vee}$} , \\
    (q_{\alpha} - 1 + Z_{-\alpha} ( (q_{\alpha} q_{i})^{\half} - (q_{\alpha} q_{i}^{-1})^{\half}  ))      \frac{Z_{\lambda} - Z_{s_{\alpha} \lambda}}{1 - Z_{-2\alpha}}, &\quad\text{if $\alpha^{\vee} \in 2 \Lambda^{\vee}$}.
  \end{cases}
\end{align}
Let $\cZ$ denote the centre of $\cH(\Sigma, \{ q_{\alpha} \}, \{ q_{i}  \})$ and let $K(\cZ)$ be its field of fractions.
Set $\cH(\Sigma, \{ q_{\alpha} \}, \{ q_{i}  \})_{K(\cZ)} = \cH(\Sigma, \{ q_{\alpha} \}, \{ q_{i}  \})\otimes_{\cZ} K(\cZ)$.
\end{defn}

\begin{rmk}
The Hecke algebra with parameters $\cH(\Sigma, \{ q_{\alpha} \}, \{ q_{i}  \})$ depends on the quintuplet $(\Lambda, \Lambda^{\vee}, \Sigma, \Sigma^{\vee},\Delta)$, but we do not try to reflect this in the notation.
\end{rmk}

Next we define a certain root datum and we will show that it is the underlying root datum of the Hecke algebra structure of $\End_{\tG}(i_{\tP}^{\tG} \tEBO)$.

\begin{defn}\label{defn:root-datum-Sigma-O} Let $\Lambda_{\cO}$ be the free $\ZZ$-module contained in $a_{M}$ given by the image of $\tM^{\widetilde{\sigma}}/\tM^{1}$ under $H_{\tM}$.
For $\alpha\in \Sigma_{\cO,\mu}$, set
\begin{equation*}
  \talpha  = H_{\tM}(\tilh_{\alpha}^{t_{\alpha}}),
\end{equation*}
 where $\tilh_{\alpha}$ and $t_{\alpha}$ are defined in Section~\ref{sec:intertw-oper} (below Remark~\ref{rmk:rep-inertial-orbit}).
Let
\begin{align*}
  \Sigma_{\cO}= \{ \talpha  \spacedvert \alpha \in \Sigma_{\cO,\mu} \},\qquad
  \Delta_{\cO}= \{ \talpha  \spacedvert \alpha \in \Delta_{\cO,\mu} \}.
\end{align*}
Let $\Lambda_{\cO}^{\vee}$ be the free $\ZZ$-module contained in $a_{M}^{*}$ in perfect duality with $\Lambda_{\cO}$ with respect to the pairing between $a_{M}$ and $a_{M}^{*}$.
Let
\begin{align*}
  \Sigma_{\cO}^{\vee} = \{ \alpha^{*} \in a_{M}^{*} \spacedvert \alpha^{*} \text{ is a multiple of } \alpha \text{ and } \forma[\alpha^{*}, \talpha ] = 2\}.
\end{align*}
\end{defn}

\begin{prop}\label{prop:root-datum-Sigma-O}
  The quadruple $(\Lambda_{\cO},\Lambda_{\cO}^{\vee},\Sigma_{\cO},\Sigma_{\cO}^{\vee})$ defined in Definition~\ref{defn:root-datum-Sigma-O} is a root datum.
  The underlying root system is reduced.
  The Weyl group of $\Sigma_{\cO}$ is canonically isomorphic to $W_{\cO}$.
  The set $\Delta_{\cO}$ is a base for $\Sigma_{\cO}$.
  Let $\Sigma$ be an irreducible component of $\Sigma_{\cO,\mu}$.
  Then the set $\widetilde{\Sigma} := \{ \talpha  \spacedvert \alpha\in \Sigma \}$ is an irreducible component of $\Sigma_{\cO}$ and
  $\widetilde{\Sigma}$ is of the same type as $\Sigma$ if $\Sigma$ is of type A, B or D and $\widetilde{\Sigma}$ is of type B, if $\Sigma$ is of type C.
\end{prop}

\begin{proof}
  The proof of \cite[Proposition 6.1]{Heiermann-MR2827179} generalizes, as we deal in fact effectively only with the underlying reductive group.
  Note that when $G$ is the unitary group, we should  take a prime element of $\Fs$ in the proof of \cite[Proposition 6.1]{Heiermann-MR2827179}.
\end{proof}

\begin{rmk}\label{rmk:root-in-2-Lambda-check}
 This includes the special case $B_1$, when $\alpha_{i,k_{i}}^{*}$ lies in $2 \Lambda_{\cO}^{\vee}$.
\end{rmk}

We now go back to the analysis of the algebra of intertwining operators $\End_{\tG}(i_{\tP}^{\tG} \tEBO)$.
We note that $W(\Sigma_{\cO,\mu})$ is canonically identified with $W(\Sigma_{\cO})$.
\begin{lemma}\label{lemma:bs-is-mostly-zero}
  Let $\Sigma_{i}$ be an irreducible component of $\Sigma_{\cO}$ and assume  $\alpha  \in \Delta_{\cO} \cap \Sigma_{i}$.
  Then $a_{s_{\alpha},-} \neq 0$ implies that $\Sigma_{i}$ is of type B and $\alpha$ is the short root.
\end{lemma}

\begin{proof} By \cite[Remark 1.3]{Heiermann-notes}, one has necessarily $a_{s_{\alpha},-}=0$ in the situation when $\mu^{M_{\alpha }}$ is associated to a supercuspidal representation of a maximal Levi subgroup of a general linear group.
  The argument therein applies also to a cover of a general linear group, when the Levi subgroups decompose, which is  assumed here.
  When $\mu^{M_{\alpha }}$ is associated to other types of Levi subgroups, then $\Sigma_{i}$ is of type $B$ and $\alpha $ is the short root.
\end{proof}

By the same reasoning as in Lemme~7.2 of \cite{Heiermann-MR2827179}, we have
\begin{lemma}\label{lemma:Ts-braid}
  Let $\alpha,\alpha'\in\Delta_{\cO}$ and set $s=s_{\alpha}$ and $s'=s_{\alpha'}$.
  Let $m\in\ZZ$ such that $(ss')^{m} = 1$.
  Then
  \begin{align*}
    T_{s} T_{s'} T_{s} \cdots =     T_{s'} T_{s} T_{s'} \cdots
  \end{align*}
  where there are $m$ factors on either side.
\end{lemma}

From now on the parameters $\{q_{\alpha}\}$ and $\{ q_{i} \}$ are chosen as follows.
For $\alpha\in\Delta_{\cO}$, set
\begin{align}\label{eq:q-alpha}
  q_{\alpha} = q_{\Fz}^{a_{s_{\alpha}}+a_{s_{\alpha},-}}
\end{align}
and if $\Sigma_{\cO,i}$ is of type B and $\alpha$ is the short root, set
\begin{align}\label{eq:q-i}
  q_{i} = q_{\Fz}^{a_{s_{\alpha}}-a_{s_{\alpha},-}}.
\end{align}

We identify $C = \CC[\Lambda]$ with $\BO$ via $Z_{\talpha} \mapsto X_{\alpha}$ for $\alpha\in \Sigma_{\cO,\mu}$ and let $T_{s_{\alpha}}$ play the role of the generators $U_{s_{\alpha}}$.

\begin{prop}\label{prop:BO-Tw-is-Hecke}
  The algebra $\bigoplus_{w\in W_{\cO}} \BO T_{w}$ is isomorphic to the Hecke algebra with parameters $\cH(\Sigma_{\cO}, \{ q_{\alpha} \}, \{ q_{i}  \})$
  associated to the based root datum $(\Lambda_{\cO},  \Lambda_{\cO}^{\vee}, \Sigma_{\cO}, \Sigma_{\cO}^{\vee}, \Delta_{\cO})$ where the parameters $\{ q_{\alpha} \}, \{ q_{i}\}$ are given as in \eqref{eq:q-alpha} and \eqref{eq:q-i}.
\end{prop}

\begin{proof}
  Proposition~\ref{prop:hecke-relation} gives the Hecke relation and Lemma~\ref{lemma:Ts-braid} gives the braid relation.
  A basic computation shows that \eqref{eq:comm-rel-Z-U} is also satisfied.
\end{proof}

The following relation is easy to verify by using Proposition~\ref{lemma:Jw-properties}\eqref{item:4} and the definition of $T_{w}$. %
\begin{prop}\label{prop:comm-relation-Jr-Tw}
  For $w\in W_{\cO}$ and $r\in R(\cO)$, we have $r^{-1}w r\in W_{\cO}$ and $T_{w}J_{r} = J_{r}T_{r^{-1}w r}$.
\end{prop}

Recall from Lemma~\ref{lemma:Jw-properties}~\eqref{item:3}  that, for  $r,r'\in R(\cO)$, $J_{r}J_{r'}= \eta(r,r') J_{rr'}$ for some $2$-cocycle $\eta$.
Write $\CC[R(\cO)]_{\eta}$ for the  group algebra on $R(\cO)$ twisted by $\eta$.
Combining Theorem~\ref{thm:Jr-Tw-basis}, Proposition~\ref{prop:BO-Tw-is-Hecke} and Proposition~\ref{prop:comm-relation-Jr-Tw}, we get:
\begin{thm}\label{thm:End-iEBO-is-Hecke}
  The algebra $\End_{\tG}(i_{\tP}^{\tG} \tEBO)$ is isomorphic to the semi-direct product $\CC[R(\cO)]_{\eta}\ltimes \cH(\Sigma_{\cO}, \{ q_{\alpha} \}, \{ q_{i}  \})$
  where the parameters $\{ q_{\alpha} \}, \{ q_{i}\}$ are given as in \eqref{eq:q-alpha} and \eqref{eq:q-i}.
\end{thm}

The following theorem is due to J.N.
Bernstein for $p$-adic reductive groups.
Its proof \cite[1.6]{Roche-Bernstein-decomp-MR1934305} generalizes to covering groups: one important ingredient is  Bernstein's second adjointness theorem which has been checked in \cite{Fratila-Prasad-MR4538041} to generalize to covering groups.
Other ingredients like properties of parabolic induction and Jacquet functor are known to be valid in the general context of $l$-groups \cite{Bernstein-MR771671} and thus are valid for covering groups.
The proof of the important proposition \cite[1.5]{Roche-Bernstein-decomp-MR1934305} is also contained in \cite{Heiermann-notes}.

\begin{thm}
The representation $i_{\tP}^{\tG} \tEBO$ is a projective generator of the category $\Rep_{\tG}(\cO)$.
In particular, $\Rep_{\tG}(\cO)$ is equivalent to the category of right-$\End_{\tG}(i_{\tP}^{\tG} \tEBO)$-modules.
\end{thm}

\begin{cor}\label{cor:bernstein-block-hecke-module}
  The category $\Rep_{\tG}(\cO)$ is equivalent to the category of right modules over the algebra $\CC[R(\cO)]_{\eta } \ltimes \cH(\Sigma_{\cO}, \{ q_{\alpha} \}, \{ q_{i}  \})$
  where the parameters $\{ q_{\alpha} \}, \{ q_{i}\}$ are given as in \eqref{eq:q-alpha} and \eqref{eq:q-i}.
\end{cor}

The isomorphism is compatible with parabolic inductions and Jacquet functors in the following sense.
Let $P'=M'U'$ be a parabolic subgroup of $G$ such that  $M' \supseteq M$  and $P' \supseteq P$.
We have a canonical injection $\End_{\tM'}(i_{\tP\cap \tM'}^{\tM'} \tEBO) \rightarrow \End_{\tG}(i_{\tP}^{\tG} \tEBO)$ given by $i_{\tP'}^{\tG}$.
We can define the various objects relative to the ambient group $M'$, since it is a product of general linear groups, symplectic groups, special orthogonal groups, full orthogonal groups and unitary groups.
In particular, we have the operators $T_{w}^{M'}$ and $J_{r}^{M'}$ for $w\in W_{\cO}^{M'}$ and $r\in R^{M'}(\cO)$ relative to the ambient group $M'$.
These operators are compatible with parabolic inductions.

The Jacquet functor $r_{\widetilde{\overline{P'}}}^{\tG}$ defines a functor from the category $\Rep_{\tG}(\cO)$ to $\Rep_{\tM'}(\cO)$.
The canonical injection $\End_{\tM'}(i_{\tP\cap \tM'}^{\tM'} \tEBO) \rightarrow \End_{\tG}(i_{\tP}^{\tG} \tEBO)$ gives a right $\End_{\tM'}(i_{\tP\cap \tM'}^{\tM'} \tEBO)$-module structure to every right $\End_{\tG}(i_{\tP}^{\tG} \tEBO)$-module.
By \cite[5.3]{Roche-Bernstein-decomp-MR1934305}, this commutes with the Jacquet functor.

\section{Preservation of Square-integrability and Temperedness --- decomposed Levi subgroups}
\label{sec:pres-square-integr}
Let $\cH = \End_{\tG} (i_{\tP}^{\tG}\tEBO )$, which, by  Theorem~\ref{thm:End-iEBO-is-Hecke}, is isomorphic to  the extended Hecke algebra 
$\CC[R(\cO)]_{\eta} \ltimes \cH(\Sigma_{\cO}, \{ q_{\alpha} \}, \{ q_{i}  \})$ where the affine Hecke algebra $\cH(\Sigma_{\cO}, \{ q_{\alpha} \}, \{ q_{i}  \})$ is attached to the based root datum $(\Lambda_{\cO},  \Lambda_{\cO}^{\vee}, \Sigma_{\cO}, \Sigma_{\cO}^{\vee}, \Delta_{\cO})$ defined in Section~\ref{sec:hecke-algebras}.
The aim of this section is to show (in Theorem~\ref{thm:preserv-sq-int-temper}) that  the equivalence of categories
\begin{align*}
  \Rep_{\tG}(\cO) &\rightarrow \text{ Right-}\cH\text{-}\MOD\\
  (\widetilde{\pi}, \tV) &\mapsto \Hom_{\tG}(i_{\tP}^{\tG}\tEBO, \tV) 
\end{align*}
preserves  square-integrability and temperedness (recalled below).
The proof is largely the same as in \cite{Heiermann-MR2975415}, though the definitions need to be adapted to the cover case.
Let $Z_{\tG}$ denote the centre of $\tG$.
Let $\cZ_{\tG}$ denote the image of $\tA_{G}^{\dagger}$ in $\Lambda_{\cO}$ under the Harish-Chandra map.
Then $\cZ_{\tG}$ is contained in the centre of $\cH$.

\begin{thm}\label{thm:preserv-sq-int-temper}
  Let $(\widetilde{\pi}, \tV)$ be an irreducible representation in $\Rep_{\tG}(\cO)$.
  Then the representation $(\widetilde{\pi},\tV)$ is tempered (resp. square integrable modulo $Z_{\tG}$) if and only if the $\End_{\tG} (i_{\tP}^{\tG}\tEBO )$-module $\Hom_{\tG}(i_{\tP}^{\tG}\tEBO, \tV)$ is  tempered (resp. square integrable modulo  $\cZ_{\tG}$).
\end{thm}

\subsection{Square-integrability and Temperedness of admissible representations of $\tG$}
\label{sec:square-integr-temp-rep}

Let $M'$ be a Levi subgroup of $G$.
Via the map $a_{M',\CC}^{*} \rightarrow \Xnr(\tM'), \lambda \mapsto \chi_{\lambda}$,
we define the real part of $\chi_{\lambda}\in\Xnr(\tM')$ to be the real part of $\lambda$ and denote it by  $\Re \chi_{\lambda}$.
This is well-defined.

Let $r_{\tP'}^{\tG}$ be the Jacquet functor with respect to the unipotent radical of the parabolic subgroup $P'$ of $G$.
The functor $r_{\tP'}^{\tG}$ is left-adjoint to $i_{\tP'}^{\tG}$ and $r_{\widetilde{\overline{P'}}}^{\tG}$  is right-adjoint to $i_{\tP'}^{\tG}$ by Bernstein's Second Adjointness Theorem generalised to the case of finite central extensions (\cite[Theorem~8.15]{Fratila-Prasad-MR4538041}).

We have the canonical decomposition $a_{M'}^{*} = a_{M'}^{G,*} \oplus a_{G}^{*}$.
Write $\lsup{+}{a}_{P'}^{G,*}$ for the positive open obtuse Weyl chamber of $a_{M'}^{G,*}$, namely, the linear combinations of elements in $\Delta(P')$ with  coefficients $>0$.
Write $\overline{\lsup{+}{a}_{P'}^{G,*}}$ for the positive closed obtuse Weyl chamber of $a_{M'}^{G,*}$, namely, the linear combinations of elements in $\Delta(P')$ with  coefficients $\ge 0$.

Let $(\widetilde{\pi},\tV)$ be an admissible representation of $\tM'$.
Exponents of $(\widetilde{\pi},\tV)$ are defined with respect to the action of $\tA_{M'}^{\dagger}$.
Those unramified characters $\chi$ of $\tA_{M'}^{\dagger}$ which have a  nontrivial eigenspace are called the exponents of  $(\widetilde{\pi},\tV)$.

For an  admissible representation of $\tG$, 
we have Casselman's criterion for square-integrability (resp. temperedness).
See \cite[III.1.1, III.2.2]{Waldspurger-plancherel-MR1989693} for  reductive groups and \cite[Theorem~3.4, Proposition~3.5]{Ban-Jantzen-MR3151110} for the finite central extensions.
The condition on the central character can be replaced by using the restriction to $\tA_{G}^{\dagger}$, since $Z_{\tG} \cap \tA_{G} / \tA_{G}^{\dagger}$ is  finite.

\begin{prop}
  Let $\widetilde{\pi}$ be an admissible representation of $\tG$.
  \begin{enumerate}
  \item The representation $\widetilde{\pi}$ is tempered, if and only if for all standard maximal  parabolic subgroups $P'$ of $G$, the real part of the exponents of $r_{\tP'}^{\tG}$ lie in $\overline{\lsup{+}{a}_{P'}^{G,*}}$.
  \item Assume that $\widetilde{\pi}$ admits a central character.
    The representation $\widetilde{\pi}$ is square-integrable (modulo centre), if and only if its central character is unitary and for all standard maximal  parabolic subgroups $P'$ of $G$, the real part of the exponents of $r_{\tP'}^{\tG}$ lie in $\lsup{+}{a}_{P'}^{G,*}$.
  \end{enumerate}
\end{prop}

With the same argument as \cite[Proposition~1.2]{Heiermann-MR2975415}, we have the following.
\begin{prop}\label{prop:casselman-criterion}
  Let $\widetilde{\pi}$ be an irreducible representation in $\Rep_{\tG}(\cO)$.
  \begin{enumerate}
  \item The representation $\widetilde{\pi}$ is tempered, if and only if the real parts of the exponents of $r_{\tP}^{\tG}$ lie in $\overline{\lsup{+}{a}_{P}^{G,*}}$.
  \item 
    The representation $\widetilde{\pi}$ is square-integrable, if and only if its central character is unitary and  the real part of the exponents of $r_{\tP}^{\tG}$ lie in $\lsup{+}{a}_{P}^{G,*}$.
  \end{enumerate}
\end{prop}

\subsection{Square-integrability and Temperedness of modules  over a Hecke algebra}
\label{sec:square-integr-temp-H-mod}

Let $V$ be an $\cH$-module of finite dimension.
Write $(\Lambda,\Lambda^{\vee},\Sigma,\Sigma^{\vee},\Delta)$ for the based root datum associated to $\cH$.
The exponents of $V$ are defined in terms of the action of  $\Lambda$.
Those  characters $\chi$ of $\Lambda$ with a nontrivial eigenspace are called the exponents of the $\cH$-module $V$.
A character $\chi$ of $\Lambda$ gives rise to a $\ZZ$-linear map $\Lambda \rightarrow \RR$ given by $\lambda \mapsto -\log |\chi(\lambda)|$, which we identify with an element in $\Lambda^{\vee}\otimes_{\ZZ}\RR$ and we call it the real part of $\chi$.

\begin{defn}\label{defn:H-mod-sq-int-tempered}
  Let $V$ be a right $\cH$-module of finite dimension. %
  \begin{enumerate}
  \item We say that $V$ is tempered, if the real part of the exponents of $V$ lie in the negative closed obtuse Weyl chamber of $\Lambda^{\vee}\otimes_{\ZZ}\RR$, namely, the linear combinations of elements in $\Delta^{\vee}$ with  coefficients $\le 0$.
  \item Let $Z$ be a sublattice of $\Lambda$ such that it lies in the centre of $\cH$ and it acts on $V$ by a scalar.
    We say that $V$ is square-integrable modulo $Z$, if the action of $Z$ is unitary, $Z\cup\Delta$ generates a $\ZZ$-module of finite index in $\Lambda$ and the real parts of the exponents of $V$ lie in the negative open obtuse Weyl chamber of $\Lambda^{\vee}\otimes_{\ZZ}\RR$, namely, the linear combinations of elements in $\Delta^{\vee}$ with coefficients $<0$.
  \end{enumerate}
\end{defn}

\subsection{Comparison of the real parts of the exponents}
\label{sec:preservation}

Recall that the supercuspidal representation $\widetilde{\sigma}$ of $\tM$ is the fixed base point in $\cO$.
Let $\chi_{\widetilde{\sigma}}$ denote the character which is the restriction of $\widetilde{\sigma}$ to $\tA_{M}^{\dagger}$.
The statements and proofs in Sections 3.1--3.3 of \cite{Heiermann-MR2975415} carry over once we replace the occurrences of central characters of representations of $M'$ there by characters of $\tA_{M'}^{\dagger}$,  where $M'$ is a Levi subgroup of $G$.
Let $B_{\tA_{M'}^{\dagger}}$ be the subalgebra of $\BO$ generated by $\{ b_{a'} \spacedvert a'\in A_{M'}^{\dagger} \}$ for  a Levi subgroup  $M'$ of $G$ that contains $M$.
Then we have analogous results to Propositions~3.2--3.3 and Corollary~3.4 of \cite{Heiermann-MR2975415}.

\begin{prop}\label{prop:central-action}
  Let $(\widetilde{\pi},\tV)$ be an irreducible representation contained in $\Rep_{\tG}(\cO)$.
  Let $\chi_{\widetilde{\pi}}$ be the restriction of $\widetilde{\pi}$ to $\tA_{G}^{\dagger}$.
  Then $\cZ_{\tG}$ acts by a character on $\Hom_{\tG}(i_{\tP}^{\tG}\tEBO, \tV)$ given by
  \begin{align*}
    z \mapsto (\chi_{\widetilde{\pi}} \chi_{\widetilde{\sigma}}^{-1})(z),
  \end{align*}
where  the expression   $(\chi_{\widetilde{\pi}} \chi_{\widetilde{\sigma}}^{-1})(z)$ means evaluating $\chi_{\widetilde{\pi}} \chi_{\widetilde{\sigma}}^{-1}$ on any element $\ta\in \tA_{G}^{\dagger}$ that maps to $z\in \cZ_{\tG}$ and is well-defined.
  In particular, the action of $\cZ_{\tG}$ is unitary if and only if $\chi_{\widetilde{\pi}}$ is unitary or equivalently if and only if the central character of $\widetilde{\pi}$ is unitary.
\end{prop}

\begin{prop}\label{prop:exp-Jacquet-mod-exp-BA-mod}
  Let $(\widetilde{\pi},\tV)$ be an irreducible representation of $\tG$  in $\Rep_{\tG}(\cO)$.
  Let $P'$ be a parabolic subgroup of $G$ which contains $P$ and whose standard Levi subgroup is $M'$.
  Then $\Hom_{\tG}(i_{\tP}^{\tG}\tEBO, \tV)$ decomposes into generalised eigenspaces for the action of $B_{\tA_{M'}^{\dagger}}$.
  The exponents of $\Hom_{\tG}(i_{\tP}^{\tG}\tEBO, \tV)$ with respect to $B_{\tA_{M'}^{\dagger}}$ are exactly the characters $\chi \chi_{\widetilde{\sigma}}^{-1}$ for $\chi$ running over the $\tA_{M'}^{\dagger}$-exponents of $r_{\widetilde{\overline{P'}}}^{\tG}\widetilde{\pi}$   such that 
\begin{equation*}
  \Hom_{\tM'}(i_{\tP\cap \tM'}^{\tM'}\tEBO, (r_{\widetilde{\overline{P'}}}^{\tG} \tV)_{\chi})\neq 0,
\end{equation*}
where $(r_{\widetilde{\overline{P'}}}^{\tG} \tV)_{\chi}$ denotes the $\chi$-generalised eigenspace of $r_{\widetilde{\overline{P'}}}^{\tG} \tV$.
\end{prop}

If we consider the real parts of the exponents, we can pass from characters of $B_{\tA_{M'}^{\dagger}}$ to those of $\BO$, by making use of the isomorphism $\Re \Xnr(\tM^{\widetilde{\sigma}}) \rightarrow \Re \Xnr(\tA_{M}^{\dagger})$.

\begin{cor}\label{cor:re-exp-same}
The set of the real parts of  the exponents of $\Hom_{\tG}(i_{\tP}^{\tG} \tEBO, \tV)$ with respect to  the action of $\BO$
is equal to the set of the real parts of  the exponents of the irreducible quotients of $r_{\widetilde{\overline{P}}}^{\tG}\widetilde{\pi}$ belonging to $\cO$.
\end{cor}

In general, as  elements  of $\Delta_{\cO}^{\vee}$ and those of $\Delta(P)$ are not non-zero multiples of each other, some further arguments are needed to prove Theorem~\ref{thm:preserv-sq-int-temper}.
For this, we note the following proposition, which is an analogue of \cite[Proposition~4.4]{Heiermann-MR2975415}, where we keep only the last statement which encompasses all the other statements.
Let $w\in W$ such that $wM$ is a standard Levi subgroup.
Write $P_{w}$ for the standard    parabolic subgroup with Levi factor $wM$.

\begin{prop}\label{prop:treat-other-quotients-of-Jacquet-mod}
  Let $(\widetilde{\pi}, \tV)$ be an irreducible representation in $\Rep_{\tG}(\cO)$.
  Assume that  $\widetilde{\sigma}''$ is an irreducible subquotient of $r_{\tP}^{\tG}(\widetilde{\pi})$.
 If $r\in W(M)$ satisfies $r\Delta_{\cO} = \Delta_{r\cO}$, then $r \widetilde{\sigma}''$ is an irreducible subquotient of $r_{\tP_{r}}^{\tG}\widetilde{\pi}$.
\end{prop}

We also need the analogue of Lemme~4.5 of \cite{Heiermann-MR2975415}.
\begin{lemma}\label{lemma:w-sigma}
  Let $\widetilde{\pi}$ be an irreducible representation in $\Rep_{\tG}(\cO)$.
  Assume that  $w\widetilde{\sigma}'$ is an irreducible subquotient of $r_{\tP}^{\tG}(\widetilde{\pi})$ for $w\in W(M)$ and $\widetilde{\sigma}'\in\cO$.
  Then there exist $w'\in W(M,\cO)$ and $r\in W(M)$ such that $w'\widetilde{\sigma}'$ is an irreducible subquotient of $r_{\tP}^{\tG}(\widetilde{\pi})$, $r \Delta_{\cO} = \Delta_{r\cO}$ and $w=rw'$.
\end{lemma}

\begin{proof}[Proof of Theorem~\ref{thm:preserv-sq-int-temper}]
We will indicate the idea of the proof, referring the details to \cite[Th\'eor\`eme~5]{Heiermann-MR2975415}.

  The conditions on the ``central'' characters have been treated in Proposition~\ref{prop:central-action} and we will focus on comparing the exponents.

  Assume that the $\End_{\tG} (i_{\tP}^{\tG}\tEBO )$-module $\Hom_{\tG}(i_{\tP}^{\tG}\tEBO, \tV)$ is  square integrable modulo  $\cZ_{\tG}$.
  Thus $\Delta_{\cO}^{\vee}$ must span $a_{M}^{G,*}$.
   Corollary~\ref{cor:re-exp-same}, Proposition~\ref{prop:treat-other-quotients-of-Jacquet-mod} and Lemma~\ref{lemma:w-sigma} show that the real parts of  the exponents of the irreducible quotients of $r_{\widetilde{\overline{P}}}^{\tG}\widetilde{\pi}$ are linear combinations of elements in $\Delta_{r \cO}^{\vee}$ with coefficients $<0$ for some appropriate elements $r\in W(M)$ satisfying $r\Delta_{\cO} = \Delta_{r\cO}$.
  Every element of $\Delta_{r\cO}^{\vee}$ is  a linear combination of elements in $\Delta(P)$ with coefficients $\ge 0$ and the change of basis matrix is invertible, since they both  span $a_{M}^{G,*}$.
  Thus the real parts of the exponents are  linear combinations of elements of $\Delta(P)$ with coefficients $<0$.

  For the opposite direction, assume that $(\widetilde{\pi},\tV)$ is square integrable.
  By \cite[Theorem 3.9.1]{Silberger-MR636960}\footnote{The proof of this result of Silberger is based on the Casselman criterion for square integrability and properties of the $\mu $-function both of which carry over to covers of reductive groups as already remarked several times in this paper.}, $\Delta_{\cO}^{\vee}$ must span $a_{M}^{G,*}$.
  The  real parts of the exponents of $r_{\widetilde{\overline{P}}}^{\tG}\widetilde{\pi}$ are linear combinations of $\Delta(P)$ with coefficients $<0$ and so are those of the $\End_{\tG} (i_{\tP}^{\tG}\tEBO )$-module $\Hom_{\tG}(i_{\tP}^{\tG}\tEBO, \tV)$,   by Corollary~\ref{cor:re-exp-same}.
  Write the real parts as linear combinations in terms of $\Delta_{\cO}^{\vee}$ and we must show that the coefficients are $<0$.
  By applying an appropriate $r\in W(M)$ satisfying $r\Delta_{\cO} = \Delta_{r\cO}$, one can show that every  coefficient with respect to $\Delta_{\cO}^{\vee}$ is equal to  a certain coefficient of the real part of an exponent of $r_{\widetilde{\overline{P_{r}}}}^{\tG}\widetilde{\pi}$ with respect to $\Delta_{P_{r}}$ and thus it is indeed $<0$.

  The idea of proofs for preservation of temperedness is analogous.
\end{proof}

\section{Some results of M\oe glin on the extended Langlands--Deligne parameters for metaplectic groups}
\label{sec:some-results-moe}

First we set up some notations for the metaplectic groups and then we collect  M\oe glin's results  on the extended Langlands--Deligne parameters of supercuspidal representations of the  metaplectic groups.

Now we assume, in addition, that $\Fz$ is of characteristic $0$\footnote{This is needed for the theory of Langlands--Deligne parameters in the context of classical groups and metaplectic groups, as we follow \cite{Moeglin-MR2767522}.}.
Write $\tG$ for $\Mp_{2n}$ which  is  the $8$-fold  metaplectic  cover of $\mathbb{Sp}_{2n}(\Fz)$ (with the Leray cocycle).
It is defined as follows.
Let $\bpsi$ be a non-trivial additive character $\Fz\rightarrow \CC^{1}$, where $\CC^{1}$ denotes the norm $1$ subgroup of $\CC^{\times}$.
The metaplectic group $\Mp_{2n}$ considered here is a topological  $8$-fold cover of the symplectic group $\Sp_{2n}$.
It can be regarded as  the set $\Sp_{2n}\times\bmu_{8}$ equipped with the multiplication determined by the Leray cocycle $c_{\bpsi}$ (depending on $\bpsi$) which can be found in \cite[Theorem~5.2]{GKT-theta-book-Jun26}:
\begin{align*}
  (g_{1},\zeta_{1})(g_{2},\zeta_{2}) = (g_{1}g_{2}, \zeta_{1}\zeta_{2}c_{\bpsi}(g_{1},g_{2}) ).
\end{align*}
The advantage of using this realisation of the metaplectic group is that it splits over the general linear factors of the standard Levi subgroups $M$ of $\Sp_{2n}$, so that $\tM \isom \GL_{d_{1}} \times \cdots \times \GL_{d_{r}} \times \Mp_{2n'}$ for some appropriate integers $d_{1},\ldots, d_{r}, n'$.
We will suppress $\bpsi$ from the notation, though some constructions below, such as, the Howe correspondence and the Langlands--Deligne parameters for $\Mp_{2n}$  depend on $\bpsi$.
The name, metaplectic group, also refers to the non-trivial double cover of $\Sp_{2n}$.
Let us denote it by $\Mp^{(2)}_{2n}$.
It can be regarded as the set $\Sp_{2n}\times\bmu_{2}$ equipped with the multiplication determined by the Rao cocycle $c_{\Rao}$ which is independent of $\bpsi$ (see, for example, \cite[Theorem~5.10]{GKT-theta-book-Jun26}).
The categories of genuine representations of $\Mp^{(2)}_{2n}$ and of $\Mp_{2n}$ are equivalent, which can be seen as follows.
The Leray cocycle and the Rao cocycle are cohomologous.
There exists a  function $\beta : \Sp_{2n} \rightarrow \bmu_{8}$ such that $c_{\Rao} (g_{1},g_{2}) = c_{\bpsi}(g_{1},g_{2}) \beta(g_{1}) \beta(g_{2}) \beta(g_{1}g_{2})^{-1}$.
With the realisations as above, we get an embedding
\begin{align*}
  \Mp^{(2)}_{2n} &\rightarrow \Mp_{2n}\\
  (g, \zeta) &\mapsto (g, \zeta \beta(g)) ,
\end{align*}
and we note that $\beta(1) = 1$.
The embedding is continuous by the automatic continuity theorem (\cite[Theorem~6.49]{GKT-theta-book-Jun26}), since the Leray cocycle and the Rao cocycle are measurable (\cite[Proposition~8.19]{GKT-theta-book-Jun26}).
Given a genuine representation $\widetilde{\pi}$ of $\Mp^{(2)}_{2n}$, we can extend it to a genuine representation of $\Mp_{2n}$ by requiring that $\bmu_{8} \subset \CC^\times$ acts by multiplication on the representation space of $\widetilde{\pi}$,
whereas given a genuine representation $\widetilde{\pi}$ of $\Mp_{2n}$, it restricts to a genuine representation of $\Mp^{(2)}_{2n}$.
Even though $\Mp^{(2)}_{2n}$ is more canonical in the sense that the Rao cocycle $c_{\Rao}$  is independent of $\bpsi$, it has the disadvantage that it does not split over the general linear factors of the standard Levi subgroups $M$ of $\Sp_{2n}$.

Below, all representations of $\tG$ are assumed to be genuine, namely, that $z\in\bmu_{8} \subset \CC^{1}$ acts by multiplication by $z$.
We will adopt the Bernstein--Zelevinsky--Tadi\'c notation for parabolically induced representations.
Then an induction from the standard parabolic with the standard Levi subgroup isomorphic to  $\GL_{d_{1}} \times \cdots \times \GL_{d_{r}} \times \Mp_{2n'}$ will be denoted by $\rho_{1}\times \cdots \times \rho_{r} \rtimes \widetilde{\tau}$ where $\rho_{i}$ is a representation of $\GL_{d_{i}}$ and $\widetilde{\tau}$ is a representation of $\Mp_{2n'}$.
Similar notation will also be used in the case of classical groups.

The $L$-group for $\tG$ is  $\lsup{L}{\tG} := \Sp_{2n}(\CC)$.
Let $\iota : \lsup{L}{\tG} = \Sp_{2n}(\CC) \rightarrow \GL_{2n}(\CC)$ be the standard representation of $\lsup{L}{\tG}$.
Let $W_{\Fz}$ denote the Weil group of $\Fz$ and let $I_{\Fz}$ denote the inertia group.
Fix a Frobenius automorphism $\Fr\in W_{\Fz}$ so that $W_{\Fz} = \langle \Fr \rangle \ltimes I_{\Fz}$.
A character of $W_{\Fz}$ is said to be unramified if it is trivial on $I_{\Fz}$.
By the local class field theory, such a character is identified with a character of $\Fz^{\times}$ that is trivial on $\cO_{\Fz}^{\times}$.
We identify the unramified character of $W_{\Fz}$ that sends $\Fr^{-1}$ to $q_{\Fz}$ with the unramified character $|\ |_{\Fz}$ of $\Fz^{\times}$.

 A Langlands parameter for $\tG$ is defined to be a continuous homomorphism $\rho: W_{\Fz} \rightarrow \lsup{L}{\tG}$ which sends $\Fr$ to a semisimple element.
A continuous homomorphism $\rho: W_{\Fz}\times \SL_{2}(\CC) \rightarrow \lsup{L}{\tG}$ is said to be a Langlands--Deligne parameter for $\tG$, if its restriction to $W_{\Fz}$ is a Langlands parameter and its restriction to $\SL_{2}(\CC)$ is a morphism of algebraic groups.
A Langlands parameter (resp. Langlands--Deligne parameter) is said to be discrete, if its image is not contained in a proper Levi subgroup.
Two Langlands parameters (resp.  Langlands--Deligne parameters) are said to be equivalent, if they are conjugate by an element in $\lsup{L}{\tG}$.

Let $\rho$ and $\rho'$ be  irreducible representations of $W_{\Fz}$.
We say that $\rho$ and $\rho'$ are in the same inertial class if $\rho \isom \rho'|\ |_{\Fz}^{s}$ for some $s\in\CC$.
The group of unramified characters of $W_{\Fz}$ acts on the inertial class of $\rho$ by twisting.
Let $t_{\rho}$ denote the order of the stabiliser for this action.
It is clear that $t_{\rho} = t_{\rho'}$ if $\rho$ and $\rho'$ are in the same inertial class.

Let $\rho$ be a self-dual representation of $W_{\Fz}$.
We say that $\rho$ is of type $\lsup{L}{\tG}$, if it factors through a group of type $\lsup{L}{\tG}$; otherwise we say that $\rho$ is not of type $\lsup{L}{\tG}$.
We stress that the use of either of these notions will presume that $\rho$ is self-dual.

For $a\in\ZZ_{\ge 1}$, write $\spr(a)$ for the unique irreducible representation of $\SL_{2}(\CC)$ of dimension $a$.

We state now a result on extended Langlands--Deligne parameters for metaplectic groups along the line of  \cite[Theorem~1.1]{Heiermann-MR3719526} where supercuspidal representations of quasi-split classical groups and their pure inner forms are considered.
It is a reformulation of  \cite[Th\'eor\`eme 2.5.1]{Moeglin-MR2767522}.
Even though \cite[Th\'eor\`eme 2.5.1]{Moeglin-MR2767522} excludes the metaplectic case, arguments in Sections~2.2, 3.1, 3.2 of \cite{Moeglin-MR2767522} enable one to obtain an analogous result for the metaplectic case.
After obtaining the extended Langlands--Deligne parameters for metaplectic groups, we recall how to read off the points of reducibility for representations of metaplectic groups induced from maximal parabolic subgroups.

\begin{thm}\label{thm:LD-param-and-char-comp-group}
  \begin{enumerate}
  \item
    The $L$-packet associated to a Langlands--Deligne parameter $\varphi: W_{F} \times \SL_{2}(\CC) \rightarrow \lsup{L}{\tG}$ contains a supercuspidal representation of $\tG$, if and only if
    \begin{align}\label{eq:discrete-LD-param}
      \iota \circ \varphi = \bigoplus_{\rho \text{ not of type }  \lsup{L}{\tG}} \left( \bigoplus_{k=1}^{a_{\rho}} (\rho\otimes \spr(2k)) \right) \oplus \bigoplus_{\rho \text{ of type }  \lsup{L}{\tG}} \left( \bigoplus_{k=1}^{a_{\rho}} (\rho\otimes \spr(2k-1)) \right)
    \end{align}
    where the $\rho$'s are irreducible representations and the $a_{\rho}$'s are non-negative integers.
    Such a Langlands--Deligne parameter is a discrete Langlands--Deligne parameter.
    Any Langlands--Deligne parameter $W_{\Fz} \times \SL_{2}(\CC) \rightarrow \GL_{2n}(\CC)$ of the form on the right hand side of \eqref{eq:discrete-LD-param} factors through $\lsup{L}{\tG}$.
  \item \label{item:LD-1}Given $\varphi$ satisfying \eqref{eq:discrete-LD-param}, denote by $z_{\varphi,\rho,k}$ the diagonal matrix in $\lsup{L}{\tG}$ that acts by $-1$ on the space of the direct summand $\rho\otimes \spr(2k)$ (resp. $\rho\otimes \spr(2k-1$)) of $\iota \circ \varphi$  and by $1$ elsewhere.
    Put $\scrS_{\varphi} = C_{\lsup{L}{\tG}}(\Im(\varphi)) / C_{\lsup{L}{\tG}}(\Im(\varphi))^{\circ}$.
    The elements $z_{\varphi,\rho,k}$ lie in $C_{\lsup{L}{\tG}}(\Im(\varphi))$ and their images $\overline{z}_{\varphi,\rho,k}$ generate the commutative group $\scrS_{\varphi}$.
    A pair $(\varphi,\epsilon)$ formed by a discrete Langlands--Deligne parameter as in \eqref{eq:discrete-LD-param} and a character $\epsilon$ of $\scrS_{\varphi}$ corresponds to a supercuspidal representation of  $\tG$, if and only if $\epsilon$ is alternating, i.e. $\epsilon(\overline{z}_{\varphi,\rho,k}) = (-1)^{k-1} \epsilon(\overline{z}_{\varphi,\rho,1})$ with $\epsilon(\overline{z}_{\varphi,\rho,1}) = -1$ for $\rho$ not of type $\lsup{L}{\tG}$ and $\epsilon(\overline{z}_{\varphi,\rho,1}) \in \{\pm 1\} $ for $\rho$  of type $\lsup{L}{\tG}$.
  \item
    Suppose that $\varphi$ is as in \eqref{eq:discrete-LD-param}.
    Let $t_{o}$ be the number of $\rho$ of type $\lsup{L}{\tG}$ with $a_{\rho}$ odd and let $t_{e}$ be the number of  $\rho$ of type $\lsup{L}{\tG}$ for which $a_{\rho}$ is even.
    Then there are $2^{t_{o} + t_{e}}$ non-isomorphic supercuspidal representations of $\tG$ with Langlands--Deligne parameter $\varphi$.\footnote{Of course,    $t_{o}+t_{e}$ is just the number of $\rho$ of type $\lsup{L}{G}$.  The metaplectic group has no other pure inner form and there is no extra constraint on the alternating characters. The way this is stated is meant to supplement \cite[Theorem~1.1]{Heiermann-MR3719526}.
}
   \end{enumerate}
\end{thm}

\begin{rmk}
  We follow M\oe glin's construction in \cite{Moeglin-MR2767522} which uses  theta correspondence in the stable range to transfer extended Langlands--Deligne parameters from  special odd orthogonal groups to metaplectic groups.
  We will recall her construction below.
  This is different from the transfer induced by the map  in Theorem~1.1 of \cite{Gan-Gordan-MR2999299} where equal-rank theta correspondence is used.
  The Langlands--Deligne parameters for both methods match but the characters of the component group may differ.
  The discrepancy is determined in Proposition~7.2(1) of \cite{chen2026constructionlocalarthurpackets}.
  To convert from M\oe glin's normalisation to  Gan--Savin's normalisation, in the extended Langlands--Deligne parameter of Theorem~\ref{thm:LD-param-and-char-comp-group}, $\varphi$ remains the same while $\epsilon(\overline{z}_{\varphi,\rho,k})$ is changed to its negative exactly when $\rho=\triv$.
  We note that $\triv$ is not of type $\lsup{L}{\tG}$, as $\tG=\Mp_{2n}$ and $\lsup{L}{\tG}=\Sp_{2n}(\CC)$.
  Thus if we adopted Gan--Savin's normalisation, the definition of an alternating character would have a different initial condition exactly when $\rho=\triv$, namely, $\epsilon(\overline{z}_{\varphi,\triv,1}) = 1$.
\end{rmk}

\begin{rmk}
  Before  explaining M\oe glin's construction,
    first we remark that we can remove certain constraints from \cite{Moeglin-MR2767522}.
  Recall that $\Fz$ is a $p$-adic field.
  We can remove the condition on the residual characteristic that $p\neq 2$ in \cite{Moeglin-MR2767522}, since the Howe duality conjecture for residue characteristic $2$  is now also proven \cite{Gan-Takeda-MR3454380}.

  We also remove the unitarity assumption on discrete series in \cite[Section~3]{Moeglin-MR2767522}.
  We note that the discrete series are defined in terms of the exponents of their Jacquet modules while square-integrable representations are defined in terms of their matrix coefficients.
  \cite[Theorem 3.4]{Ban-Jantzen-MR3151110}  shows that the Casselman criterion for square-integrability holds for finite central extensions of $p$-adic groups and
  \cite[Lemme III.1.3]{Waldspurger-plancherel-MR1989693}   shows that square-integrable representations are unitary for $p$-adic groups and the proof  remains valid for their finite central extensions.
  Combining these, we find that a discrete series representation of $\tG$ is unitary.
  A fortiori, a supercuspidal representation of $\tG$ is unitary.
\end{rmk}

\begin{defn}\label{defn:Salt}
  We say that a Langlands--Deligne parameter is without holes if it satisfies \eqref{eq:discrete-LD-param}.
  Given a Langlands--Deligne parameter $\varphi$ without holes, write $\hat{\scrS}_{\varphi}$ for the set of characters of the component group $\scrS_{\varphi}$ and let $\hat{\scrS}_{\varphi}^{\alt}$ be the subset consisting of alternating characters, as defined in Theorem~\ref{thm:LD-param-and-char-comp-group}\eqref{item:LD-1}. 
\end{defn}

Following \cite{Moeglin-MR2767522}, we will make use of the Howe correspondence of $\Mp_{2n}$ with odd orthogonal groups to associate an extended Langlands--Deligne parameter to a supercuspidal representation of $\Mp_{2n}$.
We introduce some notation for the Howe correspondence.
Let $m$ be a positive odd integer and $\zeta \in \{\pm 1\}$.
We sometimes just write $\zeta = +$ or $-$ without the $1$.
Let $\rO_{m}^{\zeta}$ denote the orthogonal group associated to  the quadratic space of dimension $m$ with discriminant $1$ and Hasse invariant $\zeta$.
We note that $\rO_{m}^{+}$ is split while $\rO_{m}^{-}$ is non-split.
Let $\widetilde{\sigma }$ be an irreducible supercuspidal representation of $\Mp_{2n}$ and let $\rho$ be an irreducible unitary supercuspidal representation of $\GL_{d_{\rho}}$, for some positive integer $d_{\rho}$. Denote by $\Omega_{n,m,\zeta}$ the Weil representation of $\Mp_{2n}\times\rO_{m}^{\zeta}$, which  is the pullback of the Weil representation of $\Mp_{2nm}$ along the natural homomorphism  $\Mp_{2n}\times\rO_{m}^{\zeta}\rightarrow \Mp_{2nm}$.
We have the Howe correspondence between $\Mp_{2n}$ and $\rO_{m}^{\zeta}$.
For fixed $\zeta$, the smallest value of $m$ such that there exists a nonzero representation $\pi^{\zeta}$ of $\rO_{m}^{\zeta}$ such that $\widetilde{\sigma }\otimes \pi^{\zeta}$ is a quotient of the Weil representation $\Omega_{n,m,\zeta}$
is called the first occurrence index of $\widetilde{\sigma }$ and $\pi^{\zeta}$ is called the first occurrence representation for $\widetilde{\sigma }$,
which is shown by Kudla to be supercuspidal.

\begin{lemma}[{\cite[Lemme 1.2]{Moeglin-MR2767522}}]\label{lemma:point-red-SO-Mp-relation}
  For a given $\zeta$, let $m^{\zeta}$ denote the first occurrence index of $\widetilde{\sigma }$ with respect to the tower $\{ \rO_{m}^{\zeta} \spacedvert m=1,3,5,\ldots \}$ and let $\pi^{\zeta}$ denote  the first occurrence representation for $\widetilde{\sigma }$.
  Let $\rho$ be an irreducible unitary supercuspidal representation of $\GL_{d_{\rho}}$.
Then the following hold.
\begin{enumerate}
\item If $\rho$ is not self-dual, then $\rho |\ |_{\Fz}^{x} \rtimes \widetilde{\sigma }$ is irreducible for all $x\in\RR$.
\item If $\rho$ is self-dual, then there is exactly one point of reducibility in $\RR_{\ge 0}$, which we denote by $x_{\rho,\widetilde{\sigma }}$.
  It can be described in terms of the point of reducibility $x_{\rho,\pi^{\zeta}}\in\RR_{\ge 0}$ for $\rho |\ |_{\Fz}^{x} \rtimes \pi^{\zeta}$ as follows.
  \begin{enumerate}
  \item If $\rho \neq \triv$, then  $x_{\rho,\widetilde{\sigma }} = x_{\rho,\pi^{\zeta}}$.
  \item     If $\rho = \triv$, then $x_{\rho,\widetilde{\sigma }} = | \half (2n - (m^{\zeta} - 1) + 1) |$, which is a half integer but not an integer.
  \end{enumerate}
\end{enumerate}
\end{lemma}

We now summarise the process for assigning   an extended Langlands--Deligne parameter $(\varphi,\epsilon)$ to a supercuspidal representation of $\tG=\Mp_{2n}$,  where $\varphi$ is  without holes and  $\epsilon\in \hat{\scrS}_{\varphi}^{\alt}$.
Then we outline the opposite process  for constructing a supercuspidal representation of $\tG=\Mp_{2n}$ from  such an extended Langlands--Deligne parameter.

For a Langlands--Deligne parameter $\varphi: W_{F} \times \SL_{2}(\CC) \rightarrow \lsup{L}{\tG}$, $\iota\circ\varphi$ is a direct sum of irreducible representations $\rho\otimes\spr(a)$ of $W_{\Fz} \times \SL_{2}(\CC)$ and we set $\Jord(\varphi)$ to be the multiset of pairs $(\rho,a)$.

Given a supercuspidal representation $\widetilde{\sigma }$ of $\tG$, we associate an extended Langlands--Deligne parameter $(\varphi,\epsilon)$ as follows.
Let $x_{\rho,\widetilde{\sigma }}\ge 0$ denote  the non-negative point of reducibility of $\rho|\ |^{x}\times \widetilde{\sigma }$, for an irreducible self-dual supercuspidal representation $\rho$ of $\GL_{d_{\rho}}$.
Let
\begin{align}\label{eq:Jord-pi}
  \Jord(\widetilde{\sigma }) = \bigcup_{\rho} \{ (\rho, 2x_{\rho,\widetilde{\sigma }} + 1 - 2\ell) \spacedvert \ell = 1, \ldots, \lfloor x_{\rho,\widetilde{\sigma }} \rfloor \}.
\end{align}
with $\rho$ running over irreducible self-dual supercuspidal representations of general linear groups.
We remark that in \eqref{eq:Jord-pi}, if $x_{\rho,\widetilde{\sigma }} = 0, \half$, then  it is understood that $\rho$ contributes an empty set.
Associate to $\Jord(\widetilde{\sigma })$ the representation of $W_{\Fz} \times \SL_{2}(\CC)$ given by
\begin{align}\label{eq:param-GL-assoc-to-Jord}
  \bigoplus_{(\rho, 2x_{\rho,\widetilde{\sigma }} + 1 - 2\ell)\in \Jord(\widetilde{\sigma })} \rho \otimes \spr(2x_{\rho,\widetilde{\sigma }} + 1 - 2\ell).
\end{align}
It factors through $\lsup{L}{\tG}$ and gives a discrete Langlands--Deligne parameter
\begin{align*}
 \varphi: W_{\Fz} \times \SL_{2}(\CC) \rightarrow \lsup{L}{\tG}
\end{align*}
such that $\iota\circ\varphi$ is isomorphic to \eqref{eq:param-GL-assoc-to-Jord}.
Then $\varphi$ is the Langlands--Deligne parameter associated to $\widetilde{\sigma }$.
The set of characters of the component group associated to $\varphi$ can be identified to  the set of functions $\Jord(\varphi) \rightarrow \{ \pm 1  \}$.
We will sometimes refer to such functions as characters of $\Jord(\varphi)$.
We now associate such a character to the supercuspidal representation $\widetilde{\sigma }$ of $\tG$.
Let $N$ be large enough relative to $\widetilde{\sigma }$.
This means $N>a$ where $a$ is the maximal $a'$ among all $(\triv,a') \in \Jord(\widetilde{\sigma })$.
Let $\pi_{N}^{\zeta}$ denote the image of $\widetilde{\sigma }$ under the Howe correspondence to $\rO_{2n+1+N}^{\zeta}$.
The representation $\pi_{N}^{\zeta}$ gives rise to a character $\epsilon_{N}^{\zeta}$ on the Arthur parameter  $\Jord(\varphi) \sqcup \{ (1,1,N) \}$ (see \cite[Section~3.2]{Moeglin-MR2767522} for the definitions of Arthur parameters and the factor $(1,1,N)$). Then we define $\epsilon$ to be the restriction of $\epsilon_{N}^{\zeta}$ to  $\Jord(\varphi)$.
It can be shown that $\epsilon$ does not depend on the choice of $N$ or $\zeta$.
We have associated an extended Langlands--Deligne parameter $(\varphi,\epsilon)$ to $\widetilde{\sigma }$.
It is shown in \cite[Section~3.2]{Moeglin-MR2767522} that $\varphi$ is without holes and $\epsilon$ is alternating.

We consider now the opposite direction.
Given a Langlands--Deligne parameter $\varphi: W_{\Fz} \times \SL_{2}(\CC) \rightarrow \lsup{L}{\tG}$ without holes  and a character $\epsilon \in \hat{\scrS}_{\varphi}^{\alt}$, we associate a supercuspidal representation $\widetilde{\sigma }$ of $\tG$ as follows.

Assume first that $\Jord(\varphi)$ does not contain  $(\triv,a')$ for $a'\in\ZZ_{>0}$.
Put $\zeta = \epsilon(-1)$ where $-1$ is the nontrivial element in the center of $\lsup{L}{\tG}$.
We regard $\varphi$ as a Langlands--Deligne parameter for $\SO_{2n+1}^{\zeta}$.
Then let $\pi_{0}$ be the supercuspidal representation of $\SO_{2n+1}^{\zeta}$ associated to $(\varphi,\epsilon)$.
By the conservation relation of first occurrences \cite{Sun-Zhu-MR3369906}, exactly one of the extensions $\pi$ of $\pi_{0}$ to $\rO_{2n+1}^{\zeta}$ has a nonzero image under the Howe correspondence to $\Mp_{2n}$.
Let $\widetilde{\sigma }$ denote this image which is a representation of $\Mp_{2n}$.
It is shown in \cite[Section~3.2]{Moeglin-MR2767522} that $\widetilde{\sigma }$ is supercuspidal and has $(\varphi,\epsilon)$ as its extended Langlands--Deligne parameter.

Next assume that $\Jord(\varphi)$ contains a factor of the form $(\triv,a')$ for $a'\in\ZZ_{>0}$.
Then $a'$ must be even.
Let $a_{0}$ be the largest $a'$ among all pairs $(\triv,a')$ in $\Jord(\varphi)$.
Consider $\Jord(\varphi')$ which is $\Jord(\varphi)$ but with the pair $(\triv,a_{0})$ removed and $\epsilon'$ which is the restriction of $\epsilon$ to $\Jord(\varphi')$.
Then $\varphi'$ can be regarded as a Langlands--Deligne parameter for $\SO_{2n+1-a_{0}}^{\zeta'}$ with   $\zeta' = \epsilon'(-1)$, where $-1$ is the nontrivial element in the centre of the $L$-group $\Sp_{2n-a_{0}}(\CC)$ of $\SO_{2n+1-a_{0}}^{\zeta'}$.
Let $\pi_{0}'$ be the supercuspidal representation of $\SO_{2n+1-a_{0}}^{\zeta'}$ associated to $(\varphi',\epsilon')$.
By the conservation relation of first occurrences \cite{Sun-Zhu-MR3369906}, exactly one of the extensions $\pi'$ of $\pi_{0}'$ to $\rO_{2n+1-a_{0}}^{\zeta'}$ has first occurrence index  strictly greater than $2n+1-a_{0}$.
Let $\widetilde{\sigma }$ be that first occurrence representation of $\Mp_{2n'}$ with $2n' > 2n+1-a_{0}$.
It is shown in \cite[Section~3.2]{Moeglin-MR2767522} that the extended Langlands--Deligne parameter of $\widetilde{\sigma }$ is $(\varphi,\epsilon)$, so in particular, $2n'=2n$.

We have explained that the  results of \cite{Moeglin-MR2767522} give us part (1) and (2) of Theorem~\ref{thm:LD-param-and-char-comp-group}.
Next consider part (3) of Theorem~\ref{thm:LD-param-and-char-comp-group}.
Given a Langlands--Deligne parameter satisfying \eqref{eq:discrete-LD-param}, we count the number of alternating characters.
It is clearly given by $2^{t_{o}+t_{e}}$.

We note a consequence of M\oe glin's construction, which shows that points of reducibility and Langlands--Deligne parameters determine each other.
See the  proof of Th\'eor\`eme~2.3.1 in \cite{Moeglin-MR2767522} for the case of classical groups.
See also Proposition~3.2 and Corollary~9.1 of \cite{Xu-MR3713922}.
Using Lemma~\ref{lemma:point-red-SO-Mp-relation}, we  transfer the result to the case of metaplectic groups to get the following proposition.

\begin{prop}\label{prop:Jord-point-red}
  Let $\widetilde{\sigma }$ be an irreducible supercuspidal representation of $\tG$ with Langlands--Deligne parameter $\varphi$.
  Assume that $\rho$ is a self-dual unitary supercuspidal representation of a general linear group.
 If $\rho$  occurs in $\Jord(\varphi)$, let
\begin{equation*}
  a_{\rho,\varphi} = \max \{ \ell  \spacedvert (\rho, \ell) \in \Jord(\varphi) \}.
\end{equation*}
  If $\rho$ does not occur in $\Jord(\varphi)$, then set
  \begin{align*}
    a_{\rho,\varphi} =
    \begin{cases}
      0 , &\quad\text{if $\rho$ is not of type $\lsup{L}{\tG}$};\\
      -1 ,& \quad\text{if $\rho$ is  of type $\lsup{L}{\tG}$}.
    \end{cases}
  \end{align*}
  Set $x_{\rho,\varphi} = \half(a_{\rho,\varphi} + 1)$.
  Let $x_{\rho,\widetilde{\sigma }} \ge 0$ be the non-negative point of reducibility for $\rho |\ |_{\Fz}^{x} \rtimes \widetilde{\sigma }$.
  Then $x_{\rho,\varphi} = x_{\rho,\widetilde{\sigma }}$.
\end{prop}

\begin{proof}
  Assume that $\iota\circ \varphi_{\widetilde{\sigma }}$ has no factor of the form $\triv\otimes \spr(a)$, where we write $\varphi_{\widetilde{\sigma }}$ for the Langlands--Deligne parameter associated to $\widetilde{\sigma }$.
Then there exists a supercuspidal representation $\pi$ of $\SO_{2n+1}^{\zeta}$ for some appropriate $\zeta$ that is in Howe correspondence with $\widetilde{\sigma }$ and $\varphi_{\pi} = \varphi_{\widetilde{\sigma }}$, where we write $\varphi_{\pi}$ for the Langlands--Deligne parameter associated to $\pi$.
Then by Lemma~\ref{lemma:point-red-SO-Mp-relation}, if $\rho\neq \triv$, we find $x_{\rho,\varphi_{\widetilde{\sigma }}} = x_{\rho,\varphi_{\pi}} = x_{\rho,\pi} = x_{\rho,\widetilde{\sigma }}$ and if $\rho=\triv$, we find $x_{\triv,\widetilde{\sigma }} = \half | 2n - (2n+1-1) + 1 | = \half  = x_{\triv,\varphi_{\widetilde{\sigma }}}$, noting that $\triv$ is not of type $\lsup{L}{\tG}$.
Assume that $\iota\circ \varphi_{\widetilde{\sigma }}$ has a factor of the form $\triv\otimes \spr(a)$ and we let $a_{0}$ be the largest $a$ among all such factors.
Then there exists a supercuspidal representation $\pi$ of $\SO_{2n+1-a_{0}}^{\zeta}$ for some appropriate $\zeta$ that is in Howe correspondence with $\widetilde{\sigma }$ and satisfying $\iota\circ \varphi_{\pi}\oplus \triv\otimes \spr(a_{0})  = \iota\circ\varphi_{\widetilde{\sigma }}$.
Then by Lemma~\ref{lemma:point-red-SO-Mp-relation}, if $\rho\neq \triv$, then $x_{\rho,\varphi_{\widetilde{\sigma }}} = x_{\rho,\varphi_{\pi}} = x_{\rho,\pi} = x_{\rho,\widetilde{\sigma }}$ and if $\rho=\triv$, then $x_{\triv,\widetilde{\sigma }} = \half | 2n - (2n+1-a_{0}-1) + 1 | = \half (a_{0}+1) = x_{\triv,\varphi_{\widetilde{\sigma }}}$.
\end{proof}

\section{ Comparison of  combined Bernstein blocks of metaplectic groups and   unipotent Bernstein blocks of classical groups}
\label{sec:langl-deligne-param}

\subsection{Langlands--Deligne parameters and the algebras of intertwining operators for metaplectic groups}
\label{sec:case-metapl-groups}

\begin{defn}
  For each inertial class $\cO$ of an irreducible representation of $W_{\Fz}$, we fix a base point $\rho_{\cO}$ and call it a normed representation of $W_{\Fz}$.
  It is required to be unitary. In addition, $\rho_{\cO}$ is required to be self-dual, if  $\cO$ contains a self-dual element, and it is required to be  of  type $\lsup{L}{\tG}$, if $\cO$ contains a self-dual element of  type $\lsup{L}{\tG}$.

  A Langlands parameter $\varphi_{0} : W_{\Fz} \rightarrow \lsup{L}{\tG}$ is said to be normed, if $\iota\circ \varphi_{0}$ decomposes into irreducible representations of $W_{\Fz}$ which are the chosen base points in their respective  inertial class.

  Given a Langlands parameter or a Langlands--Deligne parameter $\varphi$, denote by $m(\rho;\varphi)$ the multiplicity of $\rho$ (up to equivalence) in the decomposition of $\iota\circ\varphi$ into irreducible representations of $W_{\Fz}$ and by $m_{\cO}(\rho;\varphi)$ the number of irreducible components of $\iota\circ\varphi$ that lie in the inertial orbit of $\rho $.

The inertial class of a normed Langlands parameter $\varphi _0$ is defined to be the set of all Langlands or Langlands--Deligne parameters $\varphi $ such that $m(\rho;\varphi_0)=m_{\cO}(\rho;\varphi)$ for every irreducible normed representation $\rho$ of $W_{\Fz}$.

  We will denote by the same symbol $\rho $ the supercuspidal representation of $\GL_{d_{\rho}}(\Fz)$ that corresponds to $\rho $ by the local Langlands correspondence for the general linear groups \cite{Harris-Taylor-MR1876802,Henniart-MR1738446}.
\end{defn}

When $\cO$ is the inertial class of an irreducible representation of $W_{\Fz}$, we have a bijection
\begin{align*}
  \cO &\rightarrow \CC^{\times}\\
  \rho &\mapsto f_{\rho} := \rho(\Fr^{t_{\rho}})\rho_{\cO}(\Fr^{t_{\rho}})^{-1}
\end{align*}
and if $\rho_{\cO}$ is self-dual, then $\rho\in\cO$ is self-dual if and only if $f_{\rho}=\pm 1$ (\cite[Proposition~1.4]{Heiermann-MR3719526}).

\begin{defn}
  Let $\cO$ be an inertial orbit of an irreducible representation of $W_{\Fz}$.
  Assume that the chosen base point  $\rho$ is self-dual.
  Write $\rho_{-}$ for the element $\rho'$ in $\cO$ with $f_{\rho'} = -1$, which is the other self-dual irreducible representation in $\cO$.
\end{defn}

Let $\varphi_{0}:W_{\Fz} \rightarrow\ \lsup{L}{\tG}$ be a normed Langlands parameter.

Denote by $\supp(\varphi_{0})$ the set of irreducible representations $\rho$ of $W_{\Fz}$, up to isomorphism, with $m(\rho;\varphi_{0})\neq 0$ and by $\supp'(\varphi_{0})$ the subset consisting of those representations which are self-dual.
We put an equivalence $\sim$ on $\supp(\varphi_{0})$ defined by $\rho \sim \rho^{\vee}$.

Let $\varphi $ be a Langlands--Deligne parameter for $\tG$ that lies in the inertial orbit of $\varphi_{0}$ and that corresponds to a supercuspidal representation $\widetilde\sigma $ of a Levi subgroup $\tM $ of $\tG $ such that $$\tM \isom \GL_{d_1}\times\cdots\times \GL_{d_r} \times \tH,\qquad \widetilde{\sigma} = \rho_1 \otimes\cdots\otimes\rho_r\otimes \widetilde{\tau},$$ with $\rho _1,\dots\rho _r$ normed and $\tH$ a metaplectic group.
Then, $\varphi $ factors through $\lsup{L}{\tM}\isom \GL_{d_1}(\CC )\times\cdots\times \GL_{d_r}(\CC ) \times\ \lsup{L}{\tH}$ by
\begin{equation*}
  \rho_1\otimes\cdots\otimes\rho_r\otimes\varphi_{\widetilde{\tau}}: W_{\Fz}\times \SL_{2}(\CC)\rightarrow \GL_{d_1}(\CC )\times\cdots\times \GL_{d_r}(\CC ) \times\ \lsup{L}{\tH},
\end{equation*}
where $\varphi_{\widetilde{\tau}}$ denotes the Langlands--Deligne parameter of $\widetilde{\tau }$.
As $\widetilde{\tau }$ is supercuspidal, one can write, by Theorem~\ref{thm:LD-param-and-char-comp-group},
\begin{align*}
  \varphi_{\widetilde{\tau}}=\bigoplus_{\rho\in \supp'(\varphi_{0})} \left( \bigoplus_{k=1}^{a_{\rho,+}} (\rho \otimes \spr(2k-\kappa_{\rho})) \oplus \bigoplus_{k=1}^{a_{\rho,-}} (\rho_{-} \otimes \spr(2k-\kappa_{\rho_{-}}))  \right) ,
\end{align*}
where $\kappa_{\rho}, \kappa_{\rho_{-}}\in\{0,1\}$ are determined by the type of $\rho, \rho_-$.
Let $S_{\varphi_{\widetilde{\tau}}}$ be the set $\{(a_{\rho,+},a_{\rho,-}) \in \ZZ_{\ge 0}\times\ZZ_{\ge 0} \spacedvert \rho\in \supp'(\varphi_{0}) \}$.
Write $m^{\GL}_{\rho}$ for the multiplicity of $\rho $ in the multiset $\{ \rho_1,\dots ,\rho_r \}$. Define
\begin{equation*}
  m_{+}(\rho;\varphi_{\widetilde{\tau}}) = \sum_{k=1}^{a_{\rho,+}} (2k-\kappa_{\rho}) =
  \begin{cases}
    a_{\rho,+}^{2},  &\quad\text{if $\kappa_{\rho}= 1$};\\
    a_{\rho,+}(a_{\rho,+} + 1), &\quad\text{if $\kappa_{\rho}= 0$}.
  \end{cases}
\end{equation*}
and similarly $m_{-}(\rho;\varphi_{\widetilde{\tau}}) = \sum_{k=1}^{a_{\rho,-}} (2k-\kappa_{\rho_{-}})$.
Then
\begin{equation*}
  m_{\cO}(\rho;\varphi_{\widetilde{\tau}})=m_{+}(\rho;{\varphi_{\widetilde{\tau}}})+m_{-}(\rho;{\varphi_{\widetilde{\tau}}}),
\end{equation*}
and it is equal to $m(\rho;\varphi_0)-2m^{\GL}_{\rho}$, as $\varphi $ lies in the inertial orbit of $\varphi _0$.
Thus for all $\rho\in \supp'(\varphi_{0})$, the pairs $(a_{\rho,+},a_{\rho,-})$ satisfy the identity

\begin{equation}\label{eq:size-parity-constraint-for-a+-}
   m(\rho;\varphi_0)-2m^{\GL}_{\rho} 
=    \begin{cases}
      a_{\rho,+}(a_{\rho,+} + 1)+a_{\rho,-}(a_{\rho,-} + 1), &\text{if $\rho $ and $\rho _-$ are not of type $\lsup{L}{\tG}$};\\
      a_{\rho,+}^{2}+a_{\rho,-}(a_{\rho,-} + 1),  &\text{if $\rho $ is of type $\lsup{L}{\tG}$, but  $\rho _-$ is not};\\
      a_{\rho,+}^{2}+a_{\rho,-}^{2}, &\text{if $\rho $ and $\rho _-$ are of type $\lsup{L}{\tG}$}.
    \end{cases}
\end{equation}

Now let
\begin{multline}\label{eq:defn-S-phi0-tuples}
  \cS(\varphi _0) = \{ (a_{\rho,+},a_{\rho,-})_{\rho\in\supp'(\varphi_{0})}\spacedvert \text{ for every $\rho\in\supp'(\varphi_{0})$}, \\ (a_{\rho,+},a_{\rho,-})\in \ZZ_{\ge 0}\times\ZZ_{\ge 0}\quad\text{satisfies \eqref{eq:size-parity-constraint-for-a+-} for some $m^{\GL}_{\rho }\in\ZZ_{\ge 0}$}  \},
\end{multline}
To every $S\in \cS(\varphi _0)$, which is a set of tuples indexed by $\rho\in\supp'(\varphi_{0})$, we can associate a representation of $W_{\Fz}\times\SL_{2}(\CC)$
\begin{align*}
  \bigoplus_{\rho\in \supp'(\varphi_{0})} \left( \bigoplus_{k=1}^{a_{\rho,+}} (\rho \otimes \spr(2k-\kappa_{\rho})) \oplus \bigoplus_{k=1}^{a_{\rho,-}} (\rho_{-} \otimes \spr(2k-\kappa_{\rho_{-}}))  \right)
\end{align*}
which factors through $\Sp_{2\ell}(\CC )$ for an appropriate $\ell$.
Then,  by Theorem~\ref{thm:LD-param-and-char-comp-group}, with $\tH^{S}=\tH^{S,\epsilon }=\Mp_{2\ell}$ (independent of $\epsilon$), it factors through a unique discrete Langlands--Deligne parameter
\begin{equation*}
  \varphi^S : W_{\Fz}\times\SL_{2}(\CC) \rightarrow\ \lsup{L}{{}\tH^{S}}.
\end{equation*}
Set as a short hand
\begin{equation*}
  \hat{S} = \hat{\scrS}_{\varphi^{S}}^{\alt},
\end{equation*}
where the right-hand side is defined in   Definition~\ref{defn:Salt}.
By Theorem~\ref{thm:LD-param-and-char-comp-group}, for every $\epsilon\in  \hat{S}$, the extended Langlands--Deligne parameter $(\varphi^S ,\epsilon)$ corresponds to a unique supercuspidal representation $\widetilde{\tau }^{S,\epsilon }$ of $\tH^{S}$.

Let
\begin{equation*}
  \tM_{S} \isom \prod_{\rho\in(\supp(\varphi_{0}) \setminus \supp'(\varphi_{0}))/\sim}\GL_{d_{\rho}}^{m(\rho;\varphi_{0})} \times \prod_{\rho\in\supp'(\varphi_{0}) } \GL_{d_{\rho}}^{m_{\rho, S }^{\GL}} \times \tH^{S} ,
\end{equation*}
where $m_{\rho,S}^{\GL} = (m(\rho;\varphi_{0}) - m_{\cO}(\rho,\varphi^{S}))/2$.
It is a Levi subgroup of $\tG=\Mp_{2n}$. Let $\varphi_S: W_{\Fz}\times\SL_{2}(\CC) \rightarrow \lsup{L}{\tM}_{S}$ be the Langlands--Deligne parameter such that
\begin{align*}
  \iota\circ\varphi_{S} = \bigoplus_{\rho\in\supp(\varphi_{0})} [m(\rho;\varphi_{0}) - m_{\cO}(\rho;\varphi^S)] \rho \oplus (\iota\circ \varphi^{S}).
\end{align*}
Thus, by the local Langlands correspondence, to $(\varphi_{S},\epsilon)$ corresponds the supercuspidal representation
\begin{equation*}
  \widetilde{\sigma}_{S,\epsilon} = \bigotimes_{\rho\in(\supp(\varphi_{0}) \setminus \supp'(\varphi_{0}))/\sim}\rho^{\otimes m(\rho;\varphi_{0})} \otimes \bigotimes_{\rho\in\supp'(\varphi_{0}) } \rho^{\otimes m_{\rho, S}^{\GL}} \otimes \widetilde{\tau }^{S,\epsilon }
\end{equation*}
of $\tM_{S}$.
Let $\cO_{S,\epsilon}$ denote the inertial class of $\widetilde{\sigma}_{S,\epsilon}$.
We have shown the following theorem (see also \cite[Theorem~1.6]{Heiermann-MR3719526}):

\begin{thm}
  The family $(\tM_{S}, \cO_{S,\epsilon})_{S,\epsilon}$ with $S$ running over $\cS(\varphi_0)$ defined in \eqref{eq:defn-S-phi0-tuples} and $\epsilon$ running over $\hat{S}$ exhausts the set of inertial orbits of supercuspidal representations $\widetilde{\sigma}'$ of standard Levi subgroups $\tM'$ of $\tG$ such that $(\varphi_{\widetilde{\sigma}'})|_{W_{\Fz}}$ lies in the inertial orbit of $\varphi_{0}$.
  Furthermore, $(\tM_{S}, \cO_{S,\epsilon}) = (\tM_{S'}, \cO_{S',\epsilon'})$ if and only if $S=S'$ and $\epsilon=\epsilon'$.
\end{thm}

Before  describing the Bernstein blocks $\Rep_{\tG}(\cO_{S,\epsilon})$, 
we explicate $\Sigma_{\cO,\mu}$ and $R(\cO)$ where $\cO$ is an inertial class of a supercuspidal representation of a standard Levi subgroup $\tM$ of $\tG$.

The following two propositions are  Proposition~\ref{prop:Sigma-O-mu-cover-decomposed-Levis} specialized to the case of  metaplectic groups, where $\mu^{\tM_{\alpha _{i,1}}}$ is never regular and the effect of $\bT_{(i,j)}$ in Proposition~\ref{prop:action-on-representations} is taking contragredient.

\begin{prop} \label{prop:Sigma-O-mu-explicit-metaplectic}
  \begin{enumerate}
  \item After conjugating $(\tM,\widetilde{\sigma})$ and twisting  $\widetilde{\sigma}$ by an element in $\Xnr(\tM)$, we can put $(\tM,\widetilde{\sigma})$ into the following form
    \begin{align*}
      \tM &\isom \GL_{d_{1}} \times \cdots \times \GL_{d_{1}}
          \times \GL_{d_{2}} \times \cdots \times \GL_{d_{2}}
          \times \cdots \times \GL_{d_{r}} \times \cdots \times \GL_{d_{r}}
          \times \tH_{d} \\
      \widetilde{\sigma} &\isom \sigma_{1} \otimes \cdots \otimes \sigma_{1}
               \otimes \sigma_{2} \otimes \cdots \otimes \sigma_{2}
               \otimes \cdots \otimes \sigma_{r} \otimes \cdots \otimes \sigma_{r}
               \otimes \widetilde{\tau},
    \end{align*}
    where  $H_{d}$ is a symplectic group of rank $d = \rank G  - \sum_{i=1}^{r} k_{i} d_{i}$ with $k_{i}$ being the number of $\GL_{d_{i}}$'s and the inertial classes of $\sigma_{i}$ are pairwise distinct and if $\sigma_{i}$ and its contragredient $\sigma_{i}^{\vee}$ are in the same inertial class, then $\sigma_{i} \isom \sigma_{i}^{\vee}$.
    
  \item Identify $A_{M}$ with $\GL_{1}^{k_{1}} \times \GL_{1}^{k_{2}} \times \cdots \times \GL_{1}^{k_{r}}$.
    Then the root system $\Sigma_{\cO,\mu}$ is the direct sum of $r$ components $\Sigma_{\cO,\mu, i}$ with $i=1,\ldots , r$ which are either irreducible or empty.
    The $\Sigma_{\cO,\mu, i}$'s are determined as follows.

    For $i=1,\ldots , r$ and $j = 1,\ldots, k_{i}$, set $\alpha_{i,j}$ to be the rational character of $A_{M}$ which sends
    \begin{align*} 
      ( x_{1,1}, x_{1,2}, \ldots , x_{1,k_{1}},
      x_{2,1}, x_{2,2}, \ldots , x_{2,k_{2}},
      \ldots , x_{r,1}, x_{r,2}, \ldots , x_{r,k_{r}})  \in A_{M}
    \end{align*}
    to $x_{i,j} x_{i,j+1}^{-1}$ if $j < k_{i}$
    and $x_{i,j}$ if $j=k_{i}$.
    Consider the $\mu$-function
    \begin{equation}\label{eq:mu-for-explicit-description-of-SigmaOmu-metaplectic}
      s \mapsto \mu(\sigma_{i} |{\textstyle \det_{d_{i}}} |_{\Fz}^{s} \otimes \widetilde{\tau})
    \end{equation}
    defined relative to $\GL_{d_{i}} \times \tH_{d}$ and $\tH_{d_{i} + d}$, where if $d=0$, $\widetilde{\tau}$ should be removed.
    Then we have the following.
    \begin{enumerate}
      \item If \eqref{eq:mu-for-explicit-description-of-SigmaOmu-metaplectic} has a pole in $\CC$, then a base for $\Sigma_{\cO,\mu,i}$ is
        \begin{align*}
          \{ \alpha_{i,1}, \alpha_{i,2}, \ldots, \alpha_{i,k_{i}} \} &\quad\text{if $d > 0$} \\
          \{ \alpha_{i,1}, \alpha_{i,2}, \ldots, 2\alpha_{i,k_{i}} \} &\quad\text{if $d = 0$} .
        \end{align*}
        This is a root system of type $B_{k_{i}}$ if $d>0$ and of type  $C_{k_{i}}$ if $d = 0$.

      \item   If \eqref{eq:mu-for-explicit-description-of-SigmaOmu-metaplectic} does not have a pole in $\CC$ and $\sigma_{i} \isom \sigma_{i}^{\vee}$, then a basis for $\Sigma_{\cO,\mu,i}$ is
        \begin{equation*}
          \{ \alpha_{i,1}, \alpha_{i,2}, \ldots, \alpha_{i,k_{i}-1},  \alpha_{i,k_{i}-1} + 2\alpha_{i,k_{i}}\}.
        \end{equation*}
        This is a root system of type $D_{k_{i}}$.

      \item  If \eqref{eq:mu-for-explicit-description-of-SigmaOmu-metaplectic} does not have a pole in $\CC$ and $\sigma_{i} \not\isom \sigma_{i}^{\vee}$, then a basis for $\Sigma_{\cO,\mu,i}$ is
        \begin{equation*}
          \{ \alpha_{i,1}, \alpha_{i,2}, \ldots, \alpha_{i,k_{i}-1}\}.
        \end{equation*}
        This is a root system of type $A_{k_{i}-1}$.
      \end{enumerate}
  \end{enumerate}
\end{prop}

\begin{prop}\label{prop:RO-structure}
  The group $R(\cO)$ is the direct product over $i\in I$ of $R(\cO)_{i}$  determined as follows.
  The index set $I$ consists of exactly those indices $i$ such that $\Sigma_{\cO,\mu,i}$ in Prop.~\ref{prop:Sigma-O-mu-explicit-metaplectic} is of type $D$ and
  each $R(\cO)_{i}$ is given by $\{ 1, s_{\alpha_{i},k_{i}}  \}$ where $s_{\alpha_{i},k_{i}}$ exchanges the roots $\alpha_{i,k_{i}-1}$ and $\alpha_{i,k_{i}-1} + 2\alpha_{i,k_{i}}$ defined in Prop.~\ref{prop:Sigma-O-mu-explicit-metaplectic}.
\end{prop}

Consider the Bernstein blocks $\Rep_{\tG}(\cO_{S,\epsilon})$. 
By Corollary~\ref{cor:bernstein-block-hecke-module} $\Rep_{\tG}(\cO_{S,\epsilon})$ is equivalent to the category of right modules over $\CC[R(\cO_{S,\epsilon})] \ltimes \cH(\Sigma_{\cO_{S,\epsilon}}, \{ q_{\alpha} \}, \{q_{i}\})$.
The root datum and the $R$-group are given in Propositions~\ref{prop:Sigma-O-mu-explicit-metaplectic} and \ref{prop:RO-structure} and the parameters $\{ q_{\alpha} \}, \{q_{i}\}$ 
 are determined by the invariants $a_{s_{\alpha}}, a_{s_{\alpha,-}}$ in the $\mu$-functions (see \eqref{eq:q-alpha}, \eqref{eq:q-i} and Proposition~\ref{prop:expression-mu-function}).
See also Theorem 5.2 of  \cite{Heiermann-MR2643577} and the remark thereafter.
The invariants $a_{s_{\alpha}}, a_{s_{\alpha,-}}$ are given either by the well-known points of reducibility for general linear groups or  by the points of reducibility of $\rho \rtimes \widetilde{\tau}$.
The latter are, in turn, determined by the Langlands--Deligne parameter of $\widetilde{\tau}$, as shown in Proposition~\ref{prop:Jord-point-red}.
We also note that  the cocycle $\eta$ can be made trivial by taking the normalisation in \cite[Section~2.5]{Heiermann-MR2827179}.
(See also Remark~\ref{rmk:2-cocycle}.)
This is what we do and we have suppressed $\eta$ from the notation of the group algebra $\CC[R(\cO_{S,\epsilon})]$.

Combining  these metaplectic  ingredients, we   have proved Theorem~\ref{thm:Bernstein-block-Hecke-in-terms-of-LD-param}, which is the analogue of \cite[Theorem~1.8]{Heiermann-MR3719526}.

\begin{thm}\label{thm:Bernstein-block-Hecke-in-terms-of-LD-param}
  Let $\varphi_{0}$ be a normed Langlands parameter for $\tG=\Mp_{2n}$.
  Let $S\in \cS(\varphi_{0})$ and $\epsilon \in \hat{S}$.
  Then $\Rep_{\tG}(\cO_{S,\epsilon})$ is equivalent to the category of right modules over $\bigotimes_{\rho\in\supp(\varphi_{0})/\sim} \cH_{\tG,\varphi_{0},S,\rho}$ (taken in the category of $\CC$-algebras)  where $\cH_{\tG,\varphi_{0},S,\rho}$ is an extended Hecke algebra defined as follows.
  \begin{enumerate}
  \item If $\rho$ is not self-dual, then $\cH_{\tG,\varphi_{0},S,\rho}$ is an affine Hecke algebra with root datum associated to $\GL_{m(\rho;\varphi_{0})}$ with equal parameters $q_{\Fz}^{t_{\rho}}$.
  \item Assume that $\rho$ is self-dual. If both $\rho$ and $\rho_{-}$ are of type $\lsup{L}{\tG}$ and neither of them is in $\supp(\varphi^{S})$, then $\cH_{\tG,\varphi_{0},S,\rho}$ is the semidirect product  of an affine Hecke algebra with root datum associated to the even special orthogonal group $\SO_{m(\rho;\varphi_{0})}$ with equal parameters $q_{\Fz}^{t_{\rho}}$ with the group algebra $\CC[\ZZ/2\ZZ]$ where the action of the nontrivial element of $\ZZ/2\ZZ$ is given by the outer automorphism on the root system.
  \item   Assume that $\rho$ is self-dual. If at least one of $\rho$ or $\rho_{-}$ is not of type $\lsup{L}{\tG}$ or at least one of them is in $\supp(\varphi^{S})$, then $\cH_{\tG,\varphi_{0},S,\rho}$ is an affine Hecke algebra with root datum associated to the odd special orthogonal group $\SO_{m(\rho;\varphi_{0}) - m_{\cO}(\rho;\varphi^{S}) +1}$ with  unequal parameters
\begin{equation*}
  q_{\Fz}^{t_{\rho}}, \ldots , q_{\Fz}^{t_{\rho}} , q_{\Fz}^{t_{\rho} (a_{\rho,+} + a_{\rho,-} + 1- \frac{\kappa_{\rho} +\kappa_{\rho_{-}}}{2})} ; q_{\Fz}^{| t_{\rho} (a_{\rho,+} - a_{\rho,-}  + \frac{\kappa_{\rho_{-}} -\kappa_{\rho}}{2})|}.
\end{equation*}
  \end{enumerate}
\end{thm}

\begin{rmk}
  The notation $\kappa_{\rho}$ (resp. $\kappa_{\rho_{-}}$) does not  conform to that in \cite{Heiermann-MR3719526}.
  If we denote the $\kappa_{\rho}$ in \cite{Heiermann-MR3719526} by $\kappa^{\text{Hei}}_{\rho}$, then $\kappa_{\rho} = 1 - \kappa^{\text{Hei}}_{\rho}$.
\end{rmk}

We combine the Bernstein blocks which come from the normed Langlands parameter $\varphi_{0}$.
Let
\begin{align*}
  \cR_{\Fz}^{\varphi_{0}}(\tG) = \bigoplus_{S\in\cS(\varphi_0)} \bigoplus_{\epsilon\in \hat{S}} \Rep_{\tG}(\cO_{S,\epsilon}).
\end{align*}
We will make use of  \cite[Appendix B]{Heiermann-MR3719526} which applies to categories of modules of finite presentation over a coherent $\CC$-algebra and, in particular, to the category of finitely generated modules over a noetherian $\CC$-algebra.
For a normed representation $\rho$ of $W_{\Fz}$, 
write $\varphi_{0,\rho}$ for the $\rho$-projection of $\varphi_{0}$ and $S_{\rho}$ the $\rho$-component of $S\in\cS(\varphi_{0})$.
More precisely, $\varphi_{0,\rho}$ is the unique Langlands parameter for a metaplectic group $\tG_{\rho}$ of an appropriate rank such that $\supp(\varphi_{0,\rho}) = \{ \rho \}$ and $m(\rho;\varphi_{0,\rho}) = m(\rho; \varphi_{0})$ and $S_{\rho}$ is the singleton corresponding to the tuple in $S$ indexed by $\rho$.
Notice that $\cH_{\tG,\varphi_{0},S,\rho} = \cH_{\tG_{\rho},\varphi_{0,\rho},S_{\rho},\rho}$.
Then, by \cite[B.4]{Heiermann-MR3719526}, we have the equivalence of categories
\begin{equation}\label{eq:bernstein-blocks-Hecke-general}
  \begin{aligned}
    \cR_{\Fz}^{\varphi_{0}}(\tG)_{f} &= \bigoplus_{S\in\cS(\varphi_{0})}\bigoplus_{\epsilon\in\hat{S}} \Rep_{\tG}(\cO_{S,\epsilon})_{f}\\
    &\isom \bigoplus_{S\in\cS(\varphi_{0})}\bigoplus_{\epsilon\in\hat{S}} \bigotimes_{\rho\in\supp(\varphi_{0})/\sim}(\text{Right-}\cH_{\tG,\varphi_{0},S,\rho}\text{-}\MOD)_{f}\\
    &\isom \bigotimes_{\rho\in\supp(\varphi_{0})/\sim} \bigoplus_{S_{\rho}\in\cS(\varphi_{0,\rho})}\bigoplus_{\epsilon_{\rho}\in\hat{S}_{\rho}} (\text{Right-}\cH_{\tG_{\rho},\varphi_{0,\rho},S_{\rho},\rho}\text{-}\MOD)_{f} .
  \end{aligned}
\end{equation}
We have used the subscript $f$ to indicate the full subcategory of finitely generated representations or the full subcategory of finitely generated modules, as this is needed for the last equivalence of categories.
See  the proof of \cite[Corollary~3.5]{Heiermann-MR3719526} for more details.

\subsection{Summary of the setup for classical groups}
\label{sec:summ-setup-class}

In the next section we will compare $\cR_{\Fz}^{\varphi_{0}}(\tG)$ with analogous objects for various classical groups. Let us briefly explain the setup of the latter from \cite{Heiermann-MR3719526}.
See \cite[1.5--1.7]{Heiermann-MR3719526} for the symplectic and orthogonal case, additionally with \cite[Appendix A]{Heiermann-MR3719526} addressing the non-connectedness of the full orthogonal groups,  and \cite[Appendix C.2--4]{Heiermann-MR3719526} for the unitary case.

Let $G$ denote  one of the split classical groups $\SO_{2n+1}$, $\Sp_{2n}$ and $\rO_{2n}$ or the quasi-split classical group $\rU_{n}$.
In the unitary group case, let $\Fs$ be the quadratic extension of $\Fz$ over which $\rU_{n}$ splits.
To unify notation, for other cases, set $\Fs=\Fz$.
Write $q$ for the cardinality of the residue field of $\Fs$.
In the split case, there is always a canonical embedding $\iota $ of the dual group into some $\GL_N(\CC )$.

Let $\varphi_{0}$ be a normed Langlands parameter for $G$, $\varphi_0:W_{\Fz}\rightarrow\ \lsup{L}{G}$.
In the case of the unitary group, we will replace $\varphi_{0}$ in the following by its restriction to the Weil group $W_{\Fs}$.

We can define analogously the set $\cS(\varphi_{0})$ and associate to $S\in \cS(\varphi_{0})$ a Langlands--Deligne parameter $\varphi_{S}$ in the inertial orbit of $\varphi_{0}$ that factors through some Levi subgroup $\lsup{L}{M_S}$ of $\lsup{L}{G}$.
The set $\hat{S}$ is the set of alternating characters of the component group associated to $\varphi^{S}$.
Remark that $\lsup{L}{M_S}\simeq \GL_{d_1}(\CC )\times\dots\times \GL_{d_r}(\CC )\times\ \lsup{L}{H^S}$, where $H^S$ is some classical group of the same type as $G$ and it  depends possibly on $\epsilon$ though this is not reflected in the notation here.

Set $d_{S}=+1$ if $\det(\varphi^{S}) = \det(\varphi_{0})$ and $d_{S}=-1$ otherwise.
For $\epsilon\in \hat{S}$, set $\epsilon_{Z} = \epsilon(-1)$.
We will use $G^{\epsilon_{Z}}_{d_{S}}$ to denote a form of $G$, which will be made more precise below for each group when $\supp(\varphi_{0}) = \{ \triv \}$.
Remark that some of the groups $G_\pm^\pm$ may remain undefined in some cases and we will omit $\epsilon_{Z}$ or $d_{S}$, if they do not apply.

The Langlands--Deligne parameter $\varphi_{S}$ decomposes as $(\rho_1, \ldots, \rho_r, \varphi^{S})$, where the $\rho_i$'s are irreducible representations of $W_\Fs$  and $\varphi^{S}$ a discrete Langlands parameter relative to $\lsup{L}{H^S}$.

For $\epsilon\in\hat{S}$, there is again a supercuspidal representation $\tau^{S,\epsilon}$ of $(H^S)^{\epsilon_{Z}}_{d_{S}}$ associated to $(\varphi^{S},\epsilon)$ and a supercuspidal representation $\sigma_{S,\epsilon}$ associated to $(\varphi_{S},\epsilon)$  which is a representation of $$(M_S)^{\epsilon_{Z}}_{d_{S}}\simeq \GL_{d_1}(\Fs)\times\dots\times \GL_{d_r}(\Fs)\times (H^S)^{\epsilon_{Z}}_{d_{S}}$$ and $\sigma_{S,\epsilon}=\rho_1\otimes\dots\otimes\rho _r\otimes\tau^{S,\epsilon}$.
Then this defines a Bernstein component attached to $\cO_{S,\epsilon}$ of $G^{\epsilon_{Z}}_{d_{S}}$ which may be different from $G$.

As $S$ and $\epsilon $ vary while keeping $\epsilon_{Z}$, $d_{S}$ fixed, again $\Rep_{G_{d_{S}}^{\epsilon_{Z}}}(\cO_{S,\epsilon})$ exhaust the Bernstein components of $G^{\epsilon_{Z}}_{d_{S}}$ with Langlands--Deligne parameter whose restriction to $W_{\Fs}$  is in the inertial class of $\varphi _0$.

Define
\begin{align*}
  \cR_{\Fs}^{\varphi_{0}}(G) = \bigoplus_{S\in\cS(\varphi_0)} \bigoplus_{\epsilon\in \hat{S}} \Rep_{G_{d_{S}}^{\epsilon_{Z}}}(\cO_{S,\epsilon})
\end{align*}
and define
\begin{align*}
  \cR_{\Fs}^{\varphi_{0}}(\GL_{n}) =   \Rep_{\GL_{n}}(\cO_{\varphi_{0}}).
\end{align*}
where $\cO_{\varphi_{0}}$ is the inertial orbit of the representation corresponding to $\varphi_{0}$ of an appropriate  Levi subgroup of $\GL_{n}$.

\null In the following we consider the special case where $\varphi_{0}$ is such that $\supp(\varphi_{0}) = \{ \triv \}$ and we write $\underline{\triv}$ for $\varphi_{0}$.
Each element $S$ in $\cS(\underline{\triv})$ is a singleton $(a_{+},a_{-}) \in \ZZ_{\ge 0}\times\ZZ_{\ge 0}$ such that \eqref{eq:size-parity-constraint-for-a+-} is satisfied for some $m=m_\triv^{\GL}$.
 Here we  write $a_{+}$ (resp. $a_{-}$) for $a_{\triv,+}$ (resp. $a_{\triv,-}$) and similarly $m_{+}$ and $m_{-}$ for $m_{+}(\triv, \varphi^S)$ and $m_{-}(\triv, \varphi^S)$, which are determined by $a_{+}$ and $a_{-}$.
We summarise the results of \cite[Theorems~1.9, together with C.5 and C.6 for the unitary group]{Heiermann-MR3719526} for the different cases of $G$ in this special situation.

\null $G=\SO_{2n+1}$. Then $\lsup{L}{G}=\Sp_{2n}(\CC)$. Neither $\triv$ or $\triv_{-}$ is of type $\lsup{L}{G}$.
Denote by $\cH_{\SO_{2n+1}}(a_{+},a_{-})$ the Hecke algebra with root datum associated to $\SO_{2n+1-m_{+}-m_{-}}$ and unequal parameters
\begin{align*}
  q, \ldots , q, q^{a_{+}+a_{-}+1}; q^{|a_{+}-a_{-}|}.
\end{align*}
Then $\Rep_{\SO_{2n+1}^{\epsilon_{Z}}}(\cO_{S,\epsilon})$ is equivalent to the category of right modules over $\cH_{\SO_{2n+1}}(a_{+},a_{-})$.
Here $\SO_{2n+1}^{+}$ is the split $\SO_{2n+1}$ and $\SO_{2n+1}^{-}$ is its pure inner form.

\null $G=\rO_{2n}$. Then $\lsup{L}{G}=\rO_{2n}(\CC)$. Both $\triv$ and $\triv_{-}$ are of type $\lsup{L}{G}$.
Assume that at least one of $a_{+}$ and $a_{-}$ is strictly greater than $0$.
Denote by $\cH_{\rO_{2n}}(a_{+},a_{-})$ the Hecke algebra with root datum associated to $\SO_{2n+1-m_{+}-m_{-}}$ and unequal parameters
\begin{align*}
  q, \ldots , q, q^{a_{+}+a_{-}}; q^{|a_{+}-a_{-}|}.
\end{align*}
Denote by $\cH_{\rO_{2n}}(0,0)$ the extended Hecke algebra which is the semidirect product of $\CC[\ZZ/2\ZZ]$ with the Hecke algebra with root datum associated to $\SO_{2n}$ and equal parameters
\begin{align*}
  q, \ldots , q.
\end{align*}
Then $\Rep_{\rO_{2n,d_{S}}^{\epsilon_{Z}}}(\cO_{S,\epsilon})$ is equivalent to the category of right modules over $\cH_{\rO_{2n}}(a_{+},a_{-})$.
Here $\rO_{2n,+}^{+}$ is the split $\rO_{2n}$ and $\rO_{2n,+}^{-}$ is its pure inner form; $\rO_{2n,-}^{+}$ is the quasi-split $\rO_{2n}$ associated to the unramified quadratic extension of $\Fs$ and $\rO_{2n,-}^{-}$ is its pure inner form.

 $G=\Sp_{2n}$. Then we take here $\lsup{L}{G}=\rO_{2n+1}(\CC)$. Both $\triv$ and $\triv_{-}$ are of type $\lsup{L}{G}$.
 At least one of $a_{+}$ and $a_{-}$ is strictly greater than $0$.
Denote by $\cH_{\Sp_{2n}}(a_{+},a_{-})$ the Hecke algebra with root datum associated to $\SO_{2n+1-m_{+}-m_{-}}$ and unequal parameters
\begin{align*}
  q, \ldots , q, q^{a_{+}+a_{-}}; q^{|a_{+}-a_{-}|}.
\end{align*}
Then $\Rep_{\Sp_{2n,d_{S}}}(\cO_{S,\epsilon})$ is equivalent to the category of right modules over $\cH_{\Sp_{2n}}(a_{+},a_{-})$.
Here $\Sp_{2n,+}^{+} = \Sp_{2n,+}^{-} = \Sp_{2n,-}^{+} = \Sp_{2n,-}^{-} = \Sp_{2n}$.

$G=\rU_{n}$ which is the quasi-split  unitary group associated to the unramified quadratic extension $\Fs/\Fz$ and a Hermitian space with discriminant $1$.
Then $\lsup{L}{G}=\GL_{n}(\CC)$.
When $n$ is odd, then $\triv$ is of type $\lsup{L}{G}$ and $\triv_{-}$ is not of type $\lsup{L}{G}$.
When $n$ is even, then $\triv$ is not of type $\lsup{L}{G}$ and $\triv_{-}$ is of type $\lsup{L}{G}$.
For the unitary case, we relax the condition that the base point in an inertial orbit must be of type $\lsup{L}{G}$ if such an element exists.
Denote by $\cH_{\rU_{n}}(a_{+},a_{-})$ the Hecke algebra with root datum associated to $\SO_{n-m_{+}-m_{-}+1}$ and unequal parameters
\begin{align*}
  q, \ldots , q, q^{a_{+}+a_{-}+\half}; q^{|a_{+}-a_{-} + \half(-1)^{n}|}.
\end{align*}
Then $\Rep_{\rU_{n}^{\epsilon_{Z}}}(\cO_{S,\epsilon})$ is equivalent to the category of right modules over $\cH_{\rU_{n}}(a_{+},a_{-})$.
Here $\rU_{n}^{+} = \rU_{n}$  and $\rU_{n}^{-}$ is its pure inner form.

\null $G=\GL_{n}$. Then $\cR_{\Fs}^{\underline{\triv}}(\GL_{n})$ is equivalent to the category of right modules over the Hecke algebra with root datum associated to $\GL_{n}$ and equal parameters $q$.

\subsection{Comparison of certain combined Bernstein blocks}
\label{sec:comp-bernst-blocks}

Comparing the extended affine Hecke algebras coming from classical groups using the Langlands parameter $\underline{\triv}$ and those from the metaplectic groups using $\varphi_{0}$ which is supported on $\rho$ and $\rho^{\vee}$, we get the following result.

Let $\Fs_{t_{\rho}}$ denote the unramified extension of $\Fs$ of degree $t_{\rho}$ so that the cardinality of the residue field of $F_{t_{\rho}}$ is $q^{t_{\rho}}$.

\begin{thm}\label{thm:bernstein-block-comparison-singleton-supp}
  Let $\varphi_{0}$ be a normed Langlands parameter for $\tG=\Mp_{2n}$.
  Assume that $\supp(\varphi_{0}) = \{ \rho, \rho^{\vee}  \}$.
  Set $m= m(\rho;\varphi_{0})$.
  Then the following holds.
  \begin{enumerate}
  \item If $\rho$ is not self-dual, then the category $\cR_{\Fs}^{\varphi_{0}}(\tG)$ is equivalent to $\cR_{\Fs_{t_{\rho}}}^{\underline{\triv}}(\GL_{m})$.
  \item Assume that $\rho$ is  self-dual and not of type $\lsup{L}{\tG}$.
    Then by our choice of base points, $\rho_{-}$ is not of type $\lsup{L}{\tG}$ either.
    Thus the category $\cR_{\Fs}^{\varphi_{0}}(\tG)$ is equivalent to $\cR_{\Fs_{t_{\rho}}}^{\underline{\triv}}(\SO_{m+1})$.
  \item If $\rho$ is  self-dual and both $\rho$ and $\rho_{-}$  are of type $\lsup{L}{\tG}$, then the category $\cR_{\Fs}^{\varphi_{0}}(\tG)$ is equivalent to $\cR_{\Fs_{t_{\rho}}}^{\underline{\triv}}(\rO_{m})$ if $m$ is even and to $\cR_{\Fs_{t_{\rho}}}^{\underline{\triv}}(\Sp_{m-1})$ if $m$ is odd.
  \item If $\rho$ is  self-dual,  $\rho$ is of type $\lsup{L}{\tG}$ and $\rho_{-}$ is not  of type $\lsup{L}{\tG}$, then the category $\cR_{\Fs}^{\varphi_{0}}(\tG)$ is equivalent to $\cR_{\Fs_{t_{\rho}}}^{\underline{\triv}}(\rU_{m})$.
  \end{enumerate}
\end{thm}

\begin{proof}
  Assume that $\rho$ is not self-dual.
  Then $\cR_{\Fs}^{\varphi_{0}}(\tG) = \Rep_{\tG}(\cO)$ where $\cO$ is the inertial orbit of $\rho\otimes \cdots \otimes \rho$.
 By Theorem~\ref{thm:Bernstein-block-Hecke-in-terms-of-LD-param}, it is equivalent to the category of right modules over an affine Hecke algebra with root datum associated to $\GL_{m}$ with equal parameters $q^{t_{\rho}}$ and this is exactly equivalent to $\cR_{\Fs_{t_{\rho}}}^{\underline{\triv}}(\GL_{m})$.

  Next assume that $\rho$ is self-dual and then $\supp(\varphi_{0}) = \{ \rho  \}$.
  Then by Theorem~\ref{thm:Bernstein-block-Hecke-in-terms-of-LD-param}, $\cR_{\Fs}^{\varphi_{0}}(\tG)$ is equivalent to the category of right modules over 
  \begin{align}\label{eq:prod-Hecke}
    \bigoplus_{S\in\cS(\varphi_{0})}\bigoplus_{\epsilon\in\hat{S}} \cH_{\tG,\varphi_{0}, S, \rho}.
  \end{align}
  The Hecke algebras are defined in the statements of Theorem~\ref{thm:Bernstein-block-Hecke-in-terms-of-LD-param}.

  We compare the Hecke algebras with those in Section~\ref{sec:summ-setup-class} for the classical groups.
  Up to changing all  $q$ which occur in the parameters for the Hecke algebras on the right-hand side to $q^{t_{\rho}}$, we have
\small %
  \begin{align}\label{eq:compare-Hecke}
    \cH_{\tG,\varphi_{0},S,\rho} =
    \begin{cases}
      \cH_{\SO_{m+1}}(a_{\rho,+},a_{\rho,-}), &\hspace{-1em}\text{ if $\rho$ and $\rho_{-}$ are not of type $\lsup{L}{\tG}$}; \\
      \cH_{\rU_{m}}(a_{\rho,+},a_{\rho,-}), &\hspace{-1em}\text{ if $\rho$ is of type $\lsup{L}{\tG}$,  $\rho_{-}$ is not of type $\lsup{L}{\tG}$ and $m$ is odd};\\
      \cH_{\rU_{m}}(a_{\rho,-},a_{\rho,+}), &\hspace{-1em}\text{ if $\rho$ is of type $\lsup{L}{\tG}$,  $\rho_{-}$ is not of type $\lsup{L}{\tG}$ and $m$ is even}; \\
      \cH_{\rO_{m}}(a_{\rho,+},a_{\rho,-}), &\hspace{-1em}\text{ if $\rho$ and $\rho_{-}$ are of type $\lsup{L}{\tG}$  and $m$ is even }; \\
      \cH_{\Sp_{m-1}}(a_{\rho,+},a_{\rho,-}), &\hspace{-1em}\text{ if $\rho$ and $\rho_{-}$ are of type $\lsup{L}{\tG}$  and $m$ is odd}.
    \end{cases}
  \end{align}
  \normalsize
  Note that in the first case above, $m$ is forced to be even.
  
  Write $H_{\rho}(m)=\rU_{m},\SO_{m+1},\Sp_{m-1},\rO_{m}$ determined by the types of $\rho$ and $\rho_{-}$ according to the group in the subscript on the right-hand side of \eqref{eq:compare-Hecke}.
Then  the category $\cR_{\Fs}^{\underline{\triv}}(H_{\rho}(m))$ is equivalent to the category of right modules over
  \begin{align}\label{eq:prod-Hecke-for-triv}
    \bigoplus_{S = \{ (a_{+},a_{-})  \}\in\cS(\underline{\triv})} \bigoplus_{\epsilon\in\hat{S}} \cH_{H_{\rho}(m)}(a_{+},a_{-}) .
  \end{align}

  Since the  direct sums in \eqref{eq:prod-Hecke} and \eqref{eq:prod-Hecke-for-triv} match,    we get the desired equivalence of categories between $\cR_{\Fs}^{\varphi_{0}}(\tG)$ and $\cR_{\Fs_{t_{\rho}}}^{\underline{\triv}}(H_{\rho}(m))$.
\end{proof}

Using \eqref{eq:bernstein-blocks-Hecke-general}, we can decompose a normed Langlands parameter $\varphi_{0}$ according to its $\rho$-projections to reduce to the situation of Theorem~\ref{thm:bernstein-block-comparison-singleton-supp}.
Then the proof of \cite[Corollary~3.5]{Heiermann-MR3719526} carries over to give us the following.
The statements on square-integrability and temperedness come from Theorem~\ref{thm:preserv-sq-int-temper}.

\begin{cor}
  Let $\varphi_{0}$ be a normed Langlands parameter for $\tG=\Mp_{2n}$.
  The category $\cR_{\Fs}^{\varphi_{0}}(\tG)_{f}$ is equivalent to
  \begin{align*}
    \bigotimes_{\rho\in\supp(\varphi_{0})/\sim} \cR_{\Fs_{t_{\rho}}}^{\underline{\triv}}(H_{\rho}(m(\rho;\varphi_{0})))_{f}
  \end{align*}
  where
  \begin{equation*}
    H_{\rho}(m) =
        \begin{cases}
           \GL_{m}, &\text{ if $\rho$ is not self-dual}; \\
          \SO_{m+1}, &\text{ if $\rho$ and $\rho_{-}$ are not of type $\lsup{L}{\tG}$ ($m$ is forced to be even)}; \\
          \rU_{m}, &\text{ if $\rho$ and $\rho_{-}$ are not of the same type};\\
          \rO_{m}, &\text{ if $\rho$ and $\rho_{-}$ are of type $\lsup{L}{\tG}$  and $m$ is even }; \\
          \Sp_{m-1}, &\text{ if $\rho$ and $\rho_{-}$ are of type $\lsup{L}{\tG}$  and $m$ is odd}.
    \end{cases}
  \end{equation*}
  This equivalence of categorie preserves discrete series representations, tempered representations, parabolic induction and the Jacquet functor.
\end{cor}

We note a relation between the categories of representations for $\Mp_{2n}$ and $\SO_{2n+1}$ as they share the same $L$-group $\Sp_{2n}(\CC)$.
As we have noted that
\begin{align*}
  \cR_{\Fs}^{\varphi_{0}}(\Mp_{2n})_{f} &= \bigoplus_{S\in\cS(\varphi_{0})}\bigoplus_{\epsilon\in\hat{S}} \Rep_{\Mp_{2n}}(\cO_{S,\epsilon})_{f}\\
  &\isom \bigoplus_{S\in\cS(\varphi_{0})}\bigoplus_{\epsilon\in\hat{S}} \bigotimes_{\rho\in\supp(\varphi_{0})/\sim}(\text{Right-}\cH_{\Mp_{2n},\varphi_{0},S,\rho}\text{-}\MOD)_{f}.
\end{align*}
We also have the analogous result for $\SO_{2n+1}$ given in \cite[Corollary~1.9]{Heiermann-MR3719526} which is valid for classical groups.
Let $\bullet$ be $+$ or $-$.
We have
\begin{align*}
  \cR_{\Fs}^{\varphi_{0}}(\SO_{2n+1}^{\bullet})_{f} &:= \bigoplus_{S\in\cS(\varphi_{0})}\bigoplus_{\substack{\epsilon\in\hat{S}\\ \epsilon_{Z}=\bullet}} \Rep_{\SO_{2n+1}^{\bullet}}(\cO_{S,\epsilon})_{f} \\
  &\isom \bigoplus_{S\in\cS(\varphi_{0})}\bigoplus_{\substack{\epsilon\in\hat{S}\\ \epsilon_{Z}=\bullet}} \bigotimes_{\rho\in\supp(\varphi_{0})/\sim}(\text{Right-}\cH_{\SO_{2n+1},\varphi_{0},S,\rho}\text{-}\MOD)_{f}
\end{align*}
and
\begin{align*}
  \cR_{\Fs}^{\varphi_{0}}(\SO_{2n+1})_{f} &= \cR_{\Fs}^{\varphi_{0}}(\SO_{2n+1}^{+})_{f} \oplus \cR_{\Fs}^{\varphi_{0}}(\SO_{2n+1}^{-})_{f}\\
  & \isom \bigoplus_{S\in\cS(\varphi_{0})}\bigoplus_{\epsilon\in\hat{S}} \bigotimes_{\rho\in\supp(\varphi_{0})/\sim}(\text{Right-}\cH_{\SO_{2n+1},\varphi_{0},S,\rho}\text{-}\MOD)_{f}.
\end{align*}
The inertial orbits $\cO_{S,\epsilon}$ are for Levi subgroups of different groups which should be clear from the context.
We note that the affine Hecke algebras depend only on $\varphi_{0},S,\rho$.
In particular, they depend on the $L$-group $\Sp_{2n}$, but not on whether it is the $L$-group for $\Mp_{2n}$ or $\SO_{2n+1}$.
Thus $\cH_{\Mp_{2n},\varphi_{0},S,\rho} = \cH_{\SO_{2n+1},\varphi_{0},S,\rho}$.
We get an equivalence of categories $\cR_{\Fs}^{\varphi_{0}}(\Mp_{2n})_{f} \isom \cR_{\Fs}^{\varphi_{0}}(\SO_{2n+1})_{f}$.
We can, in addition, separate $\cR_{\Fs}^{\varphi_{0}}(\Mp_{2n})_{f} $ into two parts according to the value $\epsilon_{Z}$.

This can be summarised in the following theorem.
Write $\cR_{\Fs}(\Mp_{2n})_{f}$ for $\oplus_{\varphi_{0}} \cR_{\Fs}^{\varphi_{0}}(\Mp_{2n})_{f} $ and
$\cR_{\Fs}(\SO_{2n+1}^{\bullet})_{f}$ for $\oplus_{\varphi_{0}} \cR_{\Fs}^{\varphi_{0}}(\SO_{2n+1}^{\bullet})_{f} $ where $\varphi_{0}$ runs over all normed Langlands parameters mapping to $\Sp_{2n}(\CC)$.
The statements on square-integrability and temperedness come from Theorem~\ref{thm:preserv-sq-int-temper}.

\begin{thm}
  The category $\cR_{\Fs}(\Mp_{2n})_{f}$ of smooth finitely generated genuine complex representations of $\Mp_{2n}$ is naturally equivalent to the direct sum of the categories of smooth finitely generated complex representations of $\SO_{2n+1}^{+}$ and $\SO_{2n+1}^-$, $$\cR_{\Fs}(\Mp_{2n})_{f}\sim \cR_{\Fs}(\SO_{2n+1}^{+})_{f}\oplus \cR_{\Fs}(\SO_{2n+1}^-)_{f}.$$

 The representations of $\Mp_{2n}$ with supercuspidal support given by an extended Langlands--Deligne parameter $(\varphi,\epsilon)$ with $\epsilon_Z=1$ correspond to $\cR_{\Fs}(\SO_{2n+1}^{+})_{f}$ and those with supercuspidal support given by an extended Langlands--Deligne parameter $(\varphi,\epsilon)$ with $\epsilon_Z=-1$ correspond to $\cR_{\Fs}(\SO_{2n+1}^-)_{f}$.

  This equivalence of categories preserves discrete series representations, tempered representations, parabolic induction and the Jacquet functor.
\end{thm}

\begin{eg}
  The Weil representation $\omega$ of $\Mp_{2n}$ has a decomposition into even  and odd constituents $\omega = \omega^{+} \oplus \omega^{-}$.
 Consider the  Bernstein component of $\Mp_{2n}$ that contains $\omega^{+}$ (resp. $\omega^{-}$) and that of $\SO_{2n+1}^{\pm}$ that contains the trivial representation.

  The Bernstein component of $\SO_{2n+1}^{+}$ that contains the trivial representation is associated to the inertial orbit of the trivial representation of the maximal split torus.
  We denote this Bernstein block by $\Rep_{\SO_{2n+1}^{+}}^{\triv}$.
  The associated Langlands--Deligne parameter is   $\underbrace{\triv \oplus \cdots \oplus \triv}_{{2n\text{-times}}}$.
  By the case for odd special orthogonal groups in Section~\ref{sec:summ-setup-class}, $\Rep_{\SO_{2n+1}^{+}}^{\triv}$ is equivalent to the category of right modules over the affine Hecke algebra associated to $\SO_{2n+1}$ with unequal parameters
  \begin{align*}
    \underbrace{q, \ldots , q}_{n\text{-times}} ; q^{0}.
  \end{align*}

  The Bernstein component of $\SO_{2n+1}^{-}$ that contains the trivial representation is associated to the inertial orbit of the trivial representation of the minimal Levi subgroup.
  We denote this Bernstein block by $\Rep_{\SO_{2n+1}^{-}}^{\triv}$.
  The associated Langlands--Deligne parameter is   $\underbrace{\triv \oplus \cdots \oplus \triv}_{(n-1)\text{-times}} \oplus (\triv \otimes\spr(2)) \oplus \underbrace{\triv \oplus \cdots \oplus \triv}_{(n-1)\text{-times}} $.
  By the case for odd special orthogonal groups in Section~\ref{sec:summ-setup-class}, $\Rep_{\SO_{2n+1}^{-}}^{\triv}$ is equivalent to the category of right modules over the affine Hecke algebra associated to $\SO_{2n-1}$ with unequal parameters
  \begin{align*}
    \underbrace{q, \ldots , q}_{(n-2)\text{-times}}, q^{2} ; q.
  \end{align*}
  
The Bernstein component of $\Mp_{2n}$ that contains $\omega^{+}$ is associated to the inertial orbit of the trivial representation of the maximal split torus, by Lemma~3.2 of \cite{Takeda-Wood-MR3729496}.
  We denote this Bernstein block by $\Rep_{\Mp_{2n}}^{\omega^{+}}$.
  The associated Langlands--Deligne parameter is   $\underbrace{\triv \oplus \cdots \oplus \triv}_{{2n\text{-times}}}$.
  By Theorem~\ref{thm:Bernstein-block-Hecke-in-terms-of-LD-param}, $\Rep_{\Mp_{2n}}^{\omega^{+}}$ is equivalent to the category of right modules over the affine Hecke algebra associated to $\SO_{2n+1}$ with unequal parameters
  \begin{align*}
    \underbrace{q, \ldots , q}_{n\text{-times}} ; q^{0}.
  \end{align*}

  The Bernstein component of $\Mp_{2n}$ that contains $\omega^{-}$ is associated to the inertial orbit of $\triv \otimes \omega_{1}^{-}$ of the Levi subgroup which is isomorphic to $\GL_{1} \times \cdots \times \GL_{1} \times \Mp_{2}$, by Lemma~3.5 of \cite{Takeda-Wood-MR3729496}, where $\triv$ denotes the trivial representation of $\GL_{1} \times \cdots \times \GL_{1}$ and $\omega_{1}^{-}$ is the odd component of the Weil representation $\omega_{1}$ of  $\Mp_{2}$, which is supercuspidal.
  We denote this Bernstein block by $\Rep_{\Mp_{2n}}^{\omega^{-}}$.
  Since $\omega_{1}^{-}$ is the first occurrence representation in Howe correspondence with the sign representation of $\rO_{1}$, M\oe glin's construction in Section~\ref{sec:some-results-moe} shows that the  Langlands--Deligne parameter of $\omega_{1}^{-}$ is $(\triv \otimes\spr(2))$.
    Thus the associated Langlands--Deligne parameter for $\triv \otimes \omega_{1}^{-}$ is   $\underbrace{\triv \oplus \cdots \oplus \triv}_{(n-1)\text{-times}} \oplus (\triv \otimes\spr(2)) \oplus \underbrace{\triv \oplus \cdots \oplus \triv}_{(n-1)\text{-times}} $.
  By Theorem~\ref{thm:Bernstein-block-Hecke-in-terms-of-LD-param}, $\Rep_{\Mp_{2n}}^{\omega^{-}}$ is equivalent to the category of right modules over the affine Hecke algebra associated to $\SO_{2n-1}$ with unequal parameters
  \begin{align*}
    \underbrace{q, \ldots , q}_{(n-2)\text{-times}}, q^{2} ; q.
  \end{align*}

  Noticing that the Hecke algebras match, we get equivalences of categories
  \begin{align*}
    \Rep_{\Mp_{2n}}^{\omega^{+}} \isom \Rep_{\SO_{2n+1}^{+}}^{\triv}   \quad\text{and}\quad \Rep_{\Mp_{2n}}^{\omega^{-}} \isom \Rep_{\SO_{2n+1}^{-}}^{\triv}
  \end{align*}
  which preserve square-integrability and temperedness.

  This recovers the equivalences of categories in   \cite{Gan-Savin-MR2982417} (under the restriction that the residue characteristic is not $2$) which is generalised by \cite{Takeda-Wood-MR3729496}  to also include the case of residue characteristic $2$.
  We note that our affine Hecke algebras, which are in the Bernstein--Lusztig presentation correspond exactly to those in \cite{Gan-Savin-MR2982417} and \cite{Takeda-Wood-MR3729496} which are in the Iwahori presentation.

\end{eg}  

\section{Some remarks in the case of non decomposed Levi subgroups}
\label{sec:non-decomposed}

When the Levi subgroups of $\tG$ do not decompose, the map 
$$\widetilde{\GL_{d_1,1}(\Fs)} \times\cdots\times \widetilde{\GL_{d_r,r}(\Fs)}\times \tH_d\rightarrow \widetilde{M}$$
in Definition~\ref{defn:prod-of-covers} won't be a homomorphism anymore, only a surjective map. The proof for multiplicity $1$ for Proposition~\ref{prop:M1-restriction-mult-one} is not valid anymore and this property may well fail, which makes a direct application of the results in \cite{Heiermann-MR2827179} impossible.

An additional complication which makes the $R$-group and the overall description more involved is the following: let $\alpha $ be a root in $\Sigma(A_M)$ such that $s_{\alpha }$ leaves $M$ invariant. Denote by $M^{\alpha }$ the product of factors of $M$, on which $s_{\alpha }$ acts by the identity, $\tM^{\alpha }$ its preimage and $\widetilde{\sigma }^{\alpha }$ the restriction of $\widetilde{\sigma }$ to $\tM^{\alpha }$. Then, there is a non-genuine character $\chi^{\alpha }$ (trivial on the center) of  $\tM$ with values in $\bmu _\bm$ and trivial on the preimage $\tM_{(\alpha )}$ of the other factors of $M$ , such that $s_{\alpha }$ acts on $\tM^{\alpha }$ by $\chi^{\alpha }$ and $s_{\alpha }\widetilde{\sigma }^{\alpha}$ is isomorphic to $\widetilde{\sigma }^{\alpha}{\chi^{\alpha }}$. This character $\chi^{\alpha }$ may be unramified or not.  

If $s_{\alpha }\widetilde{\sigma }_{\vert \tM_{(\alpha )}}\isom \widetilde{\sigma }_{\vert \tM_{(\alpha )}}\otimes\chi_{(\alpha )}$ for some character $\chi_{(\alpha )}$ of  $\tM_{(\alpha )}$
then $\chi_{(\alpha )}$ can be extended to $\tM$ trivial on $\tM^{\alpha }$. If this isomorphism commutes with $\widetilde{\sigma }_{\tM^{\alpha }}$, then $s_{\alpha }\widetilde{\sigma }\isom\widetilde{\sigma } \otimes  \chi_{(\alpha )}\chi^{\alpha }$ and conversely.

It is clear that this makes things more complicate. In particular, if $s_{\alpha}=:r\in  R(\cO)$, we will not always have $r\widetilde{\sigma}\isom\widetilde{\sigma }$, but rather $r\widetilde{\sigma}\isom\widetilde{\sigma }\otimes\chi_r$ with $\chi_r$ some unramified character.
This requires a modification of the definition of $\rho _r$ and the composition properties of $J_r$ become more complicated, including the commuting rules with the affine Hecke algebra part. It will not be anymore a semi-direct product with the twisted group algebra of $R(\cO)$, but something more difficult to handle (compare with \cite[p. 468]{Solleveld-MR4432237}).

Although a comprehensive description in this more general situation does not seem impossible, it won't be much readable and it is probably better to study specific examples along these lines.

We will give in a future version more information on the non-decomposed case.

\section{Annex by Volker Heiermann: Errata from  \cite{Heiermann-MR3719526}}
\label{sec:Errata}

1. $G^-$ is defined even if $G$ is the symplectic group (see introduction of \cite{Heiermann-MR3719526}). It is just an identical copy of $G$ (compare  \cite[section 1]{Heiermann-MR3719526}). The reason is that we will take characters of centralizers in $^LG$ which is for us the full odd orthogonal group and those with value $-1$ will correspond to the copy  $G^-$.

2. If $G$ is the symplectic group, the part of a Langlands parameter $W_F\rightarrow\ ^LG$ with value in the direct factor $\pm 1$ of $^LG$ will be always supposed to be unramified, i.e. it is either trivial or equal to the quadratic unramified character.

3. If $G$ is the symplectic group, one has always $t_o\ge1$ in \cite[Theorem 1.1 (iii)]{Heiermann-MR3719526}. 

4. If $G$ is the symplectic group, $G_-^-$ has to be defined as an identical copy of $G$.

5. In \cite[Section 1.7]{Heiermann-MR3719526} (after the title "Langlands correspondence"), one has to take for $\Xi $ the irreducible representations of $C_{^LG}(Im(\phi))/C_{^LG}(Im(\phi))^0$ to get the whole category $\cR^{\varphi_0}_{\Fs}(\underline{G})$.

6. In \cite[Theorem 3.2 (i)]{Heiermann-MR3719526}, the set of triples $(s,u,\Xi )$ is understood to be taken relative to $\SO_{2d+1}(\Bbb C)$. That means that $\Xi $ is an irreducible representation of the connected centralizer of $s$ and $u$ in $\SO_{2d+1}(\Bbb C)$. When $d_-=0$, a single category of irreducible $\cH(d_+,0)$-modules appears, not two as indicated.

7. The remarks relative to the $R$-group of the special even orthogonal group at the end of \cite[Appendix A]{Heiermann-MR3719526} were corrected in \cite{Ciubotaru-Heiermann-MR3960011}, but need further correction. The correct version is contained in Proposition~\ref{prop:Sigma-O-mu-cover-decomposed-Levis}. Specializing this proposition to the case of a special even orthogonal (reductive) group, one gets the correct description of the $R$-group.
\vspace{\baselineskip}

Volker Heiermann thanks Chenyan Wu for  some of the errata.

\bibliographystyle{alpha}

\bibliography{MyBib-20220401}
\end{document}